\documentclass[a4paper,11pt]{article}

\usepackage{mathrsfs}
\usepackage{amsmath,amscd}
\usepackage{amssymb, amsmath}
\usepackage{graphicx}
\usepackage{color}

\allowdisplaybreaks

\newcommand{\newcom}{\newcommand}

\newcom{\cA}{{\mathcal A}}
\newcom{\cB}{{\mathcal B}}
\newcom{\cC}{{\mathcal C}}
\newcom{\cD}{{\mathcal D}}
\newcom{\cE}{{\mathcal E}}
\newcom{\cF}{{\mathcal F}}

\newcom{\cJ}{{\mathcal J}}

\newcom{\cL}{{\mathcal L}}
\newcom{\cM}{{\mathcal M}}
\newcom{\cP}{{\mathcal P}}
\newcom{\cS}{{\mathcal S}}
\newcom{\cQ}{{\mathcal Q}}
\newcom{\cT}{{\mathcal T}}
\newcom{\cY}{{\mathcal Y}}
\newcom{\cZ}{{\mathcal Z}}
\newcom{\R}{\mathbb R}
\newcom{\T}{\mathbb T}
\newcom{\N}{\mathbb N}
\newcom{\Z}{\mathbb Z}
\newcom{\C}{\mathbb C}
\newcom{\E}{\mathbb E}

\newcom{\e}{\epsilon}

\newcom{\al}{\alpha}
\newcom{\be}{\beta}
\newcom{\del}{\delta}

\newcom{\ga}{\gamma}
\newcom{\Ga}{\Gamma}

\newcom{\Lam}{\Lambda}
\newcom{\lam}{\lambda}

\newcom{\Om}{\Omega}
\newcom{\om}{\omega}

\newcom{\Si}{\Sigma}
\newcom{\si}{\sigma}
\newcom{\s}{\varsigma}

\newcom{\tht}{\theta}
\newcom{\dtri}{\nabla}
\newcom{\tri}{\triangle}

\newcom{\f}{\frac}
\newcom\na{\nabla}
\newcom{\Del}{\Delta}
\newcom{\ep}{\epsilon}

\newcom{\p}{\partial}

\newcom{\uep}{{\bf u}_{\epsilon}}
\newcom{\nep}{n_{\epsilon}}
\newcom{\cep}{c_{\epsilon}}

\newcom{\beq}{\begin{equation}}
\newcom{\eeq}{\end{equation}}
\newcom{\ben}{\begin{eqnarray}}
\newcom{\een}{\end{eqnarray}}
\newcom{\beno}{\begin{eqnarray*}}
\newcom{\eeno}{\end{eqnarray*}}
\newcom{\bal}{\begin{aligned}}
\newcom{\eal}{\end{aligned}}

\author{Hao Li$^1$\thanks{E-mail: lihao\_pde@163.com},
\quad Zhaoyin Xiang$^1$\thanks{E-mail: zxiang@uestc.edu.cn}, 
\quad Xiaoqian Xu$^2$\thanks{E-mail: xiaoqian.xu@dukekunshan.edu.cn}\\
\emph{\small  $1$. School of Mathematical Sciences, University of Electronic Science and Technology of China, }\\
\emph{\small Chengdu 611731, China}\\
\emph{\small  $2$. Zu Chongzhi Center for Mathematics and Computational Sciences, Duke Kunshan University, }\\
\emph{\small Suzhou 215316, China}
  }

\title{Suppression of blow-up in Patlak-Keller-Segel-Navier-Stokes system via the Poiseuille flow}

\date{}

\begin{document}
\maketitle
{\small
{\textbf{Abstract:}   In this paper, we investigate the two-dimensional Patlak-Keller-Segel-Navier-Stokes system  perturbed around the Poiseuille flow $(Ay^2, 0)^\top$ and show that the solutions to this system are global in time if the Poiseuille flow is sufficiently strong in the sense of  amplitude $A$ large enough.   This seems to be the first result showing that the Poiseuille flow can suppress the chemotactic blow-up of the solution to chemotaxis-fluid system. Our proof will be based on a weighted energy method together with the linear enhanced dissipation established by Coti Zelati-Elgindi-Widmayer (Commun. Math. Phys. 378 (2020) 987-1010).  

\smallskip 

{\textbf{Keywords:} {Patlak-Keller-Segel-Navier-Stokes system; Suppression of blow-up;  Poiseuille flow; Enhanced dissipation}

\smallskip 

{\textbf{AMS (2000) Subject Classifications:}} 35K55; 35Q92; 35Q35; 92C17 }

 \section{Introduction}
\renewcommand{\theequation}{\thesection.\arabic{equation}}
\newtheorem{Definition}{Definition}[section]
\newtheorem{Theorem}{Theorem}[section]
\newtheorem{Proposition}{Proposition}[section]
\newtheorem{Lemma}{Lemma}[section]
\newtheorem{Remark}{Remark}[section]
\newtheorem{Corollary}{Corollary}[section]
\numberwithin{equation}{section}

\quad In this paper, we consider the two-dimensional parabolic-elliptic Patlak-Keller-Segel system coupled with the Navier-Stokes equations
\begin{equation}\label{sys11}
\left\{
\begin{split}
\partial_tn  +  {\bf v}\cdot\nabla n -\Delta n   & = - \nabla\cdot(n\nabla c),\\
 -\Delta c & = n-c,\\
\partial_t{\bf v} + {\bf v}\cdot\nabla{\bf v}   + \nabla p - \Delta{\bf v}   & =  n\nabla\phi,\\
\nabla\cdot{\bf v} & = 0
\end{split}
\right.
\end{equation}
posed on the boundary-less domain $\mathbb T\times \mathbb R$, where $\mathbb T=[0,2\pi)$ is a periodic interval,  and supplemented with initial conditions 
\[
n(x, y, 0) = n_{in}(x, y), \quad {\bf v}(x, y, 0) = {\bf v_{\it in}}(x, y) \qquad (x, y)\in \mathbb T\times \mathbb R. 
\]
Here,  the unknowns $n$, $c$, ${\bf v}$ and $p$ stand for the density distribution of cells (e.g., slime mold amoebae Dictyostelium discoideum), the chemical concentration, the fluid velocity vector field and the associated scalar pressure of an incompressible fluid of uniform density, respectively, while $\phi$ is a given potential function of the gravitational field. In this paper, we will  take $\phi(x, y)=y$ for simplicity.

\subsection{Background and literature review}  

\quad If ${\bf v} = {\bf 0}$ in model \eqref{sys11},  it can be reduced to the classical Patlak-Keller-Segel  system
\begin{equation}\label{syscks}
\left\{
\begin{split}
\partial_tn  - \Delta n  & = -\nabla\cdot(n\nabla c),\\
 -\Delta c & = n-c,
\end{split}
\right.
\end{equation}
which was introduced in 1953 by Patlak \cite{P} and then  Keller-Segel \cite{KS}  in 1970 to describe the collective motion of cells.  In system \eqref{syscks}, the first equation takes into account that the motion of cells is driven by the steepest increase in the concentration of chemoattractant while follows a Brownian motion due to external interactions, and  the second one takes into account that cells are producing the chemoattractant themselves while this is diffusing into the environment. In the last five decades, the classical Patlak-Keller-Segel system \eqref{syscks}, a quite interesting mathematical model, has attracted considerable attentions and gained plenty of researches by mathematicians and biologists. One of the most well-known characteristics of system \eqref{syscks} is that the solution to this system may blow-up in finite time in dimensions larger than one.  Precisely, in one-dimensional case, the solution is global well-posedness, while in the two-dimensional case,  the global well-posedness was ensured when  the initial total mass of cells $M:=\|n_{in}\|_{L^1}$ is less than or equal to $8\pi$ \cite{BC, BD, CPZ} but otherwise the finite time blow-up happens \cite{WL,N}. Furthermore,  for any small initial mass,  solutions blowing up  in finite time can be constructed in the higher dimensional case \cite{Win2013}.

In more realistic nature, cells usually live in viscous fluids, in which  cells and chemical stimuli are transported with the fluids, and meanwhile the motion of the fluids is under the influence of gravitational forcing generated by aggregation of cells. To describe the evolution  of such coupled biological dynamics,  Tuval et al. \cite{T}  carried a detailed experiment in a water drop sitting on a glass surface containing Bacillus subtilis, where oxygen diffuses into the drop through the fluid-gas interface.  Mathematically, they proposed a convective chemotaxis system for the oxygen-consuming Bacillus subtilis coupled with the incompressible Navier-Stokes equations subjected to a gravitational force that is proportional to the relative surplus of the cell density compared to the water density.  Then in the past decade, there have been large amounts of literatures on this system (see e.g. \cite{CKL, DL, LL, LO, WWX2, W1, W2, W3, ZZ}).  Roughly speaking, these results apply to uncoupled equations in which both the chemotactic equation and the fluid equation have globally well-posedness solutions. On the other hand, some experimental observations also illustrate that in certain cases of chemotactic movement in flowing environments the mutual influence between cells and fluid may be significant. A typical example is relevant to broadcast spawning phenomena in which an effective mixing triggered by chemotaxis in flowing fluids is indispensable for successful coral fertilization (see e.g. \cite{Coll}).  The signal production system \eqref{syscks} thereby coupled with the  incompressible Navier-Stokes equations leads to model \eqref{sys11}, which covers the recent modeling approaches for biomixing performed by Espejo-Suzuki \cite{ES} on the basis of a work by Kiselev-Ryzhik \cite{KR}.  For comprehensive results regarding the Patlak-Keller-Segel-Navier-Stokes system, we refer the interested readers to \cite{GH, KMS2, LWZ, LA, WWX1, W4} and the references therein.

An interesting question thereby is whether one can suppress the finite-time chemotactic blow-up via the stabilizing effect of the moving fluid.  Recently, some important progresses have been made in showing the prevention of the chemotactic blow-up by the presence of fluid flow as the further development of Kiselev and the third author \cite{KX}, where  it was shown that for any given nonnegative smooth data  $n_{in}$ in $\mathbb{T}^d$, there exist smooth incompressible flows ${\bf v}$ such that the unique solution $n$ of system 
\begin{equation}\label{ks-fluid}
\left\{
\begin{split}
\partial_tn  +  {\bf v}\cdot\nabla n -\Delta n   & = - \nabla\cdot(n\nabla c),\\
 -\Delta c & = n - \overline{n}  
\end{split}
\right.
\end{equation}
with $\overline{n}$ denoting the average of $n$  is globally regular in time.  Indeed, this analysis was generalized in \cite{IXZ} to a broader class of fluid vector fields. Furthermore, Bedrossian-He \cite{BH} and He \cite{H} showed that strong shear flows can prevent the blow-up through a fast dimension reduction process.   The same goal was reached by exploiting the fast-spreading scenario of the hyperbolic and quenching shear flows (see \cite{HT, HTZ}).  In system \eqref{ks-fluid},  the fluid velocity fields ${\bf v}$ are passive because the motion of the fluids is not described there.  If there is active coupling between the cell dynamics and the fluid motion, the only known results owe to Zeng-Zhang-Zi \cite{ZZZ} and He \cite{H2}, where the authors showed that the solution to the two dimensional Patlak-Keller-Segel-Navier-Stokes system near the Couette flow $(Ay, 0)^\top$ in $\mathbb{T} \times\mathbb{R}$ with  $A$ large enough stays globally regular.

\subsection{Problem setting and main results}  

\quad In this paper, we investigate the two-dimensional Patlak-Keller-Segel-Navier-Stokes system \eqref{sys11} perturbed around the Poiseuille flow $(Ay^2, 0)^\top$, a stationary solution to the Navier-Stokes equations.  The motivation of this consideration is two-fold. On the one hand,  the Poiseuille flow is the simplest non-trivial example of a shear flow on  $\mathbb T\times \mathbb R$ besides the Couette flow and moreover is a prototypical example of a strictly convex and nonmonotone shear flow in  the wider physical context \cite{OK}. On the other hand,  from the biomedical aspect,  many bioconvection models are based on basic principles, such as Poiseuille blood flow through the venule, fundamental solutions of the diffusion-reaction equation for the concentration field of pathogen-released chemokines, and linear chemotaxis of the leukocytes \cite{VMG}. Indeed, the fact that cell accumulation drives the fluid to sink more quickly, resulting in focused plume structures that are evident in downwelling Poiseuille flow, has been demonstrated in a classical experiment \cite{Ke},  and thereby the stability of a two-dimensional plume in the Navier–Stokes equations coupled with a micro-organism conservation equation has been numerically investigated in e.g. \cite{GhH}.

We will show that the solutions to system \eqref{sys11} are global in time if the Poiseuille flow  $(Ay^2, 0)^\top$ is sufficiently strong in the sense of  amplitude $A$ large enough for suitably small initial vorticity but without any smallness restriction on the initial cell mass  $M=\|n_{in}\|_{L^1}$. For this purpose,  we take a similar to the work \cite{ZZZ} and first introduce the perturbation 
\[
{\bf u}(t, x, y) = {\bf v}(t, x, y) - (Ay^2, 0)^\top, \qquad \mathfrak{p} (t, x, y) = p(t, x, y) - 2Ax.
\]
Then system \eqref{sys11} can be rewritten as 
\begin{equation}\label{2}
\left\{
\begin{split}
\partial_t{n}  - \Delta n + Ay^2\partial_x n  & = -\nabla\cdot\big(n\nabla c\big) - {\bf u} \cdot\nabla n, \\
-\Delta c & =n-c,\\
\partial_t{\bf u}  - \Delta{\bf u} + A y^2\partial_x {\bf u} + (2Ayu^2,0)^\top  & =  - {\bf u} \cdot\nabla {\bf u} - \nabla \mathfrak{p}  + (0, n)^\top,\\
\nabla\cdot\ {\bf u} & = 0
\end{split}
\right.
\end{equation}
in  $\mathbb T\times \mathbb R$.  By setting ${\bf u}  := (u^1, u^2)$ and  $\nabla^\bot  :=(-\partial_y,\partial_x)$ the rotation of the gradient, we see that  the vorticity $\omega  := \nabla^\bot\cdot{\bf u} = \partial_xu^2 -\partial_y u^1 $ satisfies the scalar equation  
\[
\partial_t{\omega} - \Delta{ \omega} + Ay^2\partial_x \omega - 2Au^2 
 = - {\bf u} \cdot\nabla {\bf \omega} +  \partial_x n, \qquad 
{\bf u}  =  \nabla^\bot \Delta^{-1}\omega. 
\]
Thus after the time rescaling $t\to A^{-1}t$,  we obtain the following equivalence of system \eqref{2}:    
\begin{equation}\label{3}
\left\{
\begin{split}
\partial_{t}n  - \f{1}{A}\Delta n + y^2\partial_x n   & = - \f{1}{A}\nabla\cdot\big(n\nabla c\big) - \f{1}{A}{\bf u}\cdot\nabla n, \\
 -\Delta c & = n-c,\\
\partial_t \omega - \f{1}{A}\Delta \omega + y^2\partial_x \omega - 2\partial_x \Delta^{-1}\omega   & = - \f{1}{A}{\bf u} \cdot\nabla \omega +   \f{1}{A}\partial_x n,  \\
{\bf u} & = \nabla^\bot \Delta^{-1}\omega
\end{split}
\right.
\end{equation}
in  $\mathbb T\times \mathbb R$. System \eqref{3} will be closed by imposing the initial conditions 
\[
n(x, y, 0) = n_{in}(x, y), \quad \omega (x, y, 0) = \omega_{\it in}(x, y) \qquad (x, y)\in \mathbb T\times \mathbb R  
\]
with $\omega_{\it in} = \nabla^\bot\cdot{\bf u}_{\it in}  =  \nabla^\bot\cdot{\bf v}_{\it in} + 2Ay$. 

 To present our conclusion precisely, we introduce the weighted $L^2$ space $X$ normed by 
\[
\|f\|^2_{X}=\|f\|^2_{L^2}+\|yf\|^2_{L^2},
\]
which is  inspired by  Coti Zelati-Elgindi-Widmayer \cite{ZEW} and arises as a natural energy of the system. Our main result is as follows.
\begin{Theorem}\label{global regular}
Assume that the initial data satisfies $n_{in}\in L^1\cap X \cap L^\infty({\mathbb T}\times{\mathbb R})$ with  $\partial_xn_{in}\in X$ and  ${\bf u}_{in}\in H^1({\mathbb T}\times{\mathbb R})$ with $\omega_{in}\in X.$ Then there exists a positive constant $A_0$ relying on $\|n_{in}\|_{L^1\cap X\cap L^\infty}$, $\|\partial_xn_{in}\|_X$, $\|{\bf u}_{in}\|_{L^2}$, but not on $\|\omega_{in}\|_X$, such that if $A>A_0$ and $\|\omega_{in}\|_X\leq A^{-\f{3}{4}}$,  then the solution $(n, c, {\bf u})$  to system \eqref{3} is  global in time.
\end{Theorem}

\begin{Remark}
Theorem \ref{global regular} gives an explicit size in terms of the power of $A$, such that initial vorticity below this threshold yield global regular solutions that exhibit enhanced dissipation. To the best of authors' knowledge, this is the first result showing that the Poiseuille flow can suppress the chemotactic blow-up of the solution to Patlak-Keller-Segel-Navier-Stokes system. Moreover, our result will be proved via a very simple energy method based on the linear enhanced dissipation established by \cite{ZEW}.   
\end{Remark}

\subsection{Key steps}  

\quad  As far as we know, the properties of the $x$-independent part and the $x$-dependent part of the solutions to system \eqref{3} are quite disparate because the former one does not mix a bit. Hence, we will investigate the zero mode and the nonzero modes, separately. For this purpose, let us define
\begin{equation}\nonumber
P_0f = f_0 :=\f{1}{2\pi}\int_{\mathbb T}f(x,y)dx\qquad {\mathrm {and}}\qquad P_{\neq}f = f_{\neq}: =f-f_0,
\end{equation}
for any given function $f$, which correspond to the zero mode and the nonzero modes of $f$, respectively.  That is,  $f_0$ stands for the orthogonal projection of $f$ onto the kernel of the shear flow, while $f_{\neq}$ takes into account   the projection onto the orthogonal complement in $L^2$.   In this way, projecting orthogonally the first component of vector equation \eqref{2}$_3$ yields that
\begin{equation}\label{015}
\partial_tu_0^1-\f{1}{A}\partial_{yy}u_0^1=-\f{1}{A}\partial_y(u_{\neq}^1u_{\neq}^2)_0.
\end{equation}

 As mentioned before,  the Poiseuille flow is nonlinear and is not strictly monotone. This results in a much faster decay than the regular dissipative time-scale and in some new terms which do not appear in \cite{BH, H, H2, ZZZ}.   Thus the  method we employ here will be different from that of  \cite{BH, H, H2, ZZZ} and be based on a  weighted  energy method.  Precisely, we will bound the nonzero modes $(n_{\neq}, c_{\neq}, \omega_{\neq})$ via the  weighted energy estimates as well as the semigroup estimate established by Coti Zelati-Elgindi-Widmayer \cite{ZEW}.    Indeed, since the gradient of the Poiseuille flow is  unbounded as $|y|\to +\infty$, the weight $L^2$ norm arises as a natural energy in our analysis and ensures us to obtain a finer analysis on the nonlinear term for closing the desired energy estimate.  As usual, for instance, the estimate for the difficulty chemotaxis term $\nabla\cdot(n\nabla c)$ in equation \eqref{3}$_1$ requires us to investigate  the evolution of $\nabla c$ and to bound $\|n\|_{L^\infty L^\infty}$.  In particular, the  appearance of the linear term $\partial_xn_{\neq}$ in the equation of $n_{\neq}$ may generate large growth in $\partial_yn_{\neq}$, which amplifies the destabilization effect of the Poiseuille flow. Fortunately, this obstacle can be overcome by using  the weighted estimates  $\|\partial_xn_{\neq}\|_{L^\infty X}$.  Precisely,  $\|\partial_xn_{\neq}\|_{L^\infty X}$ and equation \eqref{3}$_2$ entail the corresponding estimates of $\nabla c$ (see Lemma \ref{LH35} and Lemma \ref{LH36}), while the $L^\infty L^\infty$ estimate of $n$ can be obtained by taking a similar strategy as \cite{ZZZ} to use the Moser-Alikakos iteration \cite{A} as long as $\|\nabla c\|_{L^\infty L^4}$ is controlled.  We remark that the weighted norm of $x$-dependent part of ${\bf u}_{\neq}$ will be controlled by $\omega_{\neq}$ thanks to the commutator estimates (Lemma \ref{LH34}).  Finally, Theorem \ref{global regular} will be proved by a bootstrap argument with the help of the enhanced dissipation of  the Poiseuille flow.

\smallskip

{\bf Notations:}

\begin{itemize}

\item[\bf(1)] Throughout the paper, we will denote by $C$ the positive constant being independent of $t$, $A$ and initial data, which might be vary from line to line.

\item[\bf(2)] Given two operators $\mathcal A$ and $\mathcal B$, the commutator relation $[\mathcal A, \mathcal B]$ will be defined by
\[
[\mathcal A,\mathcal B]f:=\mathcal A(\mathcal Bf)-\mathcal B(\mathcal Af)
\]
for suitable function $f$.

\item[\bf(3)] The Sobolev spaces are defined in a standard manner: for $1\leq p, q\leq\infty$ and $k\in{\mathbb N}$:
\[
\|f\|_{L^p}:=\|f\|_{L^p({\mathbb T}\times{\mathbb R})},\qquad 
W^{k, p}:=\{f\in L^p:\partial^\alpha f\in L^p,\,\mathrm {for} \ \mathrm {all} \ |\alpha|\leq k\}. 
\]
In particular, $H^k:=W^{k, 2}$.  For a function of space and time $f=f(t, x)$, we denote
\[
\|f\|_{L^q W^{k,p}}:=\big\|\|f\|_{W^{k, p}}\big\|_{L^q(0,T)}.
\]

\item[\bf(4)] We denote by $M$ the total mass $\|n(t)\|_{L^1}$, which is conserved because of the equation ${\eqref{sys11}}_1$. That is, 
\[
M:=\|n(t)\|_{L^1}=\|n_{in}\|_{L^1}.
\]

\end{itemize}

The rest of this paper is organized as follows. In Section \ref{prelim}, we decompose system \eqref{sys11} into two subsystems involving  the $x$-independent part and $x$-dependent part, respectively, and present some basic estimates for the solution component $c_0$, $c_{\neq}$ and  ${\bf u}_{\neq}$. The known linear enhanced dissipation estimate will also be stated in this section. Then in Section \ref{bootstrap}, we establish the key bootstrap estimates.  Finally, Theorem \ref{global regular} will be proved in Section \ref{proof}.

\section{Preliminaries}\label{prelim}

\quad In this section, we present some basic preliminaries.   Considering that the enhanced dissipation does not act in the nullspace of the advection term, we decompose system \eqref{3} into two subsystems involving  the $x$-independent part and $x$-dependent part, respectively: 
\begin{equation}\label{AH11}
\left\{
\begin{split}
\partial_t n_0 - \f{1}{A}\partial_{yy}n_0 & = -\f{1}{A}\Big(\partial_y(n_{\neq}\partial_yc_{\neq})_0 +\partial_{y}(n_0\partial_y c_0)\Big) - \f{1}{A}\partial_y(u_{\neq}^2n_{\neq})_0,   \\
-\partial_{yy}c_0 & = n_0-c_0, \\
\partial_t \omega_0-\f{1}{A}\partial_{yy}\omega_0 & = -\f{1}{A}\partial_y(u_{\neq}^2\omega_{\neq})_0, \\
{\bf u_0} & = \big(-\partial_y(\partial_{yy})^{-1}\omega_0,0\big)
\end{split}
\right.
\end{equation}
and
\begin{equation}\label{qh}
\left\{
\begin{split}
\partial_t n_{\neq} - \f{1}{A}\Delta n_{\neq} + y^2\partial_x n_{\neq} & = -\f{1}{A}\Big(\nabla\cdot(n_{\neq}\nabla c_{\neq})_{\neq}+\nabla\cdot(n_0\nabla c_{\neq})+\partial_y(n_{\neq}\partial_y c_0)\Big) \\
&\qquad  \quad  -\f{1}{A}\Big(\nabla\cdot({\bf u_{\neq}}n_{\neq})_{\neq}+\nabla\cdot({\bf u_0 }n_{\neq}) +\nabla\cdot({\bf u_{\neq}}n_0)\Big),&\\
-\Delta c_{\neq} & = n_{\neq}-c_{\neq},&\\
\partial_t\omega_{\neq} - \f{1}{A}\Delta\omega_{\neq} + y^2\partial_x \omega_{\neq} - 2\partial_x \Delta^{-1}\omega_{\neq} & = -\f{1}{A}\Big(\nabla\cdot({\bf u_{\neq}} \omega_{\neq})_{\neq}
+\nabla\cdot({\bf u_0 \omega_{\neq}}) + \nabla\cdot({\bf u_{\neq}}\omega_0)\Big) \\
&\qquad \quad  +\f{1}{A}\partial_x n_{\neq}, \\
{\bf u_{\neq}} & =  \nabla^\bot\Delta^{-1}\omega_{\neq}.
\end{split}
\right.
\end{equation}
The following lemma will provide some useful estimates on  the $x$-independent part of $c$.
\begin{Lemma}\label{LH35}
Let $c_0$ be the zero mode of $c$ and satisfy
\begin{equation}\label{c1}
-\partial_{yy} c_0+c_0=n_0.
\end{equation}
Then for all $t\geq 0$, it holds that
\begin{equation}\label{34316}
\|\partial_{yy}c_0\|_X^2 + \|\partial_{y}c_0\|_X^2 + \|c_0\|_X^2 \le C\|n_0\|_X^2
\end{equation}
for some universal positive constant $C,$ which in particular implies that
\[
\|\partial_{y}c_0\|_{L^4}\leq C\|n_0\|_{L^2}\qquad {\rm{and}}\qquad \|\partial_{yy}c_0\|_{L^\infty}+\|\partial_{y}c_0\|_{L^\infty}
\le C\big(\|n_0\|_{L^2}+\|n_0\|_{L^\infty}\big).
\]
\end{Lemma}
{\bf Proof.} Multiplying equation \eqref{c1} by $c_0,$ and using the integration by parts and the Young inequality, we have 
\[
\|\partial_{y}c_0\|^2_{L^2}+\|c_0\|^2_{L^2}=\int_{\mathbb T\times \mathbb R}n_0c_0dxdy\leq\f{1}{2}\|n_0\|^2_{L^2}+\f{1}{2}\|c_0\|^2_{L^2},
\]
which yields
\begin{equation}\label{3411}
2 \|\partial_yc_0\|^2_{L^2} + \|c_0\|^2_{L^2}\le \|n_0\|^2_{L^2}.
\end{equation}
Similarly, we can multiply \eqref{c1} by $-\partial_{yy}c_0$ and use the integration by parts to obtain
\[
\|\partial_{yy}c_0\|^2_{L^2}+\|\partial_{y}c_0\|^2_{L^2}=-\int_{\mathbb T\times\mathbb R}n_0\partial_{yy}c_0dxdy
\leq \f{1}{2}\|n_0\|^2_{L^2}+\f{1}{2}\|\partial_{yy}c_0\|^2_{L^2},
\]
which implies
\begin{equation}\label{3412}
\|\partial_{yy}c_0\|^2_{L^2} + 2\|\partial_yc_0\|^2_{L^2}\le \|n_0\|^2_{L^2}.
\end{equation}
Combining \eqref{3411} with \eqref{3412}, we can show the estimate for the non-weighted parts in  \eqref{34316}: 
\begin{equation}\label{3343}
\|\partial_{yy}c_0\|_{L^2}^2 + \|\partial_{y}c_0\|_{L^2}^2 + \|c_0\|_{L^2}^2 \le 2\|n_0\|_{L^2}^2,
\end{equation}
which together with the Gagliardo-Nirenberg inequality on $\mathbb R$ also implies that
\[
\|\partial_y c_0\|_{L^4}\leq C\|\partial_y c_0\|_{L^2}^{\frac{3}{4}}\|\partial_{yy}c_0\|_{L^2}^{\frac{1}{4}}\leq C\|n_0\|_{L^2}
\]
and
\[
\|\partial_y c_0\|_{L^\infty}\leq C\|\partial_y c_0\|_{L^2}^{\frac{1}{2}}\|\partial_{yy}c_0\|_{L^2}^{\frac{1}{2}}\leq C\|n_0\|_{L^2}.
\]
By using equation \eqref{c1}, the Minkowski inequality, the Gagliardo-Nirenberg inequality on $\mathbb R$ again and \eqref{3411}, we also have
\[
\|\partial_{yy} c_0\|_{L^\infty}
\le \|c_0\|_{L^\infty}+\|n_0\|_{L^\infty} 
\le C\|c_0\|^{\f{1}{2}}_{L^2}\|\partial_yc_0\|^{\f{1}{2}}_{L^2}+\|n_0\|_{L^\infty}
\le C\big(\|n_0\|_{L^2}+\|n_0\|_{L^\infty}\big).
\]

It remains to show the weighted estimate in  \eqref{34316}.  For this purpose, we multiply equation \eqref{c1} by $y^2c_0$ and use the integration by parts to obtain
\begin{align*}
\|y\partial_yc_0\|^2_{L^2}+\|yc_0\|^2_{L^2}
&=-2\int_{\mathbb T\times\mathbb R}yc_0\partial_yc_0dxdy+\int_{\mathbb T\times\mathbb R}y^2c_0n_0dxdy\\
& =\int_{\mathbb T\times\mathbb R}c_0^2dxdy+\int_{\mathbb T\times\mathbb R}y^2c_0n_0dxdy\\
& \le \|c_0\|^2_{L^2}+\f{1}{2}\|yc_0\|^2_{L^2}+\f{1}{2}\|yn_0\|^2_{L^2},
\end{align*}
which together with \eqref{3343} yields that
\begin{equation}\label{33431}
2\|y\partial_yc_0\|^2_{L^2}+\|yc_0\|^2_{L^2}\leq 2\|c_0\|^2_{L^2}+\|yn_0\|^2_{L^2}\leq C\|n_0\|^2_{X}.
\end{equation}
Similarly, multiplying equation \eqref{c1} by $-y^2\partial_{yy}c_0,$ we can deduce that
\begin{align*}
\|y\partial_{yy}c_0\|^2_{L^2}+\|y\partial_{y}c_0\|^2_{L^2}&=-2\int_{\mathbb T\times\mathbb R}yc_0\partial_yc_0dxdy-\int_{\mathbb T\times\mathbb R}y^2\partial_{yy}c_0n_0dxdy\\
&=\int_{\mathbb T\times\mathbb R}c_0^2dxdy-\int_{\mathbb T\times\mathbb R}y^2\partial_{yy}c_0n_0dxdy\\
&\leq \|c_0\|^2_{L^2}+\f{1}{2}\|y\partial_{yy}c_0\|^2_{L^2}+\f{1}{2}\|yn_0\|^2_{L^2},
\end{align*}
which together with \eqref{3343} again implies that
\begin{equation}\label{33432}
\|y\partial_{yy}c_0\|^2_{L^2}+2\|y\partial_yc_0\|^2_{L^2}\leq 2\|c_0\|^2_{L^2}+\|yn_0\|^2_{L^2}\leq C\|n_0\|^2_{X}.
\end{equation}
Collecting \eqref{3343}, \eqref{33431} and \eqref{33432}, we complete the proof of Lemma \ref{LH35}.\qquad $\Box$

\smallskip 

To establish similar estimates for the $x$-dependent part of $c$, we need  the following anisotropic Sobolev inequality.

\begin{Lemma}[Lemma 3.3 in \cite{ZZZ}]\label{LH32}
Let $f$ be defined on $\mathbb T\times \R$ and satisfy $f_{\neq}\in {\dot {H^{1}}(\mathbb T\times \mathbb R)}$ and $\partial_x f_{\neq}\in {\dot {H^{1}}(\mathbb T\times \mathbb R)}$.  Then for each  $\theta\in(0, 1]$, there exists a constant $C>0$ relying on $\theta$ such that
\[
\|f_{\neq}\|_{L^\infty}\leq C \|\nabla f_{\neq}\|_{L^2}^{1-\theta}\|\nabla \partial_xf_{\neq}\|_{L^2}^\theta.
\]
\end{Lemma}

\begin{Lemma}\label{LH36}
Let $c_{\neq}$ be the nonzero modes of $c$ and satisfy
\begin{equation}\label{c2}
-\Delta c_{\neq}+c_{\neq}=n_{\neq}.
\end{equation}
Then for all $t\geq 0$, it holds that
\begin{equation}\label{35320}
\|D^2 c_{\neq}\|_{X}^2 + \|\nabla c_{\neq}\|_X^2 + \|c_{\neq}\|_X^2 \le C\|n_{\neq}\|_X^2
\end{equation}
for some universal positive constant $C$,  which in particular implies that
\[
\|\nabla c_{\neq}\|_{L^4}\leq C\|n_{\neq}\|_{L^2}\qquad {\rm{and}} \qquad \|\Delta c_{\neq}\|_{L^\infty}+\|\nabla c_{\neq}\|_{L^\infty}
\le  C\big(\|n_{\neq}\|_{L^2}+\|n_{\neq}\|_{L^\infty}+\|\partial_xn_{\neq}\|_{L^2}\big). 
\]
\end{Lemma}
{\bf Proof.} We will prove the desired conclusion in a similar fashion as in the proof of Lemma \ref{LH35}. Indeed, multiplying equation \eqref{c2} by $c_{\neq}$ and using the integration by parts, we have 
\[
\|\nabla c_{\neq}\|^2_{L^2}+\|c_{\neq}\|^2_{L^2}=\int_{\mathbb T\times\mathbb R}n_{\neq}c_{\neq}dxdy\leq \f{1}{2}\|n_{\neq}\|^2_{L^2}+\f{1}{2}\|c_{\neq}\|^2_{L^2},
\]
which implies
\begin{equation}\label{3511}
2\|\nabla c_{\neq}\|^2_{L^2} + \|c_{\neq}\|^2_{L^2} \le \|n_{\neq}\|^2_{L^2}.
\end{equation}
Similarly, we can multiply equation \eqref{c2} by $-\Delta c_{\neq}$ to deduce that
\[
\|\Delta c_{\neq}\|^2_{L^2}+\|\nabla c_{\neq}\|^2_{L^2}=-\int_{\mathbb T\times\mathbb R}n_{\neq}\Delta c_{\neq}dxdy
\leq \f{1}{2}\|n_{\neq}\|^2_{L^2}+\f{1}{2}\|\Delta c_{\neq}\|^2_{L^2},
\]
which yields
\begin{equation}\label{3512}
\|\Delta c_{\neq}\|^2_{L^2} + 2\|\nabla c_{\neq}\|^2_{L^2}\le \|n_{\neq}\|^2_{L^2}.
\end{equation}
Combining \eqref{3511} with \eqref{3512}, we  obtain
\begin{equation}\label{3611}
\|\Delta c_{\neq}\|_{L^2}^2 + \|\nabla c_{\neq}\|_{L^2}^2  + \|c_{\neq}\|_{L^2}^2 \le 2\|n_{\neq}\|_{L^2}^2,
\end{equation}
which together with the fact $\|\Delta c_{\neq}\|_{L^2}=\|D^2 c_{\neq}\|_{L^2}$ also gives that
\begin{equation}\label{36111}
\|D^2 c_{\neq}\|_{L^2}^2 = \|\Delta c_{\neq}\|_{L^2}^2 \le 2\|n_{\neq}\|_{L^2}^2. 
\end{equation}
This gives the non-weighted estimate of $c_{\neq}$ in \eqref{35320}.  An application of  the Gagliardo-Nirenberg inequality on $\mathbb T\times\mathbb R$ also gives that 
\[
\|\nabla c_{\neq}\|_{L^4}\leq C\|\nabla c_{\neq}\|_{L^2}^{\frac{1}{2}}\|D^2 c_{\neq}\|_{L^2}^{\frac{1}{2}}\leq C\|n_{\neq}\|_{L^2}
\]
and that 
\[
\|\Delta c_{\neq}\|_{L^\infty}
\le \|c_{\neq}\|_{L^\infty}+\|n_{\neq}\|_{L^\infty} 
\le C\|D^2 c_{\neq}\|^{\f{1}{2}}_{L^2}\|c_{\neq}\|^{\f{1}{2}}_{L^2}+\|n_{\neq}\|_{L^\infty} 
\le C\|n_{\neq}\|_{L^2}+\|n_{\neq}\|_{L^\infty}.
\]
Furthermore, taking advantage of Lemma \ref{LH32}, equation \eqref{c2} and \eqref{3611}, we obtain 
\begin{align*}
\|\nabla c_{\neq}\|_{L^\infty}
&\leq C\|D^2 c_{\neq}\|^{\f{1}{2}}_{L^2}\|D^2\partial_xc_{\neq}\|^{\f{1}{2}}_{L^2}=C\|\Delta c_{\neq}\|^{\f{1}{2}}_{L^2}\|\Delta\partial_xc_{\neq}\|^{\f{1}{2}}_{L^2}\\
&=C \|n_{\neq}-c_{\neq}\|^{\f{1}{2}}_{L^2}\|\partial_xn_{\neq}-\partial_xc_{\neq}\|^{\f{1}{2}}_{L^2}\\
&\leq C\big (\|n_{\neq}\|_{L^2}+\|c_{\neq}\|_{L^2}\big)^{\f{1}{2}}\big(\|\partial_xn_{\neq}\|_{L^2}+\|\partial_xc_{\neq}\|_{L^2}\big)^{\f{1}{2}}\\
&\leq C\big(\|n_{\neq}\|_{L^2}+\|\partial_xn_{\neq}\|_{L^2}\big).
\end{align*}

For the weighted estimates of $c_{\neq}$  in \eqref{35320},  we first multiply equation \eqref{c2} by $y^2c_{\neq}$ and use the integration by parts to get
\begin{align*}
\|y\nabla c_{\neq}\|^2_{L^2}+\|y c_{\neq}\|^2_{L^2} 
& = -2\int_{\mathbb T\times\mathbb R}yc_{\neq}\partial_yc_{\neq}dxdy+\int_{\mathbb T\times\mathbb R}y^2c_{\neq}n_{\neq}dxdy\\
& = \int_{\mathbb T\times\mathbb R}c_{\neq}^2dxdy+\int_{\mathbb T\times\mathbb R}y^2c_{\neq}n_{\neq}dxdy\\
& \le \|c_{\neq}\|^2_{L^2}+\f{1}{2}\|y c_{\neq}\|^2_{L^2}+\f{1}{2}\|yn_{\neq}\|^2_{L^2},
\end{align*}
which together with \eqref{3611} implies that
\begin{equation}\label{36112}
2\|y\nabla c_{\neq}\|^2_{L^2}+\|yc_{\neq}\|^2_{L^2}\leq 2\|c_{\neq}\|^2_{L^2}+\|yn_{\neq}\|^2_{L^2}\leq C\|n_{\neq}\|^2_{X}.
\end{equation}
Similarly, multiplying equation \eqref{c2} by $-y^2\Delta c_{\neq}$ entails that 
\begin{align*}
\|y\Delta c_{\neq}\|^2_{L^2}+\|y\nabla c_{\neq}\|^2_{L^2}&=-2\int_{\mathbb T\times\mathbb R}yc_{\neq}\partial_yc_{\neq}dxdy-\int_{\mathbb T\times\mathbb R}y^2\Delta c_{\neq}n_{\neq}dxdy\\
&=\int_{\mathbb T\times\mathbb R}c_{\neq}^2dxdy-\int_{\mathbb T\times\mathbb R}y^2\Delta c_{\neq}n_{\neq}dxdy\\
&\leq \| c_{\neq}\|^2_{L^2}+\f{1}{2}\|y\Delta c_{\neq}\|^2_{L^2}+\f{1}{2}\|yn_{\neq}\|^2_{L^2}
\end{align*}
and thus that 
\begin{equation}\label{35332}
\|y\Delta c_{\neq}\|^2_{L^2}+2\|y\nabla c_{\neq}\|^2_{L^2}\leq 2\|c_{\neq}\|^2_{L^2}+\|yn_{\neq}\|^2_{L^2}\leq C\|n_{\neq}\|^2_{X}.
\end{equation}
Noticing that 
\begin{align*}
\|y\Delta c_{\neq}\|^2_{L^2}
& =\int_{\mathbb T\times \mathbb R}\big((y\partial_{xx}c_{\neq})^2+(y\partial_{yy}c_{\neq})^2+2y^2\partial_{xx}c_{\neq}\partial_{yy}c_{\neq}\big)dxdy\\
& =\int_{\mathbb T\times \mathbb R}\big((y\partial_{xx}c_{\neq})^2+(y\partial_{yy}c_{\neq})^2\big)dxdy+2\int_{\mathbb T\times \mathbb R}(y\partial_{xy}c_{\neq})^2dxdy+4\int_{\mathbb T\times \mathbb R}y\partial_{x}c_{\neq}\partial_{xy}c_{\neq}dxdy\\
& =\int_{\mathbb T\times \mathbb R}\big((y\partial_{xx}c_{\neq})^2+(y\partial_{yy}c_{\neq})^2\big)dxdy+2\int_{\mathbb T\times \mathbb R}(y\partial_{xy}c_{\neq})^2dxdy-2\int_{\mathbb T\times \mathbb R}(\partial_{x}c_{\neq})^2dxdy,
\end{align*}
we obtain from  \eqref{3611} and \eqref{35332}  that
\begin{align}\label{36113}
\|yD^2 c_{\neq}\|^2_{L^2}\leq
\|y\Delta c_{\neq}\|^2_{L^2}+2\|\partial_xc_{\neq}\|^2_{L^2}\leq C\|n_{\neq}\|^2_X.
\end{align}
Collecting \eqref{3611}, \eqref{36111}, \eqref{36112} and \eqref{36113}, we complete the proof of Lemma \ref{LH36}.\qquad $\Box$

\smallskip

Next, we are going to utilize the Riesz transform and commutator relation to obtain the following inequalities, which will be frequently used in the estimates of ${\bf u_{\neq}}.$
\begin{Lemma}\label{LH34}
Let $\omega_{\neq}$ be the nonzero modes of $\omega$ and ${\bf u_{\neq}}$ be determined by $\eqref{qh}_4$. Then for all $t\ge 0$,  it holds that
\begin{equation}\label{cm341}
\|y {\bf u_{\neq}}\|_{L^2}  + \|y \nabla {\bf u_{\neq}}\|_{L^2} 
\le  C\big(\|\omega_{\neq}\|_{L^2} + \|y\omega_{\neq}\|_{L^2}\big)
\end{equation}
for some universal positive constant $C$.
\end{Lemma}
{\bf Proof.} The key idea is to use the commutator relation
\[
[y,\Delta^{-1}]f:=y\Delta^{-1}f-\Delta^{-1}(yf)=2\Delta^{-2}\partial_yf
\]
for any given function $f,$ which follows from the fact that  $y\Delta^{-1}f=2\Delta^{-2}\partial_yf+\Delta^{-1}(yf)$  due to $\Delta(y\Delta^{-1}f)=2\Delta^{-1}\partial_yf+yf$.

Indeed, it follows from ${\bf u_{\neq}}=\nabla^\bot\Delta^{-1}\omega_{\neq}$ that
\begin{align*}
\|y{\bf u_{\neq}}\|_{L^2}
= \big\|y\Delta^{-1}\nabla^\bot\omega_{\neq}\big\|_{L^2} 
& = \big\|[y,\Delta^{-1}]\nabla^\bot\omega_{\neq} + \Delta^{-1}(y\nabla^\bot\omega_{\neq})\big\|_{L^2} \nonumber \\
& \le \big\|2\Delta^{-2}\partial_y(\nabla^\bot\omega_{\neq})\big\|_{L^2} + \big\|\Delta^{-1}(y\nabla^\bot\omega_{\neq})\big\|_{L^2}   \nonumber\\
& \le 2\|\omega_{\neq}\|_{L^2} + \big\|\Delta^{-1}(y\nabla^\bot\omega_{\neq})\big\|_{L^2}
\end{align*}
which together with the estimate 
\[
\big\|\Delta^{-1}(y\nabla^\bot\omega_{\neq})\big\|_{L^2}
= \big\|\Delta^{-1}\nabla^\bot(y\omega_{\neq})-\Delta^{-1}\big((\nabla^\bot y)\omega_{\neq}\big)\big\|_{L^2}
\le \|y\omega_{\neq}\|_{L^2} + \|\omega_{\neq}\|_{L^2},
\]
implies that 
\[
\|y {\bf u_{\neq}}\|_{L^2} \le 3\|\omega_{\neq}\|_{L^2} + \|y\omega_{\neq}\|_{L^2}.
\]

Taking a similar procedure, we can also deduce 
\begin{align*} 
\|y\partial_y{\bf u_{\neq}}\|_{L^2}
= \big\|y\Delta^{-1}\nabla^\bot\partial_y\omega_{\neq}\big\|_{L^2}
& \le \big\|2\Delta^{-2}\partial_y(\nabla^\bot\partial_y\omega_{\neq})\big\|_{L^2} + \big\|\Delta^{-1}(y\nabla^\bot\partial_y\omega_{\neq})\big\|_{L^2}     \nonumber\\
& \le 2\|\omega_{\neq}\|_{L^2} +\big\|\Delta^{-1}(y\nabla^\bot\partial_y\omega_{\neq})\big\|_{L^2}
\end{align*}
and 
\begin{align*} 
\big\|\Delta^{-1}(y\nabla^\bot\partial_y\omega_{\neq})\big\|_{L^2}
& = \big\|\Delta^{-1}\nabla^\bot\partial_y(y\omega_{\neq})-\Delta^{-1}\nabla^\bot\omega_{\neq}-\Delta^{-1}\big((\nabla^\bot y)\partial_y\omega_{\neq}\big)\big\|_{L^2} \nonumber\\
& \le \big\|\Delta^{-1}\nabla^\bot\partial_y(y\omega_{\neq})\big\|_{L^2} + \big\|\Delta^{-1}\nabla^\bot\omega_{\neq}\big\|_{L^2} + \big\|\Delta^{-1}\partial_y
\omega_{\neq}\big\|_{L^2} \nonumber\\
& \le \|y\omega_{\neq}\|_{L^2}+2\|\omega_{\neq}\|_{L^2}\nonumber
\end{align*}
that
\[
\|y \partial_y{\bf u_{\neq}}\|_{L^2}\leq 4\|\omega_{\neq}\|_{L^2}+\|y\omega_{\neq}\|_{L^2}.
\]

Similarly, we can obtain 
\begin{align}
\|y \partial_x{\bf u_{\neq}}\|_{L^2}
= \big\|\partial_x(y\Delta^{-1}\nabla^\bot\omega_{\neq})\big\|_{L^2}
& \le \big\|\partial_x\big(2\Delta^{-2}\partial_y(\nabla^\bot\omega_{\neq})\big)\big\|_{L^2} + \big\|\partial_x\Delta^{-1}(y\nabla^\bot\omega_{\neq})\big \|_{L^2} \nonumber\\
& \le 2\|\omega_{\neq}\|_{L^2} + \big\|\partial_x\big(\Delta^{-1}\nabla^\bot(y\omega_{\neq})\big)-\partial_x\Delta^{-1}\big((\nabla^\bot y)\omega_{\neq}\big)\big\|_{L^2}\nonumber\\
& \le 3\|\omega_{\neq}\|_{L^2}+\|y\omega_{\neq}\|_{L^2}.\nonumber
\end{align}

Collecting the above estimates, we complete the proof of Lemma \ref{LH34}.  \qquad $\Box$

\smallskip 

 We will end  this section by recalling the linear enhanced dissipation, which will be used to bound the nonzero modes of solutions to system \eqref{sys11}.  For simplicity, we denote the linear operators $\widetilde{\mathcal L}$ and $\mathcal L$ by 
 \[
 {\widetilde{\mathcal L}}:=\f{1}{A}\Delta -y^2\partial_x + 2\partial_x\Delta^{-1}  \qquad \mathrm{and}\qquad {\mathcal L}:=\f{1}{A}\Delta -y^2\partial_x. 
 \]

\begin{Lemma}[Corollary 1.2 in \cite{ZEW}]\label{LH31}
Let $A>1$ and $f_{in}\in X$, and assume that for almost every $y\in \mathbb R$, we have
\[
\int_{\mathbb T} f_{in}(x,y)dx=0.
\]
Then
\[
\big\|e^{\widetilde{\mathcal L}t}f_{in} \big\|_X
\le C_0e^{-\f{\widetilde{\epsilon}_0 t}{A^{\f{1}{2}}(1+\mathrm{log}\,A)}}\| f_{in}\|_{X} \qquad {\rm for\,\,all\,\,\,\,} t\ge 0,
\]
where the constants $ \widetilde{\epsilon}_0\geq \f{1}{20}$ and $C_0\in(1,10)$ can be explicitly computed. In particular, whenever $A>e$, we have
\begin{equation}\label{322}
\big\|e^{\widetilde{\mathcal L}t}f_{in}\big\|_X
\le C_0e^{-{\epsilon}_0\lambda_At}\| f_{in}\|_{X} \qquad {\rm for\,\,all\,\,\,\,} t \ge 0,
\end{equation}
where ${\epsilon}_0 = 2\widetilde{\epsilon}_0 \ge \f{1}{10}$ and 
\[
 \lambda_A=\f{1}{A^{\f{1}{2}}\mathrm{log}\,A}. 
 \] 
\end{Lemma}

\begin{Remark}\label{LH31remark}
  Following the proof of Corollary 1.2 in \cite{ZEW} line by line, we can  verify that the semigroup estimate \eqref{322} still holds for the linear operator $\mathcal L$.  In this case, indeed, we can get the faster dissipation indicator $A^{\f{1}{2}}$. 
\end{Remark}

\section{Bootstrap estimates}\label{bootstrap}

\quad  In this section, we establish the bootstrap estimates. Let us denote $T$ as the end-point of the largest interval $[0,T]$ such that the following assumptions hold:\\

(A-1) Weighted $L^{2}\dot{H^{1}}$ estimate for nonzero modes of $n$:
\[
\f{1}{A}\int_0^t\|\nabla n_{\neq}(\cdot,\tau)\|_X^2d\tau\leq 4 \| (n_{in})_{\neq}\|_X^2;
\]

(A-2) Enhanced dissipation  estimate  for nonzero modes of $n$:
\[
\| n_{\neq}\|_X\leq 4C_0e^{-\epsilon_0\lambda_A t}\| (n_{in})_{\neq}\|_X;
\]

(A-3) Weighted $L^{\infty} L^2$ estimate of $\partial_x n$:
\[
\| \partial_x n\|^2_{L^{\infty}X}=\| \partial_x n_{\neq}\|^2_{L^{\infty}X}\leq 4\|( \partial_x n_{in})_{\neq}\|^2_X;
\]

(A-4) $L^{\infty}L^{\infty}$ estimate of solution $n$:
\[
\|  n\|_{L^{\infty}L^{\infty}}\leq 4C_\infty;
\]

(A-5) Weighted $L^{2}\dot{H^{1}}$ estimate for nonzero modes  of $\omega$:
\[
\f{1}{A}\int_0^t\|\nabla \omega_{\neq}(\cdot,\tau)\|_X^2d\tau\leq 4 \Big(\| (\omega_{in})_{\neq}\|_X^2+{ A^{-\frac{3}{4}} }\Big);
\]

(A-6) Enhanced dissipation  estimate  for nonzero modes of $\omega$:
\[
\| \omega_{\neq}\|_X\leq 4C_0e^{-\epsilon_0\lambda_A t}\Big(\|(\omega_{in})_{\neq}\|_X+{ A^{-\frac{3}{4}} }\Big)
\]
for all $0\leq t \leq T$,  where the constants $ \epsilon_0 \geq \f{1}{10}$ and $C_0\in(1,10)$ are determined by Lemma \ref{LH31}, while the constant $C_\infty\geq1$ is determined by Lemma \ref{LH47}.

\begin{Remark}
Without loss of generality, we will assume that {$0<T\leq\lambda_A^{-\f{1}{4}}$} throughout this section.  
\end{Remark}

Our purpose is to establish some refined bounds from the above assumptions. Precisely, we will show the following proposition. 
\begin{Proposition}\label{LH1}
Assume that the initial data $(n_{in}, {\bf u}_{in})$ satisfy the assumptions of Theorem \ref{global regular}. If the solution $(n, c, {\bf u})$ of system \eqref{2} possesses the bounds {\rm (A-1)-(A-6)}, then there exists a positive constant  $A_0$ depending only on $C_0$, $C_\infty$,  $\|n_{in}\|_{L^1\cap X\cap L^\infty}$,  $\|(\partial_xn_{in})_{\neq}\|_X$  and  $\|{\bf u}_{in}\|_{L^2}$ such that  the following refined bounds hold:  \\ 

{\rm (B-1)} Weighted $L^{2}\dot{H^{1}}$ estimate for nonzero modes of $n$:
\[
\f{1}{A}\int_0^t\|\nabla n_{\neq}(\cdot,\tau)\|_X^2d\tau\leq 2\| (n_{in})_{\neq}\|_X^2; 
\]

{\rm (B-2)} Enhanced dissipation  estimate  for nonzero modes of $n$:
\[
\| n_{\neq}\|_X\leq 2C_0e^{-\epsilon_0\lambda_A t}\| (n_{in})_{\neq}\|_X;
\]

{\rm (B-3)} Weighted $L^{\infty} L^2$ estimate of $\partial_x n$:
\[
\|\partial_x n\|^2_{L^{\infty}X}=\|\partial_x n_{\neq}\|^2_{L^{\infty}X}\leq 2\| (\partial_x n_{in})_{\neq}\|^2_X;
\]

{\rm (B-4)} $L^{\infty}L^{\infty}$ estimate of solution $n$:
\[
\|  n\|_{L^{\infty}L^{\infty}}\le 2C_\infty;
\]

{\rm (B-5)} Weighted $L^{2}\dot{H^{1}}$ estimate for nonzero modes  of $\omega$:
\[
\f{1}{A}\int_0^t\|\nabla \omega_{\neq}(\cdot,\tau)\|_X^2d\tau\leq 2\Big(\|(\omega_{in})_{\neq} \|_X^2+{ A^{-\frac{3}{4}} }\Big);
\]

{\rm (B-6)} Enhanced dissipation  estimate  for nonzero modes of $\omega$:
\[
\| \omega_{\neq}\|_X\leq 2C_0e^{-\epsilon_0\lambda_A t}\Big(\| (\omega_{in})_{\neq}\|_X+{ A^{-\frac{3}{4}} }\Big)
\]
for all $0\leq t \leq T$ whenever $A>A_0$ and $A^{\f{3}{4}}\|\omega_{in}\|_X\leq 1$. 
\end{Proposition}

\begin{Remark}
Due to the restriction $A^{\f{3}{4}}\|\omega_{in}\|_X\leq 1$ in Theorem \ref{global regular}, the assumptions {\rm (A-5)} and {\rm (A-6)} can be replaced by  
\[
\int_0^t\|\nabla \omega_{\neq}\|_X^2d\tau\le { 4A^{\f{1}{4}} } \qquad  \mathrm{and} \qquad 
\| \omega_{\neq}\|_X\le { 4C_0e^{-\epsilon_0\lambda_A t} A^{-\frac{3}{4}} }, 
\]
respectively. Here we remain the forms of {\rm (A-5)} and {\rm (A-6)} since  the optimality of  exponent $\f{3}{4}$ is unclear as mentioned before.
\end{Remark}

\subsection{Zero mode estimates}

\quad The following lemma gives the basic zero mode estimates of $n_0$, $\omega_0$ and $u_0^1$. 
\begin{Lemma}\label{LH41}
Under the assumptions {\rm(A-2)} and {\rm (A-4)-(A-6)}, there exists a universal positive constant $C_\ast$ and a constant $A_1>1$ relying on $C_0$, $C_\infty$ and $\|n_{in}\|_X$ such that if $A>A_1$ and $A^{\f{3}{4}}\|\omega_{in}\|_X\le 1$,  then it holds
\begin{align} 
\|n_0\|^2_{L^\infty X}+\f{1}{A}\|\partial_yn_0\|^2_{L^2 X} 
& \le \|(n_{in})_0\|^2_{X} +  \f{C_*}{A^{\f{1}{16}}} \big(\|(n_{in})_0\|^2_{X}+1\big)
 \le 2 C_\ast\big(\|(n_{in})_0\|^2_{X}+1\big), \label{4142} \\
\|\omega_0\|^2_{L^\infty X}+\f{1}{A}\|\partial_y \omega_0\|^2_{L^2 X}
&  \le \|(\omega_{in})_0\|^2_{X} + \f{C_*}{A^{\f{3}{2}}},   \label{414}  \\
\|u_0^1\|^2_{L^\infty L^\infty}+\|u^1_0\|^2_{L^\infty L^2}
& \le C_\ast \Big(\|{\bf u}_{in}\|^2_{L^2} + \|(\omega_{in})_0\|^2_{L^2} +  \f{1}{A^{\f{3}{2}}} \Big) 
\le  2C_\ast \big(\|{\bf u}_{in}\|^2_{L^2} + 1 \big).   \label{413}
\end{align}
\end{Lemma}
{\bf Proof.}  {\bf Step 1.  Estimates of $n_0$.}  To investigate the non-weighted estimate, we first take the $L^2$ inner product of equation ${\eqref{AH11}}_1$ with $n_0$ to obtain
\begin{align*}
&\frac{1}{2}\f{d}{dt}\|n_0\|^2_{L^2}+\f{1}{A}\|\partial_yn_0\|^2_{L^2}\nonumber\\
& = \f{1}{A}\int_{\mathbb T\times\mathbb R} \Big((n_{\neq}\partial_yc_{\neq})_0 + n_0\partial_yc_0 + (u^2_{\neq}n_{\neq})_0\Big) \partial_yn_0dxdy    \nonumber\\
& \le \f{1}{2A}\|\partial_yn_0\|^2_{L^2} + \f{3}{2A}\Big(\big\|(n_{\neq}\partial_yc_{\neq})_0\big\|_{L^2}^2 + \big\|n_0\partial_yc_0\big\|_{L^2}^2 + \big\|(u^2_{\neq}n_{\neq})_0\big\|_{L^2}^2\Big)\nonumber\\
& \le \f{1}{2A}\|\partial_yn_0\|^2_{L^2} + \f{C}{A} \Big(\|n_{\neq}\|^2_{L^\infty}\|\partial_yc_{\neq}\|^2_{L^2} + \|n_{0}\|^2_{L^\infty} \|\partial_yc_{0}\|^2_{L^2} + \|n_{\neq}\|^2_{L^\infty}\|u^2_{\neq}\|^2_{L^2}\Big)
\end{align*}
for some  universal constant $C$.  Then we see from \eqref{3343}, \eqref{3611} and the fact 
\begin{equation}\label{4121}
\|{\bf u}_{\neq}\|_{L^2}=\|\Delta^{-1}\nabla^\bot\omega_{\neq}\|_{L^2}\leq C\|\omega_{\neq}\|_{L^2}
\end{equation} 
that
\[
\f{d}{dt}\|n_0\|^2_{L^2}+\f{1}{A}\|\partial_yn_0\|^2_{L^2}
\leq \f{C}{A} \Big(\|n_{\neq}\|^2_{L^\infty}\|n_{\neq}\|^2_{L^2} + \|n_0\|^2_{L^\infty}\|n_0\|^2_{L^2} + \|n_{\neq}\|^2_{L^\infty}\|\omega_{\neq}\|^2_{L^2}\Big).
\]
Thus by using the assumptions (A-2),   (A-4),  (A-6) and $\|\omega_{in}\|_X\le A^{-\f{3}{4}}$,  we deduce that
\begin{align}\label{4411}
\f{d}{dt}\|n_0\|^2_{L^2}+\f{1}{A}\|\partial_yn_0\|^2_{L^2}
& \le \f{CC^2_\infty}{A}\big(\|n_{0}\|^2_{L^2}+\|n_{\neq}\|^2_{L^2}+\|\omega_{\neq}\|^2_{L^2}\big)\nonumber\\
& \le \f{CC^2_\infty}{A}\|n_0\|^2_{L^2}+\f{CC_0^2C^2_\infty}{A}\Big(\|(n_{in})_{\neq}\|^2_X+\|(\omega_{in})_{\neq}\|^2_X+{ A^{-\f{3}{2}} }
\Big)\nonumber\\
& \le \f{CC^2_\infty}{A}\|n_0\|^2_{L^2}+\f{CC_0^2C^2_\infty}{A}{ \Big(\|(n_{in})_{\neq}\|^2_X+A^{-\f{3}{2}}
\Big)}\nonumber\\
& \le \f{C}{A^{\f{3}{4}}} \|n_0\|^2_{L^2} + \f{C}{A^{\f{3}{4}}} \big(\|(n_{in})_{\neq}\|^2_X+1\big)
\end{align}
provided that $A>A^*_1$, where $A^*_1$ satisfies ${A^*_1}^{\f{1}{4}}\ge C_0^2C_\infty^2$. It then follows from the Gronwall inequality that 
\[
\|n_0(t)\|^2_{L^2}\leq e^{CA^{-\f{3}{4}} t}\Big(\|(n_{in})_0\|^2_{L^2} + CA^{-\f{3}{4}}t\big(\|(n_{in})_{\neq}\|^2_X+1\big)\Big)
\]
for all $t\in[0, T]$.  Owing to  { $0<T\leq\lambda_A^{-\f{1}{4}}=A^{\f{1}{8}}\mathrm{log}^{\f{1}{4}}A$ }, we can obtain
\begin{align*}
\|n_0(t)\|^2_{L^2}
&\le e^{CA^{-\f{3}{4}} { \lambda_A^{-\f{1}{4}} }}\Big(\|(n_{in})_0\|^2_{L^2} + CA^{-\f{3}{4}}{ \lambda_A^{-\f{1}{4}} }\big(\|(n_{in})_{\neq}\|^2_X+1\big)\Big)\\
&= e^{C { A^{-\f{5}{8}}\mathrm{log}^{\f{1}{4}}A}}\Big(\|(n_{in})_0\|^2_{L^2} + C{ A^{-\f{5}{8}}\mathrm{log}^{\f{1}{4}}A}
\big(\|(n_{in})_{\neq}\|^2_X+1\big)\Big)\\
&\le C\Big(\|(n_{in})_0\|^2_{L^2} + { A^{-\f{3}{8}}}\big(\|(n_{in})_{\neq}\|^2_X+1\big)\Big).
\end{align*}
If we further choose $A^{*}_1$ sufficiently large such that
\begin{equation}\label{ninnonx}
{ {A^{*}_1}^{\f{3}{8}} } > \|(n_{in})_{\neq}\|^2_X+1,
\end{equation}
then we have 
\begin{equation}\label{424481}
\|n_0(t)\|^2_{L^2}\leq C\big(\|(n_{in})_0\|^2_{L^2}+1\big)
\end{equation}
for all $A>A^{*}_1$.  Integrating inequality \eqref{4411} with respect to $t$ on $[0, T]$ and using \eqref{ninnonx} and \eqref{424481}, we obtain 
\begin{align}\label{41413}
\|n_0\|^2_{L^\infty L^2}+\f{1}{A}\|\partial_yn_0\|^2_{L^2L^2}
& \le \|(n_{in})_0\|^2_{L^2} + \f{C}{A^{\f{3}{4}}}\int_0^T\Big(\|n_0\|^2_{L^2} + \|(n_{in})_{\neq}\|^2_X+1\Big)d\tau \nonumber\\
& \le \|(n_{in})_0\|^2_{L^2} + \f{C}{A^{\f{3}{4}}}\big(\|(n_{in})_0\|^2_{L^2} + { A^{\f{3}{8}} }\big)T\nonumber\\
& \le \|(n_{in})_0\|^2_{L^2} +  \f{C}{A^{\f{3}{4}}}\big(\|(n_{in})_0\|^2_{L^2} +{ A^{\f{3}{8}} } \big) { A^{\f{1}{8}}\mathrm{log}^{\f{1}{4}}A}\nonumber\\
& \le \|(n_{in})_0\|^2_{L^2} + { \f{C}{A^{\f{1}{8}}} }\big(\|(n_{in})_0\|^2_{L^2} + 1 \big).
\end{align}

Next, we  deduce the weighted estimates of $n_0$. We multiply equation \eqref{AH11}$_1$ by $y^2n_0$ and use the integration by parts to obtain
\begin{align}\label{4115}
&\frac{1}{2}\f{d}{dt}\|yn_0\|^2_{L^2}+\f{1}{A}\|y\partial_yn_0\|^2_{L^2}-\f{1}{A}\|n_0\|^2_{L^2}\nonumber\\
& = \f{1}{A} \int_{\mathbb T\times\mathbb R} \Big((n_{\neq}\partial_yc_{\neq})_0 + n_0\partial_yc_0\Big) \partial_y(y^2n_0)dxdy + \f{1}{A}\int_{\mathbb T\times\mathbb R}( u^2_{\neq}n_{\neq})_0\partial_y(y^2n_0)dxdy.
\end{align}
It follows from the H{\"o}lder inequality and the Young inequality that
\begin{align*}
& \int_{\mathbb T\times\mathbb R} \Big((n_{\neq}\partial_yc_{\neq})_0 + n_0\partial_yc_0\Big) \partial_y(y^2n_0)dxdy \\
& = \int_{\mathbb T\times\mathbb R}\Big((n_{\neq}\partial_yc_{\neq})_0 + n_0\partial_yc_0\Big)  y^2\partial_yn_0dxdy + 2\int_{\mathbb T\times\mathbb R} \Big((n_{\neq}\partial_yc_{\neq})_0 + n_0\partial_yc_0\Big)  yn_0dxdy\\
& \le \f{1}{3} \|y\partial_yn_0\|^2_{L^2} + \f{3}{2}\big(\big\|y(n_{\neq}\partial_yc_{\neq})_0\big\|^2_{L^2} + \big\|yn_0\partial_yc_0\big\|^2_{L^2} \big) + 2 \big(\big\|y(n_{\neq}\partial_yc_{\neq})_0\big\|_{L^2} + \big\|yn_0\partial_yc_0\big\|_{L^2}\big)\|n_0\|_{L^2} \\
& \le \f{1}{3} \|y\partial_yn_0\|^2_{L^2} + C\big(\|n_{\neq}\|^2_{L^\infty}\|y\partial_yc_{\neq}\|^2_{L^2} +  \|n_0\|^2_{L^\infty}\|y\partial_yc_0\|^2_{L^2} + \|n_0\|^2_{L^2}\big) \\
& \le \f{1}{3} \|y\partial_yn_0\|^2_{L^2} + C\big(\|n_{\neq}\|^2_{L^\infty}\|n_{\neq}\|^2_{X} +  \|n_{0}\|^2_{L^\infty}\|n_0\|^2_{X} + \|n_0\|^2_{L^2}\big)
\end{align*}
thanks to \eqref{34316} and \eqref{35320}. Similarly,  we have 
\begin{align*}
\int_{\mathbb T\times\mathbb R}( u^2_{\neq}n_{\neq})_0\partial_y(y^2n_0)dxdy
& = \int_{\mathbb T\times\mathbb R}(u^2_{\neq}n_{\neq})_0y^2\partial_yn_0dxdy + 2\int_{\mathbb T\times\mathbb R}(u^2_{\neq}n_{\neq})_0yn_0dxdy\\
& \le \f{1}{6}\|y\partial_yn_0\|^2_{L^2}+\f{3}{2}\|y(u^2_{\neq}n_{\neq})_0\|^2_{L^2} +
 2 \|y(u^2_{\neq}n_{\neq})_0\|_{L^1}\|n_0\|_{L^{\infty}} \\
& \le \f{1}{6}\|y\partial_yn_0\|^2_{L^2}  + C\|n_{\neq}\|^2_{L^\infty}\|yu_{\neq}^2\|^2_{L^2} + 
C\|yu^2_{\neq}\|_{L^2}\|n_{\neq}\|_{L^2} \|n_{0}\|_{L^\infty}\\
& \le \f{1}{6}\|y\partial_yn_0\|^2_{L^2}+C\|n_{\neq}\|^2_{L^\infty}\|\omega_{\neq}\|^2_{X}+
C\|\omega_{\neq}\|_{X}\|n_{\neq}\|_{L^2}\|n_{0}\|_{L^\infty}
\end{align*}
due to Lemma \ref{LH34}. Substituting the above inequalities into \eqref{4115} and using the assumptions  (A-2),  (A-4),  (A-6) and the non-weighted estimate \eqref {41413}, we  obtain
\begin{align}\label{41311}
&\f{d}{dt}\|yn_0\|^2_{L^2}+\f{1}{A}\|y\partial_yn_0\|^2_{L^2}\nonumber\\
& \le \f{C}{A}\|n_0\|^2_{L^2} + \f{C}{A} \|n_{\neq}\|^2_{L^\infty} \big(\|n_{\neq}\|^2_{X}
 + \|\omega_{\neq}\|^2_{X}\big) + \f{C}{A}  \|n_0\|^2_{L^\infty}\|n_0\|^2_{X} 
+ \f{C}{A}\|\omega_{\neq}\|_{X}\|n_{\neq}\|_{L^2}\|n_{0}\|_{L^\infty}\nonumber\\
& \le  \f{C}{A}\|n_0\|^2_{L^2}+\f{CC^2_\infty}{A}\big(\|n_{\neq}\|^2_{X}
+ \|\omega_{\neq}\|^2_X +\|n_0\|^2_X \big)\nonumber\\
& \le \f{C}{A}\big(\|(n_{in})_0\|^2_{L^2}+1\big)+\f{CC_0^2C^2_\infty}{A}\big(\|(n_{in})_{\neq}\|^2_X+\|(\omega_{in})_{\neq}\|^2_X+
{ A^{-\f{3}{2}} }\big)
+\f{CC^2_\infty}{A}\|n_0\|^2_{X}\nonumber\\
& \le \f{CC_{\infty}^2}{A}\big(\|(n_{in})_0\|^2_{L^2}+1\big)+\f{CC_0^2C^2_\infty}{A}\big(\|(n_{in})_{\neq}\|^2_X+1\big) + \f{CC^2_\infty}{A}\|yn_0\|^2_{L^2}\nonumber\\
& \le  \f{C}{A^{\f{3}{4}}}+\f{C}{A^{\f{3}{4}}}\|yn_0\|^2_{L^2}
\end{align}
with some universal constant $C$ for all $t\in[0, T]$ provided that  $A>A^{**}_1$, where $A^{**}_1$ satisfies ${A^{**}_1}^{\f{1}{4}} \ge  C^2_{\infty}\|(n_{in})_0\|^2_{L^2} + C_0^2C^2_{\infty}{A^{*}_1}^{\f{1}{4}}$.  Then by the Gronwall inequality and { $0<T\leq\lambda_A^{-\f{1}{4}}=A^{\f{1}{8}}\mathrm{log}^{\f{1}{4}}A$}, we obtain
\begin{align}\label{424482}
\|yn_0(t)\|^2_{L^2}
& \le e^{CA^{-\f{3}{4}}t} \Big( \|y(n_{in})_0\|^2_{L^2} + CA^{-\f{3}{4}}t\Big)   \nonumber\\
& \le  e^{C { A^{-\f{5}{8}}\mathrm{log}^{\f{1}{4}}A} }\Big( \|y(n_{in})_0\|^2_{L^2} + C { A^{-\f{5}{8}}\mathrm{log}^{\f{1}{4}}A }\Big) \nonumber\\
& \le C\big(\|y(n_{in})_0\|^2_{L^2} + { A^{-\f{3}{8}} }\big).
\end{align}
Integrating  inequality \eqref{41311} with respect to $t$ on  $[0, T]$ and using \eqref{424482}, we obtain
\begin{align*}
\|yn_0\|^2_{L^\infty L^2}+\f{1}{A}\|y\partial_yn_0\|^2_{L^2L^2}
& \le \|y(n_{in})_0\|^2_{L^2} + \f{CT}{A^\f{3}{4}}+\f{CT}{A^\f{3}{4}}\big(\|y(n_{in})_0\|^2_{L^2}+1\big) \\
& \le \|y(n_{in})_0\|^2_{L^2} +  { \f{C}{A^{\f{3}{8}}} } \big(\|y(n_{in})_0\|^2_{L^2}+1\big).
\end{align*}
This together with \eqref{41413} completes the proof of \eqref{4142}.

\smallskip 

{\bf Step 2.  Estimates of $\omega_0$.}   For the non-weighted part, we first multiply equation $\eqref{AH11}_3$ by $\omega_0$ and use the integration by parts to  obtain 
\begin{align*}
\frac{1}{2}\f{d}{dt}\|\omega_0\|^2_{L^2}+\f{1}{A}\|\partial_y \omega_0\|^2_{L^2}
& =\f{1}{A}\int_{\mathbb T\times\mathbb R}\partial_y\omega_0(u^2_{\neq}\omega_{\neq})_0 dxdy\\
&\le \f{1}{2A}\|\partial_y\omega_0\|^2_{L^2} + \f{1}{2A}\big\|(u^2_{\neq}\omega_{\neq})_0\big\|^2_{L^2} \\
&\leq\f{1}{2A}\|\partial_y\omega_0\|^2_{L^2}+\f{C}{A}\|\omega_{\neq}\|^2_{L^2}\|u^2_{\neq}\|^2_{L^\infty}\\
&\leq\f{1}{2A}\|\partial_y\omega_0\|^2_{L^2}+\f{C}{A}\|\omega_{\neq}\|^3_{L^2}\|\nabla\omega_{\neq}\|_{L^2}.
\end{align*}
Here, in the last inequality, we used the fact
\begin{equation}\label{43453}
\|{\bf u_{\neq}}\|_{L^\infty}=\|\nabla^\bot\Delta^{-1}\omega_{\neq}\|_{L^\infty}\leq C \|\nabla\nabla^\bot\Delta^{-1}\omega_{\neq}\|_{L^2}^{\frac{1}{2}}\|\nabla\partial_x\nabla^\bot\Delta^{-1} \omega_{\neq}\|_{L^2}^{\frac{1}{2}}\leq C \|\omega_{\neq}\|_{L^2}^{\frac{1}{2}}\|\nabla\omega_{\neq}\|_{L^2}^{\frac{1}{2}}
\end{equation}
due to Lemma \ref{LH32}. Thus, we obtain
\begin{equation*}
\f{d}{dt}\|\omega_0\|^2_{L^2}+\f{1}{A}\|\partial_y \omega_0\|^2_{L^2}
\leq \f{C}{A}\|\omega_{\neq}\|^3_{L^2}\|\nabla\omega_{\neq}\|_{L^2}.
\end{equation*}
Integrating this inequality with respect to $t$ and using the assumptions (A-5) and (A-6), we have
\begin{align}\label{414452}
\|\omega_0\|^2_{L^{\infty}L^2}+\f{1}{A}\|\partial_y \omega_0\|^2_{L^2L^2} 
& \le \|(\omega_{in})_0\|^2_{L^2}+\f{C}{A}\int_0 ^T\|\omega_{\neq}\|^3_{ L^2}\|\nabla\omega_{\neq}\|_{L^2}d\tau\nonumber\\
& \le  \|(\omega_{in})_0\|^2_{L^2}+\f{C}{A}\|\omega_{\neq}\|^3_{L^\infty L^2}\Big(\int_0^T\|\nabla\omega_{\neq}\|^2_{L^2}d\tau\Big)^{\f{1}{2}}T^{\f{1}{2}} \nonumber\\
& \le  \|(\omega_{in})_0\|^2_{L^2}+\f{CC_0^3}{A}\Big(\|(\omega_{in})_{\neq}\|^3_X+{ A^{-\f{9}{4}} }\Big)\Big(\|(\omega_{in})_{\neq}\|_X+{ A^{-\f{3}{8}} }
\Big)A^{\f{1}{2}}T^{\f{1}{2}}\nonumber\\
& \le \|(\omega_{in})_0\|^2_{L^2}+\f{CC_0^3}{A^{\f{1}{2}}} { A^{-\f{9}{4}} }\big(A^{-\f{3}{4}}+{ A^{-\f{3}{8}} }\big)
{ A^{\f{1}{16}}\mathrm{log}^{\f{1}{8}}A }\nonumber\\
& \le \|(\omega_{in})_0\|^2_{L^2}+\f{C}{A^{\f{3}{2}}} 
\end{align}
due to { $0<T\leq\lambda_A^{-\f{1}{4}}=A^{\f{1}{8}}\mathrm{log}^{\f{1}{4}}A$ }and $A^{\f{3}{4}}\|\omega_{in}\|_X\leq 1$ for all $A>A^{\circ}_1$ fulfilling { ${A^{\circ}_1}^{\f{23}{16}} \ge C^3_0$.  }

We next to establish the weighted estimates of $\omega_0$.  Similar to the above procedure, we multiply equation $\eqref{AH11}_3$ by $y^2\omega_0$ to deduce that 
\begin{align*}
&\frac{1}{2}\f{d}{dt}\|y\omega_0\|^2_{L^2}+\f{1}{A}\|y\partial_y \omega_0\|^2_{L^2} \\
& = -\f{2}{A}\int_{\mathbb T\times\mathbb R}y\omega_0\partial_y\omega_0dxdy+\f{1}{A}\int_{\mathbb T\times\mathbb R}\partial_y(y^2\omega_0)({u}^2_{\neq}\omega_{\neq})_0dxdy \\
& = \f{1}{A}\|\omega_0\|^2_{L^2}+\f{1}{A}\int_{\mathbb T\times\mathbb R}y^2\partial_y\omega_0({u}^2_{\neq}\omega_{\neq})_0dxdy
+\f{2}{A}\int_{\mathbb T\times\mathbb R}y\omega_0({u}^2_{\neq}\omega_{\neq})_0dxdy \\
& \le \f{1}{A}\|\omega_0\|^2_{L^2}+\f{1}{2A}\|y\partial_y\omega_0\|^2_{L^2}+\f{1}{2A}\|y(u^2_{\neq}\omega_{\neq})_0\|^2_{L^2} + \f{C}{A} \|y\omega_{\neq}\|_{L^2}\|u^2_{\neq}\|_{L^\infty}\|\omega_0\|_{L^2} \\
& \le \f{2}{A}\|\omega_0\|^2_{L^2}+\f{1}{2A}\|y\partial_y\omega_0\|^2_{L^2}+\f{C}{A}\|y\omega_{\neq}\|^2_{L^2}\|u^2_{\neq}\|^2_{L^\infty}, 
\end{align*}
which together with \eqref{43453}  implies that
\begin{align*}
\f{d}{dt}\|y\omega_0\|^2_{L^2}+\f{1}{A}\|y\partial_y \omega_0\|^2_{L^2}
& \le  \f{C}{A}\|\omega_0\|^2_{L^2}+\f{C}{A}\|y\omega_{\neq}\|^2_{L^2}\|\omega_{\neq}\|_{L^2}\|\nabla\omega_{\neq}\|_{L^2} \\
& \le \f{C}{A}\|\omega_0\|^2_{L^2}+\f{C}{A}\|\omega_{\neq}\|^3_{X}\|\nabla\omega_{\neq}\|_{L^2}.
\end{align*}
Integrating the above inequality with respect to $t$ and combining \eqref{414452} with the assumptions  (A-5) and (A-6), we see that
\begin{align*}
\|y\omega_0\|^2_{L^{\infty}L^2} + \f{1}{A} \|y\partial_y \omega_0\|^2_{L^2L^2}  
& \le \|y(\omega_{in})_0\|^2_{L^2} + \f{C}{A} \|\omega_0\|^2_{L^2L^2} + \f{C}{A}\int_0^T \|\omega_{\neq}\|^3_{ X}\|\nabla\omega_{\neq}\|_{L^2}d\tau  \\
& \le  \|y(\omega_{in})_0\|^2_{L^2}+\f{C}{A}\|\omega_0\|^2_{L^\infty L^2}T + \f{C}{A}\|\omega_{\neq}\|^3_{L^\infty X}\Big(\int_0^T\|\nabla\omega_{\neq}\|^2_{L^2}d\tau\Big)^{\f{1}{2}}T^{\f{1}{2}}    \\
& \le \|y(\omega_{in})_0\|^2_{L^2}+\f{C}{A}\big(\|(\omega_{in})_0\|^2_{L^2}+A^{-\f{3}{2}}\big){ A^{\f{1}{8}}\mathrm{log}^{\f{1}{4}}A } \\
& \qquad + \f{CC_0^3}{A}\big(\|(\omega_{in})_{\neq}\|_X^3+{ A^{-\f{9}{4}} }\big) \big(\|(\omega_{in})_{\neq}\|_X+ { A^{-\f{3}{8}}\big) A^{\f{1}{2}}A^{\f{1}{16}}\mathrm{log}^{\f{1}{8}}A   }\\
& \le \|y(\omega_{in})_0\|^2_{L^2} + \f{C}{A} A^{-\f{3}{2}} { A^{\f18}\mathrm{log}^{\f{1}{4}}A + \f{CC_0^3}{A} A^{-\f{9}{4}} A^{-\f{3}{8}} A^{\f{1}{2}}A^{\f{1}{16}} \mathrm{log}^{\f{1}{8}}A }\\
& \le \|y(\omega_{in})_0\|^2_{L^2} +\f{C}{A^{\f{3}{2}}}
\end{align*} 
provided that  $A>A^{\circ}_1$.  Then combining this inequality with  \eqref{414452}, we obtained \eqref{414}.

\smallskip 

{\bf Step 3.  Estimates of $u_0^1$.}   We take the $L^2$ inner product of equation \eqref{015}  with $u_0^1$ and use the integration by parts to obtain
\begin{align*}
\frac{1}{2}\frac{d}{dt}\|u_0^1\|^2_{L^2}+\f{1}{A}\|\partial_yu_0^1\|^2_{L^2}
& = \f{1}{A}\int_{\mathbb T\times\mathbb R}(u_{\neq}^1u_{\neq}^2)_0\partial_yu_0^1dxdy \\
& \le  \frac{1}{2A}\|\partial_yu_0^1\|^2_{L^2}+\f{1}{2A}\|(u_{\neq}^1u_{\neq}^2)_0\|^2_{L^2} \\
& \le \frac{1}{2A}\|\partial_yu_0^1\|^2_{L^2}+\f{C}{A}\|u^1_{\neq}\|_{L^\infty}^2\|u^2_{\neq}\|_{L^2}^2, 
\end{align*}
which together with \eqref{43453} and \eqref{4121}  implies that 
\begin{equation}\label{411}
\frac{d}{dt}\|u_0^1\|^2_{L^2}+\f{1}{A}\|\partial_yu_0^1\|^2_{L^2}
\leq\f{C}{A}\|\omega_{\neq}\|^3_{L^2}\|\nabla\omega_{\neq}\|_{L^2}.
\end{equation}
 Integrating inequality \eqref{411} with respect to $t$ on $[0, T]$ and using assumptions (A-5) and (A-6), we deduce that
\begin{align}\label{412}
\|u_0^1\|^2_{L^{\infty}L^2}
& \le \|(u_{in}^1)_0\|^2_{L^2} + \f{C}{A}\|\omega_{\neq}\|^3_{L^\infty L^2}\Big(\int_0^T\|\nabla\omega_{\neq}\|_{L^2}^2d\tau\Big)^{\frac{1}{2}}T^{\frac{1}{2}}      \nonumber\\
& \le \|(u_{in}^1)_0\|^2_{L^2}+\f{CC_0^3}{A}\big(\|(\omega_{in})_{\neq}\|^3_X+ { A^{-\f{9}{4}}\big)\big(\|(\omega_{in})_{\neq}\|_X+A^{-\f{3}{8}}
\big)A^{\f{1}{2}}A^{\f{1}{16}}\mathrm{log}^{\f{1}{8}}A }    \nonumber\\
& \le  \|(u_{in}^1)_0\|^2_{L^2} + { \f{CC_0^3}{A^{\f12}}A^{-\f{9}{4}}\big(A^{-\f{3}{4}}+A^{-\f{3}{8}}\big) A^{\f{1}{16}} \mathrm{log}^{\f{1}{8}}A    }       \nonumber\\
& \le \|{\bf u}_{in}\|^2_{L^2} + \f{C}{A^{\f{3}{2}}} 
\end{align}
for all $A>A^{\circ}_1$.  Taking advantage of the Gagliardo-Nirenberg inequality and of the fact $\partial_yu_0^1=-\omega_0,$ we also obtain
\begin{align*} 
\|u_0^1\|^2_{L^\infty L^\infty}
& \le C\|u_0^1\|_{L^\infty L^2}\|\partial_yu_0^1\|_{L^\infty L^2} = C\|u_0^1\|_{L^\infty L^2}\|\omega_0\|_{L^\infty L^2}\nonumber\\
& \le C\Big(\|{\bf u}_{in}\|_{L^2} +  \f{C}{A^{\f{3}{4}}} \Big)\Big(\|(\omega_{in})_0\|_{L^2} +  \f{C}{A^{\f{3}{4}}} \Big)\nonumber\\
& \le C\Big(\|{\bf u}_{in}\|^2_{L^2} + \|(\omega_{in})_0\|^2_{L^2} +  \f{1}{A^{\f{3}{2}}} \Big)
\end{align*}
due to \eqref{412} and \eqref{414452}. Thus we obtain the inequality \eqref{413} for all  $A>A^{\circ}_1$. 

Summarily, taking $A_1=\max\big\{A^{*}_1, A^{**}_1, A^{\circ}_1\big\}$,  we complete the proof of Lemma \ref{LH41}.           \qquad $\Box$

\subsection{Nonzero modes estimates}

\quad      We now  establish the refined bounds (B-1)-(B-6) one by one. 

\begin{Lemma} \label{LH42}
Under the assumptions {\rm(A-2)} and {\rm (A-4)-(A-6)}, there exists a positive constant $A_2$ relying on $C_0$,  $C_\infty$,  $\|n_{in}\|_X$ and $\|{\bf u}_{in}\|_{L^2}$ such that if $A>A_2$ and $A^{\f{3}{4}}\|\omega_{in}\|_X\leq1$, then it holds 
\begin{equation}\label{boota1}
\f{1}{A}\int_0^t\|\nabla n_{\neq}(\cdot,\tau)\|_X^2d\tau\leq 2\|( n_{in})_{\neq}\|_X^2  \qquad {\rm  for \,\, all \,\,\,}  t\in[0, T]. 
\end{equation}
\end{Lemma}
{\bf Proof.} Multiplying equation \eqref{qh}$_1$ by $n_{\neq}$,  using the integration by parts over $\mathbb T\times \mathbb R$ and noticing that $$\int_{\mathbb T\times\mathbb R}y^2\partial_xn_{\neq}n_{\neq}dxdy=\f{1}{2}\int_{\mathbb R}y^2\Big(\int_{\mathbb T}\partial_xn^2_{\neq}dx\Big)dy=0,$$ we have
\begin{align}\label{42426}
&\f{1}{2}\f{d}{dt}\|n_{\neq}\|^2_{L^2}+\f{1}{A}\|\nabla n_{\neq}\|^2_{L^2}     \nonumber\\
& = \f{1}{A}\int_{\mathbb T\times\mathbb R}\Big((n_{\neq}\nabla c_{\neq})_{\neq}+(n_{0}\nabla c_{\neq})\Big)\cdot\nabla n_{\neq}dxdy
+\f{1}{A}\int_{\mathbb T\times\mathbb R}(n_{\neq}\partial_yc_{0})\partial_yn_{\neq}dxdy\nonumber\\
& \qquad +\f{1}{A}\int_{\mathbb T\times\mathbb R}\Big(({\bf u}_{\neq} n_{\neq})_{\neq}+({\bf u}_{\neq} n_{0})\Big)\cdot\nabla n_{\neq}dxdy
+ \f{1}{A}\int_{\mathbb T\times\mathbb R}({\bf u}_{0} n_{\neq})\cdot\nabla n_{\neq}dxdy.
\end{align}
We now estimate the right hand side of \eqref{42426} one by one. Indeed, by using the Young inequality and the H{\"o}lder inequality, we obtain
\begin{align*}
\int_{\mathbb T\times\mathbb R}\Big((n_{\neq}\nabla c_{\neq})_{\neq}+(n_{0}\nabla c_{\neq})\Big)\cdot\nabla n_{\neq}dxdy
\leq&\|\nabla n_{\neq}\|_{L^2}\big(\|(n_{\neq}\nabla c_{\neq})_{\neq}\|_{L^2}+\|n_{0}\nabla c_{\neq}\|_{L^2}\big)\nonumber\\
\leq&\f{1}{6}\|\nabla n_{\neq}\|^2_{L^2}+C\big(\|n_{\neq}\nabla c_{\neq}\|^2_{L^2}+\|n_{0}\nabla c_{\neq}\|^2_{L^2}\big)\nonumber\\
\leq&\f{1}{6}\|\nabla n_{\neq}\|^2_{L^2}+C\big(\|n_{\neq}\|^2_{L^\infty}+\|n_{0}\|^2_{L^\infty}\big)\|\nabla c_{\neq}\|^2_{L^2}
\end{align*}
and 
\begin{align*}
\int_{\mathbb T\times \mathbb R}n_{\neq}\partial_yc_{0}\partial_yn_{\neq}dxdy\leq\f{1}{12}\|\partial_y n_{\neq}\|^2_{L^2}+C\|n_{\neq}\partial_yc_{0}\|^2_{L^2}\leq\f{1}{12}\|\nabla n_{\neq}\|^2_{L^2}+C\|n_{\neq}\|^2_{L^\infty}\|\partial_yc_{0}\|^2_{L^2}.
\end{align*}
For the terms related to the velocity, we can take a similar procedure and use \eqref{4121} to obtain
\begin{align*}
\int_{\mathbb T\times \mathbb R}\Big(({\bf u}_{\neq}n_{\neq})_{\neq}+{\bf u}_{\neq}n_{0}\Big)\cdot\nabla n_{\neq}dxdy
\leq&\|\nabla n_{\neq}\|_{L^2}\big(\|({\bf u}_{\neq}n_{\neq})_{\neq}\|_{L^2}+\|{\bf u}_{\neq}n_{0}\|_{L^2}\big)\\
\leq&\f{1}{6}\|\nabla n_{\neq}\|^2_{L^2}+C\big(\|{\bf u}_{\neq}n_{\neq}\|^2_{L^2}+\|{\bf u}_{\neq}n_{0}\|^2_{L^2}\big)\\
\leq&\f{1}{6}\|\nabla n_{\neq}\|^2_{L^2}+C\big(\|n_{\neq}\|^2_{L^\infty}+\|n_0\|^2_{L^\infty}\big)\|{\bf u}_{\neq}\|^2_{L^2}\\
\leq&\f{1}{6}\|\nabla n_{\neq}\|^2_{L^2}+C\big(\|n_{\neq}\|^2_{L^\infty}+\|n_0\|^2_{L^\infty}\big)\|\omega_{\neq}\|^2_{L^2}
\end{align*}
and
\begin{align*}
\int_{\mathbb T\times \mathbb R}({\bf u}_{0}n_{\neq})\cdot\nabla n_{\neq}dxdy=&\int_{\mathbb T\times \mathbb R}u^1_{0}n_{\neq}  \partial_xn_{\neq}dxdy\\
\leq&\f{1}{12}\|\partial_xn_{\neq}\|^2_{L^2}+C\| u_0^1n_{\neq}\|^2_{L^2}\\
\leq&\f{1}{12}\|\nabla n_{\neq}\|^2_{L^2}+C\|u^1_0\|^2_{L^\infty}\|n_{\neq}\|^2_{L^2}.
\end{align*}
Substituting the above estimates into \eqref{42426} and then using Lemma \ref{LH35} and Lemma \ref{LH36}, we obtain
\begin{align*}
&\f{d}{dt}\|n_{\neq}\|^2_{L^2}+\f{1}{A}\|\nabla n_{\neq}\|^2_{L^2}\\
& \le \f{C}{A}\big(\|n_{\neq}\|^2_{L^\infty}+\|n_{0}\|^2_{L^\infty}\big)\big(\|\nabla c_{\neq}\|^2_{L^2}+\|\omega_{\neq}\|^2_{L^2}\big)
+\f{C}{A}\|n_{\neq}\|^2_{L^\infty}\|\partial_y c_0\|^2_{L^2}+\f{C}{A}\|u^1_0\|^2_{L^\infty} \|n_{\neq}\|^2_{L^2} \\
& \le \f{C}{A} \|n\|^2_{L^\infty} \big(\|n_{\neq}\|^2_{X} + \|\omega_{\neq}\|^2_{L^2}\big)
+ \f{C}{A}\|n\|^2_{L^\infty}\| n_0\|^2_{X} + \f{C}{A}\|u^1_0\|^2_{L^\infty} \|n_{\neq}\|^2_{L^2} 
\end{align*}
for all $t\in[0, T]$.  Then a direct integration together with \eqref{4142} and \eqref{413} yields that
\begin{align*}
&\|n_{\neq}(t)\|^2_{L^2}+\f{1}{A}\int_{0}^{t}\|\nabla n_{\neq}\|^2_{L^2}d\tau \\
& \le \|( n_{in})_{\neq}\|^2_{L^2} + \f{C}{A}\int_0^t \|n\|^2_{L^\infty} \big(\|n_{\neq}\|^2_{X} + \|\omega_{\neq}\|^2_{L^2} + \|n_0\|^2_{X} \big)d\tau  +\f{C}{A}\int_0^t\|n_{\neq}\|^2_{L^2}\| u^1_{0}\|^2_{L^\infty}d\tau  \\
& \le \|( n_{in})_{\neq}\|^2_{L^2}+\f{C}{A}\|n\|^2_{L^\infty L^\infty}\big(\|n_{\neq}\|^2_{L^\infty X} + \|\omega_{\neq}\|^2_{L^\infty L^2} + \|n_0\|^2_{L^\infty X} \big)T  + \f{C}{A}\|n_{\neq}\|^2_{L^\infty L^2}\|u_0^1\|^2_{L^\infty L^\infty}T  \\
& \le \|( n_{in})_{\neq}\|^2_{L^2}+\f{C}{A}\|n\|^2_{L^\infty L^\infty}\big(\|n_{\neq}\|^2_{L^\infty X} + \|\omega_{\neq}\|^2_{L^\infty L^2} + \|( n_{in})_0\|^2_X +1 \big)T  \\
& \qquad + \f{C}{A}\|n_{\neq}\|^2_{L^\infty L^2}\big(\|{\bf u}_{in}\|^2_{L^2} + \|(\omega_{in})_0\|^2_{L^2} +1 \big)T.  
\end{align*}
According  to the assumptions (A-2), (A-4),  (A-6), { $0 < T \le \lambda_A^{-\f{1}{4}} = A^{\f{1}{8}}\mathrm{log}^{\f{1}{4}}A$ } and $\|\omega_{in}\|_X\le A^{-\f{3}{4}} < 1$, we can further deduce that 
\begin{align*}
&\|n_{\neq}(t)\|^2_{L^2}+\f{1}{A}\int_{0}^{t}\|\nabla n_{\neq}\|^2_{L^2}d\tau  \\
& \le \|( n_{in})_{\neq}\|^2_{L^2}+\f{CC^2_0C^2_\infty}{A}\Big(\|( n_{in})_{\neq}\|^2_X + \|( \omega_{in})_{\neq}\|^2_X   + \|( n_{in})_0\|^2_X + 1\Big) { A^{\f{1}{8}}\mathrm{log}^{\f{1}{4}}A } \\ 
& \qquad + \f{CC^2_0}{A}\|( n_{in})_{\neq}\|^2_{X}\Big(\|{\bf u}_{in}\|^2_{L^2}+\|( \omega_{in})_{0}\|^2_{X}+1\Big) { A^{\f{1}{8}}\mathrm{log}^{\f{1}{4}}A } \\
& \le \|( n_{in})_{\neq}\|^2_{L^2} + { \f{CC_0^2C^2_\infty}{A^{\f{5}{8}}} }\Big(\|( n_{in})_{\neq}\|^2_X\big(\|{\bf u}_{in}\|^2_{L^2} + 1 \big) + \|( n_{in})_0\|^2_X + 1 \Big).
\end{align*}
Therefore, choosing $A^{*}_2$ large enough such that 
\[
{ {A^{*}_2}^{\f{5}{8}} }\geq \f{CC_0^2C^2_\infty\Big(\|( n_{in})_{\neq}\|^2_{X}\big(\|{\bf u}_{in}\|^2_{L^2}+1\big)+\|( n_{in})_{0}\|^2_{X}+1\Big)}{\|(n_{in})_{\neq}\|^2_{L^2}}+1,
\]
 we can conclude that
\begin{equation}\label{42421}
\f{1}{A}\int_0^t \|\nabla n_{\neq}(\cdot,\tau)\|^2_{L^2}d\tau\leq2\|( n_{in})_{\neq}\|^2_{L^2},
\end{equation}
for all $A>A^{*}_2$. 

To get the weighted estimate of $n_{\neq}$ in \eqref{boota1}, we take a similar procedure by multiplying equation ${\eqref{qh}}_1$ by $y^2n_{\neq}$ and  using the integration by parts and the fact 
\[
\int_{\mathbb T\times\mathbb R}y^4\partial_xn_{\neq}n_{\neq}dxdy=\f{1}{2}\int_{\mathbb R}y^4\Big(\int_{\mathbb T}\partial_xn^2_{\neq}dx\Big)dy=0
\]
to obtain
\begin{align}\label{42438}
&\f{1}{2}\f{d}{dt}\|yn_{\neq}\|^2_{L^2}+\f{1}{A}\|y\nabla n_{\neq}\|^2_{L^2}-\f{1}{A}\|n_{\neq}\|^2_{L^2}\nonumber\\
& = \f{1}{A}\int_{\mathbb T\times\mathbb R}\Big((n_{\neq}\nabla c_{\neq})_{\neq}+n_{0}\nabla c_{\neq}\Big)\cdot \nabla(y^2n_{\neq})dxdy
+\f{1}{A}\int_{\mathbb T\times\mathbb R}(n_{\neq}\partial_yc_{0})\partial_y(y^2n_{\neq})dxdy\nonumber\\
& \qquad + \f{1}{A}\int_{\mathbb T\times\mathbb R}\Big(({\bf u}_{\neq} n_{\neq})_{\neq}+{\bf u}_{\neq} n_{0}\Big)\cdot\nabla(y^2n_{\neq})dxdy
+\f{1}{A}\int_{\mathbb T\times\mathbb R}({\bf u}_{0} n_{\neq})\cdot\nabla(y^2n_{\neq})dxdy.
\end{align}
We estimate the right hand side of \eqref{42438} one by one. Indeed, for the first two terms, we have
\begin{align*}
&\int_{\mathbb T\times\mathbb R}\Big((n_{\neq}\nabla c_{\neq})_{\neq}+n_0\nabla c_{\neq}\Big)\cdot \nabla(y^2n_{\neq})dxdy \\
& = \int_{\mathbb T\times \mathbb R}\Big((n_{\neq}\nabla c_{\neq})_{\neq} + n_0\nabla c_{\neq}\Big) \cdot y^2\nabla n_{\neq}dxdy+2\int_{\mathbb T\times \mathbb R}\Big((n_{\neq}\partial_y c_{\neq})_{\neq}+n_0\partial_y c_{\neq}\Big)yn_{\neq}dxdy  \\
& \le \|y\nabla n_{\neq}\|_{L^2}\big(\|y(n_{\neq}\nabla c_{\neq})_{\neq}\|_{L^2}+\|yn_0\nabla c_{\neq}\|_{L^2}\big)
+2\big(\|(n_{\neq}\partial_yc_{\neq})_{\neq}\|_{L^2}+\|n_0\partial_yc_{\neq}\|_{L^2}\big) \|yn_{\neq}\|_{L^2}  \\
& \le \f{1}{6}\|y\nabla n_{\neq}\|^2_{L^2}+C\big(\|n_{\neq}\|^2_{L^\infty}+\|n_0\|^2_{L^\infty}\big)\|y\nabla c_{\neq}\|^2_{L^2}+C\big(\|n_{\neq}\|_{L^\infty}+\|n_0\|_{L^\infty}\big)\|\partial_yc_{\neq}\|_{L^2}\|yn_{\neq}\|_{L^2} \\
& \le \f{1}{6}\|y\nabla n_{\neq}\|^2_{L^2}+C\|n\|^2_{L^\infty}\|y\nabla c_{\neq}\|^2_{L^2}+C\|n\|_{L^\infty}\|\partial_yc_{\neq}\|_{L^2}\|n_{\neq}\|_{X}
\end{align*}
and
\begin{align*}
\int_{\mathbb T\times\mathbb R}(n_{\neq}\partial_yc_0)\partial_y(y^2n_{\neq})dxdy 
& =  \int_{\mathbb T\times \mathbb R}n_{\neq}\partial_yc_0 y^2\partial_yn_{\neq}dxdy+2\int_{\mathbb T\times \mathbb R}n_{\neq}\partial_y c_0 y n_{\neq}dxdy\\
& \le \f{1}{12}\|y\partial_y n_{\neq}\|^2_{L^2} + 3\|yn_{\neq}\|^2_{L^2}\|\partial_yc_{0}\|^2_{L^\infty} + 2\|n_{\neq}\|_{L^2}\|\partial_yc_{0}\|_{L^\infty}\|yn_{\neq}\|_{L^2}\\
& \le \f{1}{12}\|y\nabla n_{\neq}\|^2_{L^2}+C\|n_{\neq}\|^2_{X}\big(\|\partial_yc_{0}\|^2_{L^\infty}+\|\partial_yc_{0}\|_{L^\infty}\big).
\end{align*}
Similarly,  we also have
\begin{align*}
&\int_{\mathbb T\times\mathbb R}\Big(({\bf u}_{\neq} n_{\neq})_{\neq}+{\bf u}_{\neq} n_0\Big)\cdot\nabla(y^2n_{\neq})dxdy       \\
& = \int_{\mathbb T\times \mathbb R}\Big(({\bf u}_{\neq}n_{\neq})_{\neq}+{\bf u}_{\neq}n_0\Big)\cdot y^2\nabla n_{\neq}dxdy+2\int_{\mathbb T\times \mathbb R}\Big(( u^2_{\neq}n_{\neq})_{\neq}+ u^2_{\neq}n_0\Big)yn_{\neq}dxdy        \\
& \le \|y\nabla n_{\neq}\|_{L^2}\big(\|y({\bf u}_{\neq}n_{\neq})_{\neq}\|_{L^2}+\|y{\bf u}_{\neq}n_0\|_{L^2}\big)+2\|yn_{\neq}\|_{L^2}\big(\|u^2_{\neq}n_{\neq}\|_{L^2}+\|u^2_{\neq}n_0\|_{L^2}\big)   \\
& \le \f{1}{6}\|y\nabla n_{\neq}\|^2_{L^2} + C\|y{\bf u}_{\neq}\|^2_{L^2}\big(\|n_{\neq}\|^2_{L^\infty}+\|n_0\|^2_{L^\infty}\big)+C\|u^2_{\neq}\|_{L^2}\big(\|n_{\neq}\|_{L^\infty}
+\|n_0\|_{L^\infty}\big)\|yn_{\neq}\|_{L^2}\\
& \le \f{1}{6}\|y\nabla n_{\neq}\|^2_{L^2}+C\|\omega_{\neq}\|^2_{X}\|n\|^2_{L^\infty}+C\|\omega_{\neq}
\|_{L^2}\|n\|_{L^\infty}\|n_{\neq}\|_{X}
\end{align*}
thanks to \eqref{cm341},  and
\begin{align*}
\int_{\mathbb T\times\mathbb R}({\bf u}_{0} n_{\neq})\cdot\nabla(y^2n_{\neq})dxdy
& = \int_{\mathbb T\times \mathbb R}u_0^1n_{\neq}y^2\partial_xn_{\neq}dxdy\\
& \le \f{1}{12}\|y\partial_xn_{\neq}\|^2_{L^2}+C\|u_0^1yn_{\neq} \|^2_{L^2}\\
& \le \f{1}{12}\|y\nabla n_{\neq}\|^2_{L^2}+C\|u^1_0\|^2_{L^\infty}\|n_{\neq}\|^2_{X}.
\end{align*}
Substituting the above estimates into \eqref{42438} and using Lemma \ref{LH35} and Lemma \ref{LH36}, we obtain 
\begin{align*} 
& \f{d}{dt}\|yn_{\neq}\|^2_{L^2}+\f{1}{A}\|y\nabla n_{\neq}\|^2_{L^2}    \\ 
& \le \f{C}{A}\|n\|^2_{L^\infty}\big(\|y\nabla c_{\neq}\|^2_{L^2}+\|\omega_{\neq}\|^2_{X}\big)  + \f{C}{A}\|n\|_{L^\infty}\big(\|\partial_yc_{\neq}\|_{L^2}+\|\omega_{\neq}\|_{L^2}\big)
\|n_{\neq}\|_X\\
&\qquad +\f{C}{A}\big(\|\partial_yc_0\|^2_{L^\infty} + \|\partial_yc_0\|_{L^\infty} + \|u^1_0\|^2_{L^\infty}\big)\|n_{\neq}\|^2_X \\
& \le \f{C}{A}\|n\|^2_{L^\infty}\big(\|n_{\neq}\|^2_{X}+\|\omega_{\neq}\|^2_{X}\big)  + \f{C}{A}\|n\|_{L^\infty}\big(\|n_{\neq}\|_{X}+\|\omega_{\neq}\|_{L^2}\big)
\|n_{\neq}\|_X\\
&\qquad +\f{C}{A}\big(\|n_0\|^2_{L^2} + \|n_0\|^2_{L^\infty} + \|n_0\|_{L^2} + \|n_0\|_{L^\infty} +  \|u^1_0\|^2_{L^\infty}\big) \|n_{\neq}\|^2_X  \\
& \le \f{C}{A}\|n\|^2_{L^\infty}\big(\|n_{\neq}\|^2_{X}+\|\omega_{\neq}\|^2_{X}\big)  
  + \f{C}{A}\big(\|n_0\|^2_{L^2} + \|n_0\|^2_{L^\infty}  +  \|u^1_0\|^2_{L^\infty} + 1\big) \|n_{\neq}\|^2_X.
\end{align*}
Then taking an integration with respect to $t$, we see that 
\begin{align*}
 \|yn_{\neq}(t)\|^2_{L^2}+\f{1}{A}\int_{0}^{t}\|y\nabla n_{\neq}\|^2_{L^2}d\tau  
& \le \|y( n_{in})_{\neq}\|^2_{L^2}+\f{C}{A}\|n\|^2_{L^\infty L^\infty}\big(\|n_{\neq}\|^2_{L^\infty X}+\|\omega_{\neq}\|^2_{L^\infty X}\big)T   \\
& \qquad  + \f{C}{A}\big(\|n_{0}\|^2_{L^2}+\|n_{0}\|^2_{L^\infty}+\|u_0^1\|^2_{L^\infty L^\infty}+1\big)\|n_{\neq}\|^2_{L^\infty X}T
\end{align*}
for all $t\in[0, T]$. It follows from the assumptions  (A-2), (A-4), (A-6),  { $0 < T \le \lambda_A^{-\f{1}{4}} = A^{\f{1}{8}}\mathrm{log}^{\f{1}{4}}A$ }, and the estimates \eqref{4142} and \eqref{413}, we can further deduce that 
\begin{align*}
&\|yn_{\neq}(t)\|^2_{L^2}+\f{1}{A}\int_{0}^{t}\|y\nabla n_{\neq}\|^2_{L^2}d\tau \\
& \le \|y( n_{in})_{\neq}\|^2_{L^2}+\f{CC_0^2C^2_\infty}{A}\Big(\|( n_{in})_{\neq}\|^2_{X}+\|( \omega_{in})_{\neq}\|^2_{X}+{ A^{-\f{3}{2}}\Big)A^{\f{1}{8}}\mathrm{log}^{\f{1}{4}}A   }  \\
& \qquad +\f{CC_0^2C^2_\infty}{A}\Big(\|( n_{in})_{0}\|^2_X + \|{\bf u}_{in}\|^2_{L^2} + \|( \omega_{in})_0\|^2_X + 1\Big)\|( n_{in})_{\neq}\|^2_X  { A^{\f{1}{8}}\mathrm{log}^{\f{1}{4}}A } \\
&  \le \|y( n_{in})_{\neq}\|^2_{L^2}+{ \f{CC_0^2C^2_\infty}{A^{\f{5}{8}}} }\Big(\big(\|( n_{in})_{0}\|^2_{X}+\|{\bf u}_{in}\|^2_{L^2}+\|( \omega_{in})_{0}\|^2_{X} + 1\big)\|( n_{in})_{\neq}\|^2_{X}+\|( \omega_{in})_{\neq}\|^2_{X} + 1\Big)  \\
& \le \|y( n_{in})_{\neq}\|^2_{L^2}
+{ \f{CC_0^2C^2_\infty}{A^{\f{5}{8}}} }\Big(\big(\|( n_{in})_{0}\|^2_{X}+\|{\bf u}_{in}\|^2_{L^2}
 + 1 \big)\|( n_{in})_{\neq}\|^2_X + 1\Big).
\end{align*}
By choosing $A^{\circ}_2$ large enough such that 
\[
{ {A^{\circ}_2}^{\f{5}{8}} }\geq\f{CC_0^2C^2_\infty\Big(\big(\|( n_{in})_0\|^2_{X}+\|{\bf u}_{in}\|^2_{L^2}
+ 1 \big)\|( n_{in})_{\neq}\|^2_X + 1\Big)}{\|y(n_{in})_{\neq}\|^2_{L^2}},
\]
 we obtain
\begin{equation}\label{42448}
\f{1}{A}\int_{0}^{t}\|y\nabla n_{\neq}(\cdot,\tau)\|^2_{L^2}d\tau\leq2\|y( n_{in})_{\neq}\|^2_{L^2}
\end{equation}
for all $A>A^{\circ}_2$.   Taking $A_2:=\max \big\{A^{*}_2, A^{\circ}_2\big\}$ and combining \eqref{42421} with \eqref{42448}, we conclude that
\begin{equation}\nonumber
\f{1}{A}\int_{0}^{t}\|\nabla n_{\neq}(\cdot,\tau)\|^2_{X}d\tau\leq2\|( n_{in})_{\neq}\|^2_{X}.
\end{equation}
This completes the proof of Lemma \ref{LH42}. \qquad $\Box$

\begin{Lemma}\label{LH43}
Under the assumptions {\rm(A-1)-(A-6)},  there exists a positive constant $A_3$ relying on $C_0$,  $C_\infty$,  $\|n_{in}\|_X$, $\|\partial_xn_{in}\|_X$ and $\|{\bf u}_{in}\|_{L^2}$ such that if $A>A_3$ and $A^{\f{3}{4}}\|\omega_{in}\|_X\leq1$, then it holds 
\[
\| n_{\neq}(t)\|_X\leq 2C_0e^{-\epsilon_0\lambda_A t}\| ( n_{in})_{\neq}\|_X \qquad {\rm  for \,\, all \,\,\,}  t\in[0, T]. 
\]
\end{Lemma}
{\begin{Remark}
In \cite{ZZZ}, Zhang et al. obtained the corresponding enhanced dissipation estimates through a key lemma related to the Couette flow, which is invalid for the Poiseuille flow. To overcome this difficulty, there are two effective methods. The first one is using integration by parts directly, while the second one is utilizing the semigroup estimate of the Poiseuille flow established by Coti Zelati et al. \cite{ZEW}. Although one-order derivative can be reduced by using the former, the power of $\f1A$ is smaller when dealing with the enhanced dissipation estimates of the weighted part. However, the latter does not reduce the order of the derivative, but the power of $\f1A$ is larger for the weighted enhanced dissipation estimates. Here, we utilize the second method since $\|D^2 c\|_X$ can be controlled easily due to $\eqref{3}_2$.  
\end{Remark} }
{\bf Proof.} Recalling ${\mathcal L}= \f{1}{A}\Delta - y^2\partial_x$, we can deduce from  the Duhamel formula
\begin{align*}
n_{\neq}=&e^{{\mathcal L}t}( n_{in})_{\neq}-\f{1}{A}e^{{\mathcal L}t}\int_0^te^{-{\mathcal L}\tau}\Big(\nabla\cdot(n_{\neq}\nabla c_{\neq})_{\neq}+\nabla\cdot(n_0\nabla c_{\neq})+\partial_y(n_{\neq}\partial_y c_0)\Big)d\tau\\
&-\f{1}{A}e^{{\mathcal L}t}\int_0^te^{-{\mathcal L}\tau}\Big(\nabla\cdot(n_{\neq}{\bf u_{\neq}})_{\neq}+\nabla\cdot(n_{\neq}{\bf u_0})+\nabla\cdot(n_0 {\bf u_{\neq}})\Big)d\tau
\end{align*}
and the fact  $\int_{\mathbb T}f_{\neq}(x,y)dx=0$ for each $y\in{\mathbb R}$ as well as  Remark \ref{LH31remark}  that
\begin{align}\label{43450}
\|n_{\neq}(t)\|_{X}
& \le \sqrt{2}\big\|e^{{\mathcal L}t}( n_{in})_{\neq}\big\|_{X} + \f{2}{A}\Big\|e^{{\mathcal L}t}\int_0^te^{-{\mathcal L}\tau}\big(\nabla\cdot(n_{\neq}\nabla c_{\neq})_{\neq}+\nabla\cdot(n_0\nabla c_{\neq})+\partial_y(n_{\neq}\partial_y c_0)\big)d\tau\Big\|_{X}\nonumber\\
& \qquad +\f{2}{A}\Big\|e^{{\mathcal L}t}\int_0^te^{-{\mathcal L}\tau}\big(\nabla\cdot(n_{\neq}{\bf u_{\neq}})_{\neq}+\nabla\cdot(n_{\neq}{\bf u_0})+\nabla\cdot(n_0 {\bf u_{\neq}})\big)d\tau\Big\|_{X}\nonumber\\
& \le  C_0 e^{-\epsilon_0\lambda_A t}\Big(\sqrt{2}\|( n_{in})_{\neq}\|_{X} + \f{2}{A}\int_0^te^{\epsilon_0\lambda_A \tau}\big\|\nabla\cdot(n_{\neq}\nabla c_{\neq})_{\neq}+\nabla\cdot(n_0\nabla c_{\neq})+\partial_y(n_{\neq}\partial_y c_0)\big\|_{X}d\tau  \nonumber\\
& \qquad +\f{2}{A} \int_0^te^{\epsilon_0\lambda_A \tau}\big\|\nabla\cdot(n_{\neq}{\bf u_{\neq}})_{\neq}+\nabla\cdot(n_{\neq}{\bf u_0})+\nabla\cdot(n_0 {\bf u_{\neq}})\big\|_{X}d\tau\Big).
\end{align}
We now estimate the integrals on the right hand side of \eqref{43450}. For the first one,  we use the fact 
\[
\big\| \nabla \cdot f_{\neq}\big\|_{X}
=\big\| \nabla \cdot f - \nabla \cdot f_0\big\|_{X}
= \big\| \nabla \cdot f - (\nabla \cdot f)_0\big\|_{X}
\le \big\| \nabla \cdot f \big\|_{X} + \big\| (\nabla \cdot f)_0\big\|_{X}
\le C \big\| \nabla \cdot f \big\|_{X} 
\]
and the H\"{o}lder inequality to obtain
\begin{align*}
&\int_0^te^{\epsilon_0\lambda_A \tau}\big\|\nabla\cdot(n_{\neq}\nabla c_{\neq})_{\neq}+\nabla\cdot(n_0\nabla c_{\neq})+\partial_y(n_{\neq}\partial_y c_0)\big\|_{X}d\tau \\
& \le C\int_0^t \Big(\|\nabla n_{\neq}\cdot\nabla c_{\neq}\|_{X} + \|n_{\neq}\Delta c_{\neq}\|_{X} + \|\partial_yn_0\partial_yc_{\neq}\|_{X} + \|n_0\Delta c_{\neq}\|_{X} \\
& \qquad  + \|\partial_yn_{\neq}\partial_y c_0\|_{X} + \|n_{\neq}\partial_{yy} c_{0}\|_{X} \Big)d\tau   \\
& \le C \big( \|\nabla c_{\neq}\|_{L^\infty L^\infty} + \|\partial_yc_0\|_{L^\infty L^\infty} \big) \int_0^t \big( \|\nabla n_{\neq}\|_{X}  + \|\partial_y n_{0}\|_{X} \big) d\tau \\
& \qquad + C \big(\|\Delta c_{\neq}\|_{L^\infty L^\infty} + \|\partial_{yy} c_0\|_{L^\infty L^\infty} \big)\int_0^t \big( \|n_{\neq}\|_{X} + \|n_0\|_{X}\big) d\tau\\
& \le C \big( \|\nabla c_{\neq}\|_{L^\infty L^\infty} + \|\partial_yc_0\|_{L^\infty L^\infty} \big)  \Big(\int_0^t \big( \|\nabla n_{\neq}\|^2_{X} + \|\partial_y n_0\|^2_{X} \big)d\tau\Big)^{\f{1}{2}} t^{\f{1}{2}} \\
& \qquad + C \big(\|\Delta c_{\neq}\|_{L^\infty L^\infty} + \|\partial_{yy} c_0\|_{L^\infty L^\infty} \big) \big(\|n_{\neq}\|_{L^\infty X} + (\|n_0\|_{L^\infty X} \big) t. 
\end{align*}
Here in the first inequality,  we used the fact that { $\lambda_A\tau\leq\lambda_AT\leq\lambda_A^{\f{3}{4}}\leq1$ }and thus that $0\leq e^{\epsilon_0\lambda_A \tau}\leq C$ for all $\tau\in[0, t]$.  
For the terms related to the velocity, we take a similar procedure and deduce from \eqref{43453} and  the H\"{o}lder inequality  that 
\begin{align*}
&\int_0^te^{\epsilon_0\lambda_A \tau}\big\|\nabla\cdot(n_{\neq}{\bf u_{\neq}})_{\neq}+\nabla\cdot(n_{\neq}{\bf u_0})+\nabla\cdot(n_0 {\bf u_{\neq}})\big\|_{X} d\tau \\
& \le C \int_0^t  \Big( \|\nabla n_{\neq}\cdot{\bf u_{\neq}}\|_{X} + \|\partial_xn_{\neq}u_0^1\|_{X} +  \|\partial_yn_0u^2_{\neq}\|_{X}  \Big) d\tau\\
& \le  C\int_0^t  \big(\|\nabla n_{\neq}\|_{X} + \|\partial_y n_0\|_{X}\big)\big( \|{\bf u}_{\neq}\|_{L^\infty} + \|u_0^1\|_{L^\infty} \big) d\tau        \\
& \le  C\int_0^t  \big(\|\nabla n_{\neq}\|_{X} + \|\partial_y n_0\|_{X}\big) \big(\|\omega_{\neq}\|_{L^2} + \|\nabla \omega_{\neq}\|_{L^2} + \|u_0^1\|_{L^\infty} \big)d\tau  \\
& \le C \big(\|\omega_{\neq}\|_{L^\infty L^2}+ \|u_0^1\|_{L^\infty L^\infty}\big)\int_0^t  \big(\|\nabla n_{\neq}\|_{X} + \|\partial_y n_0\|_{X}\big) d\tau \\
& \qquad +  C \int_0^t  \big(\|\nabla n_{\neq}\|_{X} + \|\partial_y n_0\|_{X}\big) \|\nabla \omega_{\neq}\|_{L^2}d\tau\\
& \le C\big(\|\omega_{\neq}\|_{L^\infty L^2}+ \|u_0^1\|_{L^\infty L^\infty}\big) \Big(\int_0^t  \big(\|\nabla n_{\neq}\|^2_{X} + \|\partial_y n_0\|^2_{X}\big) d\tau\Big)^{\f{1}{2}} t^{\f{1}{2}}      \\
& \qquad + C\Big(\int_0^t \big(\|\nabla n_{\neq}\|^2_{X} +  \|\partial_y n_0\|^2_{X}\big) d\tau\Big)^{\f{1}{2}}\Big(\int_0^t\|\nabla \omega_{\neq}\|^2_{L^2}d\tau\Big)^{\f{1}{2}}. 
\end{align*}

Substituting the above estimates into \eqref{43450}, we can deduce from Lemma \ref{LH35} and \ref{LH36} that
\begin{align*}
& \|n_{\neq}(t)\|_{X}  \\
& \le C_0e^{-\epsilon_0\lambda_A t}\Bigg(\sqrt{2}\big\|( n_{in})_{\neq}\|_{X}    \\
& \qquad + \f{C}{A} \big( \|\nabla c_{\neq}\|_{L^\infty L^\infty} + \|\partial_yc_0\|_{L^\infty L^\infty} + \|\omega_{\neq}\|_{L^\infty L^2}+ \|u_0^1\|_{L^\infty L^\infty} \big)  \Big(\int_0^t \big( \|\nabla n_{\neq}\|^2_{X} + \|\partial_y n_0\|^2_{X} \big)d\tau\Big)^{\f{1}{2}} T^{\f{1}{2}}  \\
& \qquad +  \f{C}{A} \big(\|\Delta c_{\neq}\|_{L^\infty L^\infty} + \|\partial_{yy} c_0\|_{L^\infty L^\infty} \big) \big(\|n_{\neq}\|_{L^\infty X} + (\|n_0\|_{L^\infty X} \big) T  \\
& \qquad + \f{C}{A} \Big(\int_0^t \big(\|\nabla n_{\neq}\|^2_{X} +  \|\partial_y n_0\|^2_{X}\big) d\tau\Big)^{\f{1}{2}}\Big(\int_0^t\|\nabla \omega_{\neq}\|^2_{L^2}d\tau\Big)^{\f{1}{2}} \Bigg)  \\
& \le C_0e^{-\epsilon_0\lambda_A t}\Bigg(\sqrt{2}\big\|( n_{in})_{\neq}\|_{X}   \\
& \qquad + \f{C}{A} \big( \|n_{\neq}\|_{L^\infty L^2}+\|n_{\neq}\|_{L^\infty L^\infty}+\|\partial_xn_{\neq}\|_{L^\infty L^2}+ \|n_0\|_{L^\infty L^2} + \|n_0\|_{L^\infty L^\infty}  + \|\omega_{\neq}\|_{L^\infty L^2}+ \|u_0^1\|_{L^\infty L^\infty} \big)   \\
& \qquad \cdot \Big(\int_0^t \big( \|\nabla n_{\neq}\|^2_{X} + \|\partial_y n_0\|^2_{X} \big)d\tau\Big)^{\f{1}{2}} T^{\f{1}{2}} \\
& \qquad +  \f{C}{A} \big( \|n_{\neq}\|_{L^\infty L^2} + \|n_{\neq}\|_{L^\infty L^\infty} + \|\partial_xn_{\neq}\|_{L^\infty L^2}+ \|n_0\|_{L^\infty L^2} + \|n_0\|_{L^\infty L^\infty} \big) \big(\|n_{\neq}\|_{L^\infty X} + (\|n_0\|_{L^\infty X} \big) T  \\
& \qquad + \f{C}{A} \Big(\int_0^t \big(\|\nabla n_{\neq}\|^2_{X} +  \|\partial_y n_0\|^2_{X}\big) d\tau\Big)^{\f{1}{2}}\Big(\int_0^t\|\nabla \omega_{\neq}\|^2_{L^2}d\tau\Big)^{\f{1}{2}} \Bigg). 
\end{align*}
Then taking advantage of the assumptions (A-1)-(A-6), { $0 < T \le \lambda_A^{-\f{1}{4}} = A^{\f{1}{8}}\mathrm{log}^{\f{1}{4}}A$ },  $\|\omega_{in}\|_X \le A^{-\f{3}{4}}$ and Lemma \ref{LH41}, we conclude that 
\begin{align*}
& \|n_{\neq}(t)\|_{X}  \\
& \le C_0e^{-\epsilon_0\lambda_A t}\Big(\sqrt{2}\big\|( n_{in})_{\neq}\|_{X}   \\
& \qquad + \f{CC_0C_\infty}{A} \big( \| n_{in}\|_{X} +\|( \partial_xn_{in})_{\neq}\|_{X}  + \|{\bf u}_{in}\|_{L^2}+\|\omega_{in}\|_{X} + 1 \big)  
\Big(A \| n_{in}\|_X^2   + 1\Big)^{\f{1}{2}}  { A^{\f{1}{16}}\mathrm{log}^{\f{1}{8}}A }\\
& \qquad +  \f{CC_0^2C_\infty}{A} \big( \| n_{in}\|_{X} + \|( \partial_xn_{in})_{\neq}\|_{X}   + 1  \big) \big(\|n_{in})\|_{X}+1\big) { A^{\f{1}{8}}\mathrm{log}^{\f{1}{4}}A  }\\
& \qquad + \f{C}{A} \Big(A \| n_{in} \|_X^2 + 1\Big)^{\f{1}{2}}\Big(A\big(\| (\omega_{in})_{\neq}\|_X^2+{ A^{-\frac{3}{4}}\big) \Big)^{\f{1}{2}} \Big) }\\
& \le C_0e^{-\epsilon_0\lambda_A t}\Big(\sqrt{2}\big\|( n_{in})_{\neq}\|_{X}  + { \f{CC_0C_\infty}{A^{\f{5}{16}}} }\big(\| n_{in}\|_{X} +\|( \partial_xn_{in})_{\neq}\|_{X}  + \|{\bf u}_{in}\|_{L^2}+\|\omega_{in}\|_{X} + 1 \big) \big(\| n_{in}\|_X^2   + 1\big)^{\f{1}{2}}   \\
& \qquad +  { \f{CC_0^2C_\infty}{A^{\f58}} }\big( \| n_{in}\|_{X} + \|( \partial_xn_{in})_{\neq}\|_{X}   + 1  \big) \big(\|n_{in})\|_{X}+1\big) + { \f{C}{A^{\f{3}{8}}} }\big(\| n_{in} \|_X^2 + 1\big)^{\f{1}{2}}\Big) \\
& \le C_0e^{-\epsilon_0\lambda_A t}\Big(\sqrt{2}\big\|( n_{in})_{\neq}\|_{X}  + { \f{CC_0^2C_\infty}{A^{\f{5}{16}}}  }\big(\| n_{in}\|_{X} +\|( \partial_xn_{in})_{\neq}\|_{X}  + \|{\bf u}_{in}\|_{L^2} + 1 \big) \big(\| n_{in}\|_X   + 1\big) \Big)
\end{align*}

Summarily,  we can choose $A_3$  fulfilling 
\begin{align*}
{ A_3^{\f{5}{16}}  }
& \ge \f{CC_0^2C_\infty \big(\| n_{in}\|_{X} +\|( \partial_xn_{in})_{\neq}\|_{X}  + \|{\bf u}_{in}\|_{L^2} + 1 \big) \big(\| n_{in}\|_X   + 1\big)}{(2-\sqrt{2})\|( n_{in})_{\neq}\|^2_{X}}
\end{align*}
to ensure that 
\[
\|n_{\neq}(t)\|_{X}\leq 2C_0e^{-\epsilon_0\lambda_A t}\|( n_{in})_{\neq}\|_{X} 
\]
for all $t\in[0, T]$ whenever $A>A_3$.  This completes the proof of Lemma \ref{LH43}. \qquad $\Box$

\begin{Lemma}\label{LH44}
Under the assumptions {\rm(A-2)-(A-6)},  there exists a positive constant $A_4$ relying on  $C_0$,  $C_\infty$,   $\|n_{in}\|_X$,  $\|(\partial_xn_{in})_{\neq}\|_X$ and $\|{\bf u}_{in}\|_{L^2}$ such that if  $A>A_4$ and $A^{\f{3}{4}}\|\omega_{in}\|_X\leq1$, then it holds 
\[
\| \partial_x n_{\neq}\|^2_{L^{\infty}X}\leq 2\| (\partial_x n_{in})_{\neq}\|^2_X.
\]
\end{Lemma}
{\bf Proof.} To estimate the non-weighted part,  applying $\partial_x$ to equation ${\eqref{qh}}_{1}$ and taking the $L^2$ inner product with $\partial_xn_{\neq}$, we obtain 
\begin{align}\label{4683}
& \f{1}{2}\f{d}{dt} \|\partial_xn_{\neq}\|^2_{L^2}+\f{1}{A}\|\nabla\partial_x n_{\neq}\|^2_{L^2}\nonumber\\
& = \f{1}{A} \int_{\mathbb T\times\mathbb R} \Big( \big( \partial_x(n_{\neq}\nabla c_{\neq})_{\neq} + \partial_x(n_{0}\nabla c_{\neq})\big)\cdot\nabla\partial_xn_{\neq} + \partial_{x}(n_{\neq}\partial_yc_{0})\partial_{xy}n_{\neq}\Big)dxdy \nonumber \\
& \qquad + \f{1}{A} \int_{\mathbb T\times\mathbb R} \Big( \partial_x({\bf u}_{\neq} n_{\neq})_{\neq} + \partial_x({\bf u}_{\neq} n_{0}) + \partial_x({\bf u}_{0} n_{\neq}) \Big)\cdot\nabla\partial_xn_{\neq}dxdy.
\end{align}
For the first integral  on the right hand side of \eqref{4683}, we deduce from $\partial_xf_{\neq}=\partial_xf$ that
\begin{align*}
& \int_{\mathbb T\times\mathbb R} \Big( \big( \partial_x(n_{\neq}\nabla c_{\neq})_{\neq} + \partial_x(n_{0}\nabla c_{\neq}) \big)\cdot\nabla\partial_xn_{\neq} + \partial_{x}(n_{\neq}\partial_yc_{0})\partial_{xy}n_{\neq}\Big)dxdy \\
& = \int_{\mathbb T\times \mathbb R} \Big((\partial_xn_{\neq}\nabla c_{\neq}+n_{\neq}\nabla \partial_xc_{\neq} + n_0\nabla \partial_xc_{\neq} ) \cdot\nabla \partial_xn_{\neq} + \partial_x n_{\neq}\partial_yc_0\partial_{xy}n_{\neq} \Big) dxdy\\
& \le \f{1}{4} \|\nabla \partial_xn_{\neq}\|^2_{L^2} + 3\Big(\|\partial_xn_{\neq}\|^2_{L^2}\|\nabla c_{\neq}\|^2_{L^\infty}  + \|n_{\neq}\|^2_{L^\infty}\|\nabla \partial_xc_{\neq}\|^2_{L^2}  + \|n_0\|^2_{L^\infty}\|\nabla \partial_xc_{\neq}\|^2_{L^2}  \\
& \qquad   + \| \partial_xn_{\neq}\|^2_{L^2}\|\partial_yc_{0}\|^2_{L^\infty}\Big)      \\
& \le \f{1}{4} \|\nabla \partial_xn_{\neq}\|^2_{L^2} + C\Big(\|\partial_xn_{\neq}\|^2_{L^2} \big( \|\nabla c_{\neq}\|^2_{L^\infty} + \|\partial_yc_{0}\|^2_{L^\infty} \big) + \|n\|^2_{L^\infty}\|\nabla \partial_xc_{\neq}\|^2_{L^2} \Big).      
\end{align*}
Similarly, for the terms related to the velocity, we can deduce that
\begin{align*}
& \int_{\mathbb T\times\mathbb R} \Big( \partial_x({\bf u}_{\neq} n_{\neq})_{\neq} + \partial_x({\bf u}_{\neq} n_0) + \partial_x({\bf u}_{0} n_{\neq}) \Big)\cdot\nabla\partial_xn_{\neq}dxdy     \\
& = \int_{\mathbb T\times \mathbb R} \Big( \big(\partial_x{\bf u}_{\neq}n_{\neq} + {\bf u}_{\neq}\partial_xn_{\neq} + \partial_x{\bf u}_{\neq}n_{0}\big)\cdot\nabla \partial_xn_{\neq} + u_0^1\partial_xn_{\neq}\partial_{xx}n_{\neq}  \Big) dxdy\\
& \le \f{1}{4}\|\nabla \partial_xn_{\neq}\|^2_{L^2}  + 3 \Big( \|\partial_x{\bf u}_{\neq}\|^2_{L^2}\|n_{\neq}\|^2_{L^\infty} + \|{\bf u}_{\neq}\|^2_{L^\infty}\|\partial_xn_{\neq}\|^2_{L^2} + \|\partial_x{\bf u}_{\neq}\|^2_{L^2}\|n_{0}\|^2_{L^\infty}  + \|u^1_0\|^2_{L^\infty}\|\partial_xn_{\neq}\|^2_{L^2} \Big)\\
& \le \f{1}{4}\|\nabla \partial_xn_{\neq}\|^2_{L^2}  + C \Big( \|\partial_x{\bf u}_{\neq}\|^2_{L^2}\|n\|^2_{L^\infty} + \big( \|{\bf u}_{\neq}\|^2_{L^\infty}  + \|u^1_0\|^2_{L^\infty} \big) \|\partial_xn_{\neq}\|^2_{L^2} \Big). 
\end{align*}
Substituting the above estimates into \eqref{4683} and using \eqref{43453} and Lemma \ref{LH36}, we see that 
\begin{align*}
&\f{d}{dt}\|\partial_xn_{\neq}\|^2_{L^2}+\f{1}{A}\|\nabla\partial_xn_{\neq}\|^2_{L^2}\\
& \le \f{C}{A}\Big(\|\nabla c_{\neq}\|^2_{L^\infty}+\|\partial_yc_0\|^2_{L^\infty} + \|{\bf u}_{\neq}\|^2_{L^\infty}  + \|u^1_0\|^2_{L^\infty} \Big)\|\partial_xn_{\neq}\|^2_{L^2}   +  \f{C}{A} \|n\|^2_{L^\infty}  \big( \|\nabla \partial_xc_{\neq}\|^2_{L^2}  +  \|\partial_x{\bf u}_{\neq}\|^2_{L^2} \big)\\
& \le \f{C}{A}\Big(\|\nabla c_{\neq}\|^2_{L^\infty}+\|\partial_yc_0\|^2_{L^\infty} 
+\|\omega_{\neq}\|^2_{L^2}+\|\nabla\omega_{\neq}\|^2_{L^2} +\|u_0^1\|^2_{L^\infty} \Big)\|\partial_xn_{\neq}\|^2_{L^2}  \\
& \qquad +  \f{C}{A} \|n\|^2_{L^\infty} \big(\|n_{\neq}\|^2_{X}+\|\omega_{\neq}\|^2_{L^2}\big). 
\end{align*}
Then we integrate the above inequality with respect to $t$  and use Lemma \ref{LH35} and Lemma \ref{LH36} to obtain
\begin{align*}
& \|\partial_xn_{\neq}(t)\|^2_{L^2} + \f{1}{A}\int_{0}^{t}\|\nabla \partial_xn_{\neq}\|^2_{L^2}d\tau\\
& \le \|(\partial_x n_{in})_{\neq}\|^2_{L^2}  +\f{C}{A}\int_0^t\big(\|\nabla c_{\neq}\|^2_{L^\infty}+\|\partial_yc_0\|^2_{L^\infty} + \|\omega_{\neq}\|^2_{L^2}  + \|u_0^1\|^2_{L^\infty} \big) \|\partial_xn_{\neq}\|^2_{L^2}d\tau    \\
& \qquad + \f{C}{A}\int_0^t \|n\|^2_{L^\infty} \big(\|n_{\neq}\|^2_X + \|\omega_{\neq}\|^2_{L^2}\big)d\tau  +\f{C}{A}\int_0^t \|\nabla\omega_{\neq}\|^2_{L^2} \|\partial_xn_{\neq}\|^2_{L^2}d\tau   \\
& \le \|(\partial_x n_{in})_{\neq}\|^2_{L^2} +\f{C}{A}\big(\|\nabla c_{\neq}\|^2_{L^\infty L^\infty}+\|\partial_yc_0\|^2_{L^\infty L^\infty} + \|\omega_{\neq}\|^2_{L^\infty L^2} + \|u_0^1\|^2_{L^\infty L^\infty} \big)\|\partial_xn_{\neq}\|^2_{L^\infty L^2}T\\
& \qquad +\f{C}{A} \|n\|^2_{L^\infty L^\infty} \big(\|n_{\neq}\|^2_{L^\infty X} + \|\omega_{\neq}\|^2_{L^\infty L^2}\big)T + \f{C}{A}\|\partial_xn_{\neq}\|^2_{L^\infty L^2}\int_0^t\|\nabla\omega_{\neq}\|^2_{L^2}d\tau\\
& \le \|(\partial_x n_{in})_{\neq}\|^2_{L^2}  + \f{C}{A}\big(\|n_{\neq}\|^2_{L^\infty L^2} + \|\partial_xn_{\neq}\|^2_{L^\infty L^2} + \|n_0\|^2_{L^\infty L^2} + \|n\|^2_{L^\infty L^\infty} + \|\omega_{\neq}\|^2_{L^\infty L^2}  \\
& \qquad + \|u_0^1\|^2_{L^\infty L^\infty} \big)\|\partial_xn_{\neq}\|^2_{L^\infty L^2}T 
 + \f{C}{A} \|n\|^2_{L^\infty L^\infty} \big(\|n_{\neq}\|^2_{L^\infty X} + \|\omega_{\neq}\|^2_{L^\infty L^2}\big)T \\
 & \qquad  + \f{C}{A}\|\partial_xn_{\neq}\|^2_{L^\infty L^2}\int_0^t\|\nabla\omega_{\neq}\|^2_{L^2}d\tau 
\end{align*}
for all $t\in[0, T]$,  which together with the assumptions (A-2)-(A-6),  { $0 < T \le \lambda_A^{-\f{1}{4}} = A^{\f{1}{8}}\mathrm{log}^{\f{1}{4}}A$ },  $\|\omega_{in}\|_X \le A^{-\f{3}{4}}$  and Lemma \ref{LH41}  yields that
\begin{align*}
&\|\partial_xn_{\neq}(t)\|^2_{L^2}+\f{1}{A}\int_{0}^{t}\|\nabla \partial_xn_{\neq}\|^2_{L^2}d\tau\nonumber\\
& \le \|(\partial_x n_{in})_{\neq}\|^2_{L^2}  + \f{CC_0^2C^2_\infty}{A}\big(\| n_{in}\|^2_{X} +\|(\partial_xn_{in})_{\neq}\|^2_{X}  + \|\omega_{in}\|^2_{X} + \|{\bf u}_{in}\|^2_{L^2} + 1\big)\|(\partial_xn_{in})_{\neq}\|^2_{X} { A^{\f{1}{8}}\mathrm{log}^{\f{1}{4}}A } \\
& \qquad   + \f{CC_0^2C^2_\infty}{A}\big(\|(n_{in})_{\neq}\|^2_{X}+
\|(\omega_{in})_{\neq}\|^2_{X}+ { A^{-\f{3}{2}}\big) A^{\f{1}{8}}\mathrm{log}^{\f{1}{4}}A }
  + \f{C}{A}\|(\partial_xn_{in})_{\neq}\|^2_{X}\big(\|(\omega_{in})_{\neq}\|^2_{X}+ { A^{-\f{3}{4}}\big)A } \\
& \le \|(\partial_x n_{in})_{\neq}\|^2_{L^2}  + { \f{CC_0^2C^2_\infty}{A^{\f58}} }\big( \|n_{in}\|^2_{X}+\|(\partial_xn_{in})_{\neq}\|^2_{X} + \|{\bf u}_{in}\|^2_{L^2} + 1\big)\|(\partial_xn_{in})_{\neq}\|^2_{X} \nonumber\\
& \qquad + { \f{CC_0^2C^2_\infty}{A^{\f58}} 
 }\big(\|n_{in}\|^2_{X} + 1\big) + C\|(\partial_xn_{in})_{\neq}\|^2_{X}\big(A^{-\f{3}{2}}+ { A^{-\f{3}{4}}\big) }\nonumber\\
& \le \|(\partial_x n_{in})_{\neq}\|^2_{L^2}
+{ \f{CC_0^2C^2_\infty}{A^{\f{5}{8}}} }\big( \|n_{in}\|^2_{X} +\|(\partial_xn_{in})_{\neq}\|^2_{X} 
+\|{\bf u}_{in}\|^2_{L^2} + 1\big)\big(\|(\partial_xn_{in})_{\neq}\|^2_{X}+1\big).
\end{align*}
By taking $A^{*}_4$ sufficiently large such that
\[
{ {A^{*}_4}^{\f{5}{8}} }
\ge \f{CC_0^2C^2_\infty\big( \|n_{in}\|^2_{X} +\|(\partial_xn_{in})_{\neq}\|^2_{X} 
+\|{\bf u}_{in}\|^2_{L^2} + 1\big)\big(\|(\partial_xn_{in})_{\neq}\|^2_{X}+1\big)}{\|(\partial_x n_{in})_{\neq}\|^2_{L^2}},
\]
we conclude that 
\begin{equation}\label{44438}
\|\partial_xn_{\neq}(t)\|^2_{L^2}+\f{1}{A}\int_{0}^{t}\|\nabla \partial_xn_{\neq}\|^2_{L^2}d\tau\leq2\|(\partial_x n_{in})_{\neq}\|^2_{L^2}
\end{equation}
for all $t\in[0, T]$ whenever  $A>A^{*}_4$. 

We take a similar procedure to deal with the weighted part. Precisely,  we apply  $\partial_x$ to equation ${\eqref{qh}}_{1}$  and multiply the resulting equation by $y^2\partial_xn_{\neq}$ to deduce that
\begin{align}\label{480}
&{\f{1}{2}}{\f{d}{dt}}\|y\partial_xn_{\neq}\|^2_{L^2}+\f{1}{A}\|y\nabla\partial_x n_{\neq}\|^2_{L^2} - \f{1}{A} \|\partial_xn_{\neq}\|^2_{L^2}     \nonumber\\
& = \f{1}{A} \int_{\mathbb T\times\mathbb R} \Big(\big(\partial_x(n_{\neq}\nabla c_{\neq})_{\neq} 
+ \partial_x(n_{0}\nabla c_{\neq})\big)\cdot\nabla(y^2\partial_xn_{\neq}) 
+ \partial_{x}(n_{\neq}\partial_yc_{0})\partial_y(y^2\partial_xn_{\neq})\Big) dxdy   \nonumber\\
&\qquad  + \f{1}{A} \int_{\mathbb T\times\mathbb R} \Big(\partial_x({\bf u}_{\neq} n_{\neq})_{\neq}   + \partial_x({\bf u}_{\neq} n_{0}) + \partial_x({\bf u}_{0} n_{\neq})\Big) \cdot\nabla(y^2\partial_xn_{\neq})dxdy.
\end{align}
For the first integral on the right hand side of \eqref{480},  we have
\begin{align*}
& \int_{\mathbb T\times\mathbb R} \Big(\big(\partial_x(n_{\neq}\nabla c_{\neq})_{\neq} 
     + \partial_x(n_{0}\nabla c_{\neq})\big)\cdot\nabla(y^2\partial_xn_{\neq}) 
     + \partial_{x}(n_{\neq}\partial_yc_{0})\partial_y(y^2\partial_xn_{\neq})\Big) dxdy \\
& =\int_{\mathbb T\times \mathbb R} \Big(\big(\partial_x(n_{\neq}\nabla c_{\neq}) 
     + \partial_x(n_{0}\nabla c_{\neq})\big) \cdot y^2\nabla \partial_xn_{\neq}  
     + \partial_{x}(n_{\neq}\partial_yc_{0})  y^2\partial_{xy}n_{\neq})\Big)dxdy       \\
& \qquad + 2 \int_{\mathbb T\times \mathbb R} \Big(\partial_x(n_{\neq}\partial_y c_{\neq}) + \partial_x(n_{0}\partial_y c_{\neq}) +   \partial_x(n_{\neq}\partial_yc_{0})  \Big) y \partial_xn_{\neq} dxdy\\
& =\int_{\mathbb T\times \mathbb R} \Big(\big(\partial_x n_{\neq}\nabla c_{\neq} +   n_{\neq}\nabla \partial_x c_{\neq} + n_0\nabla \partial_x c_{\neq} \big) \cdot y^2\nabla \partial_xn_{\neq}  
       + \partial_{x}n_{\neq}\partial_yc_{0} y^2\partial_{xy}n_{\neq})\Big)dxdy       \\
& \qquad - 2 \int_{\mathbb T\times \mathbb R} \big(n_{\neq}\partial_y c_{\neq} + n_0\partial_y c_{\neq} + n_{\neq}\partial_yc_0 \big) y \partial_{xx}n_{\neq} dxdy  \\
& \le \f{1}{4}\|y\nabla \partial_xn_{\neq}\|^2_{L^2} + C \Big(\|y\partial_xn_{\neq}\nabla c_{\neq}\|^2_{L^2} + \|yn_{\neq}\nabla \partial_xc_{\neq}\|^2_{L^2} + \|yn_{0}\nabla\partial_x c_{\neq}\|^2_{L^2} + \|y\partial_x n_{\neq}\partial_yc_{0}\|^2_{L^2}  \\
& \qquad  + \|n_{\neq}\partial_yc_{\neq}\|^2_{L^2} + \|n_0\partial_yc_{\neq}\|^2_{L^2} +  \|n_{\neq}\partial_yc_{0} \|^2_{L^2} \Big)\\
& \le \f{1}{4}\|y\nabla\partial_x n_{\neq}\|^2_{L^2} + C \Big(\|y\partial_xn_{\neq}\|^2_{L^2}\|\nabla c_{\neq}\|^2_{L^\infty} + \|n\|^2_{L^\infty}\|y\nabla \partial_xc_{\neq}\|^2_{L^2}  + \| y\partial_xn_{\neq}\|^2_{L^2}\|\partial_yc_{0}\|^2_{L^\infty}  \\
& \qquad  +  \|n\|^2_{L^\infty}\| \partial_yc_{\neq}\|^2_{L^2}  + \|n_{\neq}\|^2_{L^2}\|\partial_yc_{0}\|^2_{L^\infty} \Big)\\
& \le \f{1}{4}\|y\nabla\partial_x n_{\neq}\|^2_{L^2} + C \Big(\|y\partial_xn_{\neq}\|^2_{L^2} \big(\|n_{\neq}\|^2_{L^2} + \|n_{\neq}\|^2_{L^\infty} + \|\partial_xn_{\neq}\|^2_{L^2}\big) + \|n\|^2_{L^\infty}\|n_{\neq}\|^2_{X}    \\
& \qquad  + \big(\| y\partial_xn_{\neq}\|^2_{L^2} + \|n_{\neq}\|^2_{L^2}\big) \big(\|n_0\|^2_{L^2}+\|n_0\|^2_{L^\infty}\big) \Big).
\end{align*}
Here we used  Lemma \ref{LH35} and  Lemma \ref{LH36} in the last inequality. Similarly, for the terms related to the velocity, we also have
\begin{align*}
&\int_{\mathbb T\times\mathbb R} \Big(\partial_x({\bf u}_{\neq} n_{\neq})_{\neq}  + \partial_x({\bf u}_{\neq} n_0) + \partial_x({\bf u}_{0} n_{\neq})\Big) \cdot\nabla(y^2\partial_xn_{\neq})dxdy \\
& = \int_{\mathbb T\times \mathbb R} \big(\partial_x({\bf u}_{\neq}n_{\neq}) + \partial_x{\bf u}_{\neq}n_0 + {\bf u}_0  \partial_x n_{\neq} \big) \cdot y^2\nabla \partial_xn_{\neq}dxdy \\
& \qquad  + 2 \int_{\mathbb T\times \mathbb R} \big(\partial_x(u^2_{\neq}n_{\neq}) + \partial_xu^2_{\neq}n_0\big)y\partial_xn_{\neq}dxdy\\
& = \int_{\mathbb T\times \mathbb R}\big(\partial_x{\bf u}_{\neq}n_{\neq} + {\bf u}_{\neq}\partial_xn_{\neq} + \partial_x{\bf u}_{\neq}n_0 + {\bf u}_0  \partial_x n_{\neq}  \big)\cdot y^2\nabla \partial_xn_{\neq}dxdy \\
& \qquad - 2\int_{\mathbb T\times \mathbb R}\big(u^2_{\neq}n_{\neq} + u^2_{\neq}n_0\big) y\partial_{xx}n_{\neq}dxdy\\
& \le \f{1}{4} \|y\nabla \partial_xn_{\neq}\|^2_{L^2} + C\Big(\|y\partial_x{\bf u}_{\neq}n_{\neq}\|^2_{L^2} + \|y{\bf u}_{\neq}\partial_xn_{\neq}\|^2_{L^2} + \|y\partial_x{\bf u}_{\neq}n_{0}\|^2_{L^2} 
+ \|yu_0^1\partial_xn_{\neq} \|^2_{L^2} \\
& \qquad +  \|u^2_{\neq}n_{\neq}\|^2_{L^2} + \| u^2_{\neq}n_0\|^2_{L^2} \Big) \\
& \le \f{1}{4} \|y\nabla \partial_xn_{\neq}\|^2_{L^2} + C\Big( \|y\partial_x{\bf u}_{\neq}\|^2_{L^2}\|n\|^2_{L^\infty}  + \|{\bf u}_{\neq}\|^2_{L^\infty}\|y\partial_xn_{\neq}\|^2_{L^2}   + \|u^1_0\|^2_{L^\infty}\|y\partial_xn_{\neq}\|^2_{L^2} +  \| u^2_{\neq}\|^2_{L^2}\|n\|^2_{L^\infty}\Big) \\
& \le \f{1}{4}\|y\nabla\partial_xn_{\neq}\|^2_{L^2} + C \Big(\|\omega_{\neq}\|^2_{X}\|n\|^2_{L^\infty}
 + \big(\|\omega_{\neq}\|^2_{L^2} + \|\nabla\omega_{\neq}\|^2_{L^2}  + \|u^1_0\|^2_{L^\infty}\big)\|y\partial_xn_{\neq}\|^2_{L^2} \Big)
\end{align*}
due to \eqref{4121}, \eqref{cm341} and \eqref{43453}. 

Substituting the above estimates into \eqref{480}, we obtain
\begin{align*}
&\f{d}{dt}\|y\partial_xn_{\neq}\|^2_{L^2}+\f{1}{A}\|y\nabla\partial_xn_{\neq}\|^2_{L^2}\\
& \le \f{C}{A} \Big(\big(\|n_{\neq}\|^2_{L^2}  + \|\partial_xn_{\neq}\|^2_{L^2} + \|n_0\|^2_{L^2} + \|n\|^2_{L^\infty} + \|\omega_{\neq}\|^2_{L^2}  + \|u^1_0\|^2_{L^\infty} \big)  \|y\partial_xn_{\neq}\|^2_{L^2}       \\
& \qquad + \|n\|^2_{L^\infty}\big(\|n_{\neq}\|^2_{X} + \|\omega_{\neq}\|^2_{X}  \big)  + \|n_{\neq}\|^2_{L^2} \|n_0\|^2_{L^2} + \|\partial_xn_{\neq}\|^2_{L^2} \Big)    + \f{C}{A}  \|\nabla\omega_{\neq}\|^2_{L^2}\|y\partial_xn_{\neq}\|^2_{L^2}. 
\end{align*}
We integrate the above inequality with respect to $t$ to deduce that 
\begin{align*}
& \|y\partial_xn_{\neq}(t)\|^2_{L^2}+\f{1}{A}\int_{0}^{t}\|y\nabla \partial_xn_{\neq}\|^2_{L^2}d\tau   \\
& \le  \|y(\partial_x n_{in})_{\neq}\|^2_{L^2} + \f{C}{A} \int_0^{t}\Big(\big(\|n_{\neq}\|^2_{L^2} + \|\partial_xn_{\neq}\|^2_{L^2} + \|n_0\|^2_{L^2} + \|n\|^2_{L^\infty}  + \|\omega_{\neq}\|^2_{L^2}  + \|u^1_0\|^2_{L^\infty} \big)  \|y\partial_xn_{\neq}\|^2_{L^2}       \\
& \qquad + \|n\|^2_{L^\infty}\big(\|n_{\neq}\|^2_{X} + \|\omega_{\neq}\|^2_{X}  \big)  + \|n_{\neq}\|^2_{L^2} \|n_0\|^2_{L^2} + \|\partial_xn_{\neq}\|^2_{L^2} \Big) d\tau   + \f{C}{A} \int_0^t \|\nabla\omega_{\neq}\|^2_{L^2}\|y\partial_xn_{\neq}\|^2_{L^2} d\tau \\
& \le \|y(\partial_x n_{in})_{\neq}\|^2_{L^2} + \f{C}{A}\Big(\|n_{\neq}\|^2_{L^\infty L^2} + \|\partial_xn_{\neq}\|^2_{L^\infty L^2} + \|n_0\|^2_{L^\infty L^2} +  \|n\|^2_{L^\infty L^\infty} 
+  \|\omega_{\neq}\|^2_{L^\infty L^2}  \\
& \qquad + \|u_0^1\|^2_{L^\infty L^\infty}\Big)\|y\partial_xn_{\neq}\|^2_{L^\infty L^2} T 
    + \f{C}{A} \|n\|^2_{L^\infty L^\infty} \big(\|n_{\neq}\|^2_{L^\infty X} + 
\|\omega_{\neq}\|^2_{L^\infty X}\big)T  \\
& \qquad  + \f{C}{A}\big(\|n_{\neq}\|^2_{L^\infty L^2} \|n_{0}\|^2_{L^\infty L^2} + \|\partial_xn_{\neq}\|^2_{L^\infty L^2} \big)T \\
& \qquad  +\f{C}{A}\|\omega_{\neq}\|^2_{L^\infty L^2}\|y\partial_xn_{\neq}\|^2_{L^\infty L^2}T 
+ \f{C}{A}\|y\partial_xn_{\neq}\|^2_{L^\infty L^2}\int_0^t\|\nabla\omega_{\neq}\|^2_{L^2}d\tau 
\end{align*}
for all $t\in[0, T]$.   It then folows the assumptions (A-2)-(A-6),  { $0 < T \le \lambda_A^{-\f{1}{4}} = A^{\f{1}{8}}\mathrm{log}^{\f{1}{4}}A$ },  $\|\omega_{in}\|_X \le A^{-\f{3}{4}}$  and Lemma \ref{LH41} that 
\begin{align*}
& \|y\partial_xn_{\neq}(t)\|^2_{L^2}+\f{1}{A}\int_{0}^{t}\|y\nabla \partial_xn_{\neq}\|^2_{L^2}d\tau  \\
& \le \|y(\partial_x n_{in})_{\neq}\|^2_{L^2}  + \f{CC_0^2C^2_\infty}{A}\big(\|( n_{in})_{\neq}\|^2_X + \|(\partial_x n_{in})_{\neq}\|^2_X + \|( n_{in})_0\|^2_X  + \|\omega_{in}\|^2_X +\|{\bf u}_{in}\|^2_{L^2} +  1\big)   \\
& \qquad \cdot \|(\partial_x n_{in})_{\neq}\|^2_X { A^{\f{1}{8}}\mathrm{log}^{\f{1}{4}}A } + \f{CC_0^2C^2_\infty}{A}\big(\|(n_{in})_{\neq}\|^2_{X}+\|(\omega_{in})_{\neq}\|^2_{X}+ { A^{-\f{3}{2}}\big)A^{\f{1}{8}}\mathrm{log}^{\f{1}{4}}A  } \\
& \qquad  + \f{CC_0^2}{A}\big(\|(n_{in})_{\neq}\|^2_X \big(\|(n_{in})_0\|^2_X + 1 \big) +  \|(\partial_x n_{in})_{\neq}\|^2_{X} \big) { A^{\f{1}{8}}\mathrm{log}^{\f{1}{4}}A  }\\
&\qquad +\f{CC_0^2}{A}\big(\|(\omega_{in})_{\neq}\|^2_{X}+ { A^{-\f{3}{2}} }\big)\|(\partial_x n_{in})_{\neq}\|^2_{X} { A^{\f{1}{8}}\mathrm{log}^{\f{1}{4}}A  }     +\f{C}{A}\|(\partial_x n_{in})_{\neq}\|^2_{X}\big(\|(\omega_{in})_{\neq}\|^2_{X}+ { A^{-\f{3}{4}}\big)A  }\\
& \le \|y(\partial_x n_{in})_{\neq}\|^2_{L^2}  + { \f{CC_0^2C^2_\infty}{A^{\f58}} }\big(\| n_{in}\|^4_X + \|(\partial_x n_{in})_{\neq}\|^2_X   +\|{\bf u}_{in}\|^2_{L^2} +  1\big) \big(\|(\partial_x n_{in})_{\neq}\|^2_X+1\big).
\end{align*}
Choosing $A^\circ_4$  fulfilling 
\[
{ {A^\circ_4}^{\f{5}{8}}  }
\ge 
\f{CC_0^2C^2_\infty\big(\| n_{in}\|^4_X + \|(\partial_x n_{in})_{\neq}\|^2_X   +\|{\bf u}_{in}\|^2_{L^2} +  1\big) \big(\|(\partial_x n_{in})_{\neq}\|^2_X+1\big)}{\|y(\partial_xn_{in})_{\neq}\|^2_{L^2}},
\]
we see that 
\begin{equation}\label{44492}
\|y\partial_xn_{\neq}(t)\|^2_{L^2}+\f{1}{A}\int_{0}^{t}\|y\nabla \partial_xn_{\neq}\|^2_{L^2}d\tau\leq2\|y(\partial_x n_{in})_{\neq}\|^2_{L^2}
\end{equation}
for all $t\in[0, T]$ whenever  $A>A^\circ_4$.  

Summarily, we can take $A_4:=\max\big\{A^{*}_4, A^\circ_4\big\}$ and then  see from \eqref{44438} and \eqref{44492} that
\begin{equation}\nonumber
\|\partial_xn_{\neq}(t)\|^2_{X}+\f{1}{A}\int_{0}^{t}\|\nabla \partial_xn_{\neq}(\cdot,\tau)\|^2_{X}d\tau\leq2\|(\partial_x n_{in})_{\neq}\|^2_{X}
\end{equation}
provided that $A>A_4$.  This completes the proof of Lemma \ref{LH44}. \qquad $\Box$

\begin{Lemma}\label{LH45}
Under the assumptions {\rm  (A-3)}, {\rm (A-5)} and {\rm (A-6)},  there exists a positive constant $A_5$ relying on $C_0$, $\|(\partial_xn_{in})_{\neq}\|_X$  and $\|{\bf u}_{in}\|_{L^2}$ such that if  $A>A_5$ and $A^{\f{3}{4}}\|\omega_{in}\|_X\leq1$, then it holds 
\[
\f{1}{A}\int_0^t\|\nabla \omega_{\neq}(\cdot,\tau)\|_X^2d\tau\leq 2\Big(\| (\omega_{in})_{\neq}\|_X^2+ { A^{-\frac{3}{4}}\Big) }   \qquad {\rm  for \,\, all \,\,\,}  t\in[0, T]. 
\]
\end{Lemma}
{\bf Proof.} For the non-weighted part, multiplying equation ${\eqref{qh}}_3$ by $\omega_{\neq},$ integrating the equation over $\mathbb T\times \mathbb R,$ and using the integration by parts, we obtain that 
\begin{align}\label{494}
& \f{1}{2}{\f{d}{dt}}\|\omega_{\neq}\|^2_{L^2}+\f{1}{A}\|\nabla \omega_{\neq}\|^2_{L^2} \nonumber\\
& = \f{1}{A} \int_{\mathbb T\times \mathbb R} \Big( ({\bf u}_{\neq} \omega_{\neq})_{\neq} + {\bf u}_0 \omega_{\neq} + {\bf u}_{\neq} \omega_0 \Big)\cdot\nabla\omega_{\neq}dxdy + \f{1}{A}\int_{\mathbb T\times\mathbb R}\partial_xn_{\neq}\omega_{\neq}dxdy.
\end{align}
Here, we also used the fact 
\[
2\int_{\mathbb T\times\mathbb R}u^2_{\neq}\omega_{\neq}dxdy 
=2\int_{\mathbb T\times\mathbb R}\partial_x\Delta^{-1}\omega_{\neq}\omega_{\neq}dxdy
= \int_{\mathbb T\times\mathbb R}\partial_x\big((-\Delta)^{-\f{1}{2}}\omega_{\neq}\big)^2dxdy=0.
\]
For the first integral on the  right hand side of \eqref{494},  we can use \eqref{43453} to obtain 
\begin{align*}
& \int_{\mathbb T\times \mathbb R} \Big( ({\bf u}_{\neq} \omega_{\neq})_{\neq} + {\bf u}_0 \omega_{\neq} + {\bf u}_{\neq} \omega_0 \Big)\cdot\nabla\omega_{\neq}dxdy \\
& \le \f{1}{2}\|\nabla \omega_{\neq}\|^2_{L^2} + C\Big( \|{\bf u}_{\neq}\omega_{\neq}\|^2_{L^2} + \|{\bf u}_0 \omega_{\neq}\|^2_{L^2}  + \|{\bf u}_{\neq} \omega_0\|^2_{L^2}\Big) \\
& \le \f{1}{2}\|\nabla \omega_{\neq}\|^2_{L^2} + C\Big( \|{\bf u}_{\neq}\|^2_{L^\infty} \big( \|\omega_{\neq}\|^2_{L^2}  + \|\omega_{0}\|^2_{L^2} \big)   + \|u^1_0\|^2_{L^\infty}\|\omega_{\neq}\|^2_{L^2}   \Big)      \\
& \le \f{1}{2}\|\nabla \omega_{\neq}\|^2_{L^2} + C\Big( \big(\|\omega_{\neq}\|^2_{L^2} + \|\nabla\omega_{\neq}\|^2_{L^2} \big)  \big(\|\omega_{\neq}\|^2_{L^2}   + \|\omega_{0}\|^2_{L^2} \big) + \|u^1_0\|^2_{L^\infty}\|\omega_{\neq}\|^2_{L^2}   \Big), 
\end{align*}
while for the second one, it is clear that 
\[
\int_{\mathbb T\times\mathbb R}\partial_xn_{\neq}\omega_{\neq}dxdy
\le \|\partial_xn_{\neq}\|_{L^2}\|\omega_{\neq}\|_{L^2}.
\]
Substituting the above inequalities into \eqref{494}, we deduce that
\begin{align*}
&{\f{d}{dt}}\|\omega_{\neq}(t)\|^2_{L^2}+\f{1}{A}\|\nabla \omega_{\neq}\|^2_{L^2}\\
& \le  \f{C}{A} \|\omega_{\neq}\|^2_{L^2}\big(\|\omega_{\neq}\|^2_{L^2} + \|\omega_0\|^2_{L^2} + \|u_0^1\|^2_{L^\infty}  \big) + \f{C}{A} \|\partial_xn_{\neq}\|_{L^2}\|\omega_{\neq}\|_{L^2} 
 +\f{C}{A}\|\nabla\omega_{\neq}\|^2_{L^2}\big(\|\omega_{\neq}\|^2_{L^2}+\|\omega_0\|^2_{L^2}\big).
\end{align*}
Then a direct integration with respect to $t$  yields that
\begin{align*}
&\|\omega_{\neq}(t)\|^2_{L^2}+\f{1}{A}\int_{0}^{t}\|\nabla \omega_{\neq}\|^2_{L^2}d\tau\\
& \le \|(\omega_{in})_{\neq}\|^2_{L^2}+\f{C}{A}\int_0^t \|\omega_{\neq}\|^2_{L^2}\big(\|\omega_{\neq}\|^2_{L^2} + \|\omega_0\|^2_{L^2} + \|u_0^1\|^2_{L^\infty}  \big) d\tau \\
& \qquad + \f{C}{A}\int_0^t \|\partial_xn_{\neq}\|_{L^2}\|\omega_{\neq}\|_{L^2}  d\tau
+ \f{C}{A}\int_0^t\|\nabla\omega_{\neq}\|^2_{L^2}\big(\|\omega_{\neq}\|^2_{L^2}+\|\omega_0\|^2_{L^2}\big)d\tau\\
& \le \|(\omega_{in})_{\neq}\|^2_{L^2} + \f{C}{A}\|\omega_{\neq}\|^2_{L^\infty L^2}\big(\|\omega_{\neq}\|^2_{L^\infty L^2} + \|\omega_0\|^2_{L^\infty L^2} + \|u_0^1\|^2_{L^\infty L^\infty}  \big)T  \\
& \qquad + \f{C}{A} \|\partial_xn_{\neq}\|_{L^\infty L^2}\|\omega_{\neq}\|_{L^\infty L^2} T
 + \f{C}{A}\big(\|\omega_{\neq}\|^2_{L^\infty L^2}+\|\omega_0\|^2_{L^\infty L^2}\big)\int_0^t\|\nabla\omega_{\neq}\|^2_{L^2}d\tau
\end{align*}
for all $t\in[0, T]$.  It follows from the assumptions (A-3), (A-5), (A-6),  { $0 < T \le \lambda_A^{-\f{1}{4}} = A^{\f{1}{8}}\mathrm{log}^{\f{1}{4}}A$ }, $\|\omega_{in}\|_X\le A^{-\f{3}{4}}$ and Lemma \ref{LH41} that 
\begin{align*}
& \|\omega_{\neq}(t)\|^2_{L^2}+\f{1}{A}\int_{0}^{t}\|\nabla \omega_{\neq}\|^2_{L^2}d\tau\nonumber\\
& \le \|(\omega_{in})_{\neq}\|^2_{L^2}+\f{CC_0^4}{A} \big(\|(\omega_{in})_{\neq}\|^2_{X} +{ A^{-\f{3}{2}}\big) }\big(\|\omega_{in}\|^2_{X} + \|{\bf u}_{in}\|^2_{L^2} + 1 \big) { A^{\f{1}{8}}\mathrm{log}^{\f{1}{4}}A  } \\
& \qquad  + \f{CC_0^2}{A} \|(\partial_xn_{in})_{\neq}\|_X \big(\|(\omega_{in})_{\neq}\|_X + { A^{-\f{3}{4}}\big)  A^{\f{1}{8}}\mathrm{log}^{\f{1}{4}}A } \\
& \qquad + \f{CC_0^2}{A} \big(\|\omega_{in}\|^2_X+{ A^{-\f{3}{2}} \big)} \big(\|(\omega_{in})_{\neq}\|^2_X +{ A^{-\f{3}{4}}\big) A }  \\
& \le \|(\omega_{in})_{\neq}\|^2_{L^2} + { \f{CC_0^4}{A^{\f58}} A^{-\f{3}{2}} }\big( \|{\bf u}_{in}\|^2_{L^2} + 1 \big)   + { \f{CC_0^2}{A^{\f58}} \|(\partial_xn_{in})_{\neq}\|_X A^{-\f{3}{4}} + CC_0^2A^{-\f{3}{2}}A^{-\f{3}{4}}  }  \\
& \le \|(\omega_{in})_{\neq}\|^2_{L^2}+{ \f{CC_0^4}{A^{\f58}} }\big(\|{\bf u}_{in}\|^2_{L^2} + \|(\partial_xn_{in})_{\neq}\|_{X}+ 1\big){ A^{-\f{3}{4}}. }
\end{align*}
Therefore, if we choose  $A^{*}_5$ fulfilling that 
\[
{ {A^{*}_5}^{\f58} }\ge CC_0^4 \big(\|{\bf u}_{in}\|^2_{L^2} + \|(\partial_xn_{in})_{\neq}\|_{X}+ 1\big),
\]
then it holds that 
\begin{equation}\label{454102}
\f{1}{A}\int_{0}^{t}\|\nabla \omega_{\neq}(\cdot,\tau)\|^2_{L^2}d\tau\leq \|(\omega_{in})_{\neq}\|^2_{L^2}+ { A^{-\f{3}{4}} }
\end{equation}
for all $t\in[0, T]$ and  $A>A^{*}_5$.

For the weighted part, we multiply equation ${\eqref{qh}}_3$ by $y^2\omega_{\neq}$ and then have
\begin{align}\label{4103}
&{\f{1}{2}}{\f{d}{dt}}\|y\omega_{\neq}\|^2_{L^2}+\f{1}{A}\|y\nabla \omega_{\neq}\|^2_{L^2}-\f{1}{A}\|\omega_{\neq}\|^2_{L^2} \nonumber\\
& = \f{1}{A} \int_{\mathbb T\times\mathbb R} \Big( ({\bf u}_{\neq} \omega_{\neq})_{\neq} + {\bf u}_{\neq} \omega_0 + {\bf u}_0 \omega_{\neq}  \Big) \cdot\nabla(y^2\omega_{\neq})dxdy    \nonumber\\
& \qquad + \f{1}{A}\int_{\mathbb T\times\mathbb R}\partial_xn_{\neq}y^2\omega_{\neq}dxdy 
    + 2 \int_{\mathbb T\times\mathbb R}u^2_{\neq}y^2\omega_{\neq}dxdy.
\end{align}
We now  estimate the integrals on the right hand side of \eqref{4103}. For the first one, we have
\begin{align*}
&\int_{\mathbb T\times\mathbb R} \Big( ({\bf u}_{\neq} \omega_{\neq})_{\neq} + {\bf u}_{\neq} \omega_0 + {\bf u}_0 \omega_{\neq}  \Big) \cdot\nabla(y^2\omega_{\neq})dxdy  \\
& = \int_{\mathbb T\times \mathbb R}\big(\big(({\bf u}_{\neq}\omega_{\neq})_{\neq} + {\bf u}_{\neq}\omega_0 \big)\cdot y^2\nabla \omega_{\neq} +  (u^1_0\omega_{\neq})y^2\partial_x\omega_{\neq}   \big) dxdy + 2 \int_{\mathbb T\times \mathbb R} (u^2_{\neq}\omega_{\neq})_{\neq} + u^2_{\neq}\omega_{0})y\omega_{\neq} dxdy\\
& \le \f{1}{2} \|y\nabla \omega_{\neq}\|^2_{L^2} + C\big( \|y{\bf u}_{\neq}\omega_{\neq}\|^2_{L^2} + \|y{\bf u}_{\neq}\omega_{0}\|^2_{L^2} + \|yu_0^1\omega_{\neq} \|^2_{L^2} \big)    + C\big(\|u^2_{\neq}\omega_{\neq}\|_{L^2} + \|u^2_{\neq}\omega_0\|_{L^2} \big) \|y\omega_{\neq}\|_{L^2}\\
  & \le \f{1}{2} \|y\nabla \omega_{\neq}\|^2_{L^2} + C\Big( \|{\bf u}_{\neq}\|^2_{L^\infty} \big(\|y\omega_{\neq}\|^2_{L^2} + \|y\omega_{0}\|^2_{L^2}\big) + \|u^1_0\|^2_{L^\infty} \|y\omega_{\neq} \|^2_{L^2} \Big) \\
& \qquad   + C\|u^2_{\neq}\|_{L^\infty} \big(\|\omega_{\neq}\|_{L^2} + \|\omega_0\|_{L^2} \big)\|y\omega_{\neq}\|_{L^2}\\
 & \le \f{1}{2} \|y\nabla \omega_{\neq}\|^2_{L^2} + C\Big( \|\omega_{\neq}\|_{L^2}\|\nabla\omega_{\neq}\|_{L^2} \big(\|\omega_{\neq}\|^2_X + \|\omega_0\|^2_X\big) + \|u^1_0\|^2_{L^\infty} \|\omega_{\neq} \|^2_X \Big) \\
& \qquad   + C  \|\omega_{\neq}\|_{L^2}^{\f12}\|\nabla\omega_{\neq}\|_{L^2}^{\f12}\big(\|\omega_{\neq}\|_{L^2} + \|\omega_0\|_{L^2} \big)\|\omega_{\neq}\|_X     \\
& \le \f{1}{2} \|y\nabla \omega_{\neq}\|^2_{L^2} + C\Big(\|\nabla\omega_{\neq}\|_{L^2} \big(\|\omega_{\neq}\|^3_X + \|\omega_0\|^3_X\big) +  \|\nabla\omega_{\neq}\|_{L^2}^{\f12}\big(\|\omega_{\neq}\|_{X}^{\f52} + \|\omega_0\|_{X}^{\f52}\big) + \|u^1_0\|^2_{L^\infty} \|\omega_{\neq} \|^2_X \Big), 
\end{align*}
while for the last two integrals, it follows from the H{\"o}lder inequality that
\[
\int_{\mathbb T\times\mathbb R}\partial_xn_{\neq}y^2\omega_{\neq}dxdy
\le \|y\partial_xn_{\neq}\|_{L^2}\|y\omega_{\neq}\|_{L^2}
\le \|\partial_xn_{\neq}\|_{X}\|\omega_{\neq}\|_{X}
\]
and that 
\[
2\int_{\mathbb T\times\mathbb R}u^2_{\neq}y^2\omega_{\neq}dxdy
\le 2 \|yu^2_{\neq}\|_{L^2}\|y\omega_{\neq}\|_{L^2}
\le C\|\omega_{\neq}\|^2_X 
\]
due to \eqref{cm341}.    Substituting the above estimates into \eqref{4103}, we obtain that
\begin{align*}
&{\f{d}{dt}}\|y\omega_{\neq}(t)\|^2_{L^2}+\f{1}{A}\|y\nabla \omega_{\neq}\|^2_{L^2}\\
& \le  \f{C}{A} \|u_0^1\|^2_{L^\infty} \|\omega_{\neq}\|^2_{X} +\f{C}{A}\|\partial_xn_{\neq}\|_X\|\omega_{\neq}\|_{X} + C \|\omega_{\neq}\|^2_{X} \\
& \qquad + \f{C}{A} \|\nabla\omega_{\neq}\|_{L^2} \big(\|\omega_{\neq}\|^3_X + \|\omega_0\|^3_X\big) + \f{C}{A}  \|\nabla\omega_{\neq}\|_{L^2}^{\f12}\big(\|\omega_{\neq}\|_{X}^{\f52} + \|\omega_0\|_{X}^{\f52}\big). 
\end{align*}
Integrating with respect to $t$, we deduce that 
\begin{align*}
&\|y\omega_{\neq}(t)\|^2_{L^2} + \f{1}{A}\int^t_0\|y\nabla \omega_{\neq}\|^2_{L^2}d\tau\\
& \le \|y(\omega_{in})_{\neq}\|^2_{L^2} + \f{C}{A}\int^t_0 \|u_0^1\|^2_{L^\infty} \|\omega_{\neq}\|^2_{X}  d\tau    + \f{C}{A}\int^t_0\|\partial_xn_{\neq}\|_X\|\omega_{\neq}\|_{X}d\tau+C\int^t_0\|\omega_{\neq}\|^2_{X}d\tau   \\
& \qquad + \f{C}{A}\int^t_0 \|\nabla\omega_{\neq}\|_{L^2} \big(\|\omega_{\neq}\|^3_X + \|\omega_0\|^3_X\big) d\tau
  + \f{C}{A}\int^t_0 \|\nabla\omega_{\neq}\|_{L^2}^{\f12}\big(\|\omega_{\neq}\|_{X}^{\f52} + \|\omega_0\|_{X}^{\f52}\big) d\tau \\
& \le \|y(\omega_{in})_{\neq}\|^2_{L^2}  
 + \f{C}{A} \|u_0^1\|^2_{L^\infty L^\infty} \|\omega_{\neq}\|^2_{L^\infty X}T   +\f{C}{A}\|\partial_xn_{\neq}\|_{L^\infty X}\|\omega_{\neq}\|_{L^\infty X}T+C\|\omega_{\neq}\|^2_{L^\infty X}T \\
& \qquad + \f{C}{A}\big(\|\omega_{\neq}\|^3_{L^\infty X} + \|\omega_0\|^3_{L^\infty X}\big)\Big(\int_0^t\|\nabla\omega_{\neq}\|^2_{L^2}d\tau\Big)^{\f{1}{2}} T^{\f{1}{2}} \\
& \qquad  + \f{C}{A} \big(\|\omega_{\neq}\|^{\f52}_{L^\infty X} + \|\omega_0\|^{\f52}_{L^\infty X}\big) \Big(\int_0^t\|\nabla\omega_{\neq}\|^2_{L^2}d\tau\Big)^{\f{1}{4}} T^{\f{3}{4}}
\end{align*}
for all $t\in[0, T]$.  It then follows from the assumptions (A-3), (A-5), (A-6), { $0 < T \le \lambda_A^{-\f{1}{4}} = A^{\f{1}{8}}\mathrm{log}^{\f{1}{4}}A$ }, $\|\omega_{in}\|_X \le A^{-\f{3}{4}}$ and Lemma \ref{LH41}  that
\begin{align*}
&\|y\omega_{\neq}(t)\|^2_{L^2}+\f{1}{A}\int^t_0\|y\nabla \omega_{\neq}\|^2_{L^2}d\tau\nonumber\\
& \le \|y(\omega_{in})_{\neq}\|^2_{L^2}
  +\f{CC_0^2}{A}\big(\|{\bf u}_{in}\|^2_{L^2}+\|(\omega_{in})_{0}\|^2_{X} + 1 \big)\big(\|(\omega_{in})_{\neq}\|^2_{X}+ { A^{-\f{3}{2}}\big) A^{\f{1}{8}}\mathrm {log}^{\f{1}{4}}A }  \\
& \qquad + \f{CC_0}{A} \|(\partial_xn_{in})_{\neq}\|_{X}\big(\|(\omega_{in})_{\neq}\|_{X}+ { A^{-\f{3}{4}}\big) A^{\f{1}{8}}\mathrm {log}^{\f{1}{4}}A } + 
CC_0^2\big(\|(\omega_{in})_{\neq}\|^2_{X}+{ A^{-\f{3}{2}}\big) A^{\f{1}{8}}\mathrm {log}^{\f{1}{4}}A } \\
&\qquad + \f{CC_0^3}{A}\big(\| \omega_{in}\|^3_{X} + { A^{-\f{9}{4}} \big) }
\big(\|(\omega_{in})_{\neq}\|_{X}+{ A^{-\f{3}{8}}\big)A^{\f12}  A^{\f{1}{16}}\mathrm {log}^{\f{1}{8}}A   }\\
&\qquad + \f{CC_0^{\f53}}{A}\big(\| \omega_{in}\|^{\f52}_{X} + { A^{-\f{15}{8}}\big) }
\big(\|(\omega_{in})_{\neq}\|_{X}+{ A^{-\f{3}{8}}\big)^{\f12} A^{\f14}  A^{\f{3}{32}}\mathrm {log}^{\f{3}{16}}A  } \\ 
& \le \|y(\omega_{in})_{\neq}\|^2_{L^2}
  +{ \f{CC_0^2}{A^{\f58}} }\big(\|{\bf u}_{in}\|^2_{L^2}+\|(\omega_{in})_{0}\|^2_{X} + 1 \big) { A^{-\f{3}{2}}  + \f{CC_0}{A^{\f58}}\|(\partial_xn_{in})_{\neq}\|_{X} A^{-\f{3}{4}}    }  \\
& \qquad  +  CC_0^2 { A^{-\f{3}{2}}  A^{\f{1}{8}}\mathrm {log}^{\f{1}{4}}A + \f{CC_0^3}{A^{\f{5}{16}}} A^{-\f{9}{4}} A^{-\f{3}{8}}   + \f{CC_0^{\f53}}{A^{\f{15}{32}}} A^{-\f{15}{8}}  A^{-\f{3}{16}}  } \\ 
& \le \|y(\omega_{in})_{\neq}\|^2_{L^2}
  +{  \f{CC_0^3}{A^{\f38}} }\big(\|{\bf u}_{in}\|^2_{L^2}+\|(\omega_{in})_0\|^2_X + \|(\partial_xn_{in})_{\neq}\|_X + 1 \big) { A^{-\f{3}{4}}. }
\end{align*}
Therefore, by choosing  $A_5^\circ$ such that 
\[
{ {A_5^\circ}^{\f{3}{8}}  } \ge CC_0^3 \big(\|{\bf u}_{in}\|^2_{L^2}+\|(\omega_{in})_0\|^2_X + \|(\partial_xn_{in})_{\neq}\|_X + 1 \big),
\]
we obtain 
\begin{equation}\label{454112}
\f{1}{A}\int_{0}^{t}\|y\nabla \omega_{\neq}(\cdot,\tau)\|^2_{L^2}d\tau\leq \|y(\omega_{in})_{\neq}\|^2_{L^2}+{ A^{-\f{3}{4}}  }
\end{equation}
for all $t\in[0, T]$ and  $A>A_5^\circ$.  It then follows from \eqref{454102} and \eqref{454112} that
\[
\f{1}{A}\int_{0}^{t}\|\nabla \omega_{\neq}(\cdot,\tau)\|^2_{X}d\tau\leq 2\Big(\|(\omega_{in})_{\neq}\|^2_{X}+ { A^{-\f{3}{4}}\Big) }
\]
for all $A>A_5$ with $A_5:=\max\big\{A_5^{*}, A_5^\circ\}$.  This completes the proof of Lemma \ref{LH45}. \qquad $\Box$

\begin{Lemma}\label{LH46}
Under the assumptions {\rm  (A-3)}, {\rm (A-5)} and {\rm (A-6)},  there exists a positive constant $A_6$ relying on $C_0$, $\|(\partial_xn_{in})_{\neq}\|_X$ and $\|{\bf u}_{in}\|_{L^2}$ such that if $A>A_6$ and $A^{\f{3}{4}}\|\omega_{in}\|_X\leq1$, then it holds 
\[
\|\omega_{\neq}(t)\|_X\leq 2C_0e^{-\epsilon_0\lambda_At}\big(\|(\omega_{in})_{\neq}\|_X+{ A^{-\frac{3}{4}}\big) } \qquad {\rm  for \,\, all \,\,\,}  t\in[0, T]. 
\]
\end{Lemma}
{\bf Proof.} Recalling ${\widetilde{\mathcal L}} = \f{1}{A}\Delta - y^2\partial_x + 2\partial_x\Delta^{-1}$, we can use  the Duhamel principle to assert that
\[
\omega_{\neq}=e^{\widetilde{\mathcal L}t}(\omega_{in})_{\neq}-\f{1}{A}e^{\widetilde{\mathcal L}t}\int_0^te^{-\widetilde{\mathcal L}\tau}\big(\nabla\cdot(\omega_{\neq}{\bf u_{\neq}})_{\neq}+\nabla\cdot(\omega_{\neq}{\bf u_0})+\nabla\cdot(\omega_0 {\bf u_{\neq}})-\partial_xn_{\neq}\big)d\tau.
\]
It then follows from Lemma \ref{LH31} and the fact  the fact $0\leq e^{\epsilon_0\lambda_A \tau}\leq C$ that 
\begin{align}\label{464115}
\|\omega_{\neq}(t)\|_{X}
& \le C_0e^{-\epsilon_0\lambda_A t}\Big(\sqrt{2}\|(\omega_{in})_{\neq}\|_{X}  \nonumber\\
& \qquad +\f{\sqrt{2}}{A} \int_0^te^{\epsilon_0\lambda_A \tau}\big\|\nabla\cdot(\omega_{\neq}{\bf u_{\neq}})_{\neq} + \nabla\cdot(\omega_{\neq}{\bf u_0}) + \nabla\cdot(\omega_0 {\bf u_{\neq}}) - \partial_xn_{\neq}\big\|_{X}d\tau\Big) \nonumber \\
& \le C_0e^{-\epsilon_0\lambda_A t}\Big(\sqrt{2}\|(\omega_{in})_{\neq}\|_{X}  \nonumber\\
& \qquad +\f{C}{A} \int_0^t \big\|\nabla\cdot(\omega_{\neq}{\bf u_{\neq}})_{\neq} + \nabla\cdot(\omega_{\neq}{\bf u_0}) + \nabla\cdot(\omega_0 {\bf u_{\neq}}) - \partial_xn_{\neq}\big\|_{X}d\tau\Big). 
\end{align}
For the integral on the right hand side of \eqref{464115},  we deduce from the H\"{o}lder inequality and \eqref{43453} that
\begin{align*}
&\int_0^t \big\|\nabla\cdot(\omega_{\neq}{\bf u_{\neq}})_{\neq} + \nabla\cdot(\omega_{\neq}{\bf u_0}) + \nabla\cdot(\omega_0 {\bf u_{\neq}}) - \partial_xn_{\neq}\big\|_{X}d\tau  \\
& \le C\int_0^t \Big( \|\nabla \omega_{\neq}\cdot{\bf u_{\neq}}\|_{X} + \|\partial_x\omega_{\neq}u_0^1\|_{X} +  \|\partial_y\omega_0u^2_{\neq}\|_{X}  +  \|\partial_xn_{\neq}\|_{X}  \Big) d\tau\\
& \le  C\int_0^t \Big(\|\nabla \omega_{\neq}\|_{X}\|{\bf u}_{\neq}\|_{L^\infty}   + \|\partial_x \omega_{\neq}\|_{X}\|u^1_0\|_{L^\infty}   +   \|\partial_y \omega_{0}\|_{X}\|u^2_{\neq}\|_{L^\infty}    +  \|\partial_xn_{\neq}\|_{X}  \Big) d\tau\\
& \le  C\int_0^t \Big( \|\nabla\omega_{\neq}\|_{X}^{\f{3}{2}}\|\omega_{\neq}\|_{L^2}^{\f{1}{2}}    +  \|\partial_x \omega_{\neq}\|_{X}\|u^1_0\|_{L^\infty}   +   \|\partial_y \omega_0\|_{X}\|\nabla\omega_{\neq}\|_{X}^{\f{1}{2}}\|\omega_{\neq}\|_{L^2}^{\f{1}{2}}     +  \|\partial_xn_{\neq}\|_{X}  \Big) d\tau    \\ 
& \le C \|\omega_{\neq}\|_{L^\infty L^2}^{\f12} \Big(\int_0^t\|\nabla\omega_{\neq}\|_{X}^2d\tau\Big)^{\f{3}{4}} T^{\f{1}{4}}  + C \|u_0^1\|_{L^\infty L^\infty} \Big(\int_0^t\|\partial_x \omega_{\neq}\|^2_{X}d\tau\Big)^{\f{1}{2}} T^{\f{1}{2}} \\
& \qquad   + C \|\omega_{\neq}\|_{L^\infty L^2}^{\f12}\Big(\int_0^t\|\partial_y \omega_{0}\|^2_{X}d\tau\Big)^{\f{1}{2}}\Big(\int_0^t\|\nabla \omega_{\neq}\|^2_{L^2}d\tau\Big)^{\f{1}{4}} T^{\f14} 
+ C\|\partial_xn_{\neq}\|_{L^\infty X} T. 
\end{align*}
Substituting the last inequality  into \eqref{464115}, we obtain that
\begin{align*}
\|\omega_{\neq}(t)\|_{X}
& \le  C_0e^{-\epsilon_0\lambda_A t}\Bigg( \sqrt{2} \big\|( \omega_{in})_{\neq}\|_{X}
+ \f{C}{A} \|\omega_{\neq}\|_{L^\infty L^2}^{\f12} \Big(\int_0^t\|\nabla\omega_{\neq}\|_{X}^2d\tau\Big)^{\f{3}{4}} T^{\f{1}{4}}  \\
&\qquad  + \f{C}{A} \|u_0^1\|_{L^\infty L^\infty} \Big(\int_0^t\|\nabla \omega_{\neq}\|^2_{X}d\tau\Big)^{\f{1}{2}} T^{\f{1}{2}}  \\
& \qquad  + \f{C}{A} \|\omega_{\neq}\|_{L^\infty L^2}^{\f12}\|\partial_y \omega_{0}\|_{L^2X}\Big(\int_0^t\|\nabla \omega_{\neq}\|^2_{L^2}d\tau\Big)^{\f{1}{4}} T^{\f14}      + \f{C}{A} \|\partial_xn_{\neq}\|_{L^\infty X} T\Bigg).
\end{align*}
Then by the assumptions (A-3), (A-5), (A-6),  { $0 < T \le \lambda_A^{-\f{1}{4}} = A^{\f{1}{8}}\mathrm{log}^{\f{1}{4}}A$ }, $\|\omega_{in}\|_X \le A^{-\f{3}{4}}$  and Lemma \ref{LH41},  we can deduce that
\begin{align*}
\|\omega_{\neq}(t)\|_{X}
& \le  C_0e^{-\epsilon_0\lambda_A t}\Bigg( \sqrt{2} \big\|( \omega_{in})_{\neq}\|_{X}
+ \f{CC_0^{\f12}}{A} \big(\|( \omega_{in})_{\neq}\|_{X}+{ A^{-\f{3}{4}} }\big)^{\f12}\big(\|( \omega_{in})_{\neq}\|_{X}^2 + { A^{-\f{3}{4}} 
 }\big)^{\f34} A^{\f34} { A^{\f{1}{32}}\mathrm{log}^{\f{1}{16}}A }\\
&\qquad  + \f{C}{A} \big(\|{\bf u}_{in}\|_{L^2}+\|( \omega_{in})_{0}\|_X+1\big) \big(\|( \omega_{in})_{\neq}\|_{X}^2 + { A^{-\f{3}{4}}\big)^{\f12} A^{\f12}  A^{\f{1}{16}}\mathrm{log}^{\f{1}{8}}A  }\\
& \qquad  + \f{CC_0^{\f12}}{A} \big(\|( \omega_{in})_{\neq}\|_{X}+{ A^{-\f{3}{4}}  }\big)^{\f12} \big(\|( \omega_{in})_0\|_{X} + A^{-\f{3}{4}}\big) A^{\f12}  \big(\|( \omega_{in})_{\neq}\|_{X}^2 +{  A^{-\f{3}{4}}\big)^{\f14} A^{\f14} A^{\f{1}{32}}\mathrm{log}^{\f{1}{16}}A  } \\
& \qquad    +  \f{C}{A} \|(\partial_xn_{in})_{\neq}\|_{X} { A^{\f{1}{8}}\mathrm{log}^{\f{1}{4}}A }\Bigg) \\
& \le  C_0e^{-\epsilon_0\lambda_A t}\Big( \sqrt{2} \big\|( \omega_{in})_{\neq}\|_{X}
+ \f{CC_0^{\f12}}{A} { A^{-\f{3}{8}} A^{-\f{9}{16}}  A^{\f34} A^{\f{1}{32}}\mathrm{log}^{\f{1}{16}}A }
 + \f{C}{A} \big(\|{\bf u}_{in}\|_{L^2} + 1\big) { A^{-\f{3}{8}} A^{\f12}  A^{\f{1}{16}}\mathrm{log}^{\f{1}{8}}A } \\
& \qquad  + \f{CC_0^{\f12}}{A} { A^{-\f{3}{8}} A^{-\f{3}{4}} A^{\f12}  A^{-\f{3}{16}} A^{\f14} A^{\f{1}{32}}\mathrm{log}^{\f{1}{16}}A}   +  \f{C}{A} \|(\partial_xn_{in})_{\neq}\|_{X} { A^{\f{1}{8}}\mathrm{log}^{\f{1}{4}}A \Big)} \\
& \le  C_0e^{-\epsilon_0\lambda_A t}\Big( \sqrt{2} \big\|( \omega_{in})_{\neq}\|_{X}
+ { \f{CC_0^{\f12}}{A^{\f38}} A^{-\f{3}{4}}
  + \f{C}{A^{\f{1}{32}}} \big(\|{\bf u}_{in}\|_{L^2} + 1\big) A^{-\f{3}{4}}  }\\
& \qquad  + { \f{CC_0^{\f12}}{A^{\f{23}{32}}} A^{-\f{3}{4}}  +  \f{C}{A^{\f{1}{16}}} \|(\partial_xn_{in})_{\neq}\|_{X} A^{-\f{3}{4}} \Big)} \\
& \le  C_0e^{-\epsilon_0\lambda_A t}\Big(\sqrt{2} \big\|( \omega_{in})_{\neq}\|_{X} + { \f{CC_0}{A^{\f{1}{32}}} }\big(\|(\partial_xn_{in})_{\neq}\|_{X}+\|{\bf u}_{in}\|_{L^2}+1\big) { A^{-\f{3}{4}}\Big).}
\end{align*}
Thus we only need to fix  $A_6$  such that 
\[
{ {A_6}^{\f{1}{32}}\ge }CC_0\big(\|(\partial_xn_{in})_{\neq}\|_{X}+\|{\bf u}_{in}\|_{L^2}+1\big),
\]
which clearly entails that 
\[
\|\omega_{\neq}(t)\|_{X}\leq 2C_0e^{-\epsilon_0\lambda_A t}\big(\|(\omega_{in})_{\neq}\|_{X}+{ A^{-\f{3}{4}}\big)  }
\]
for $t\in[0, T]$ and $A>A_6$.  This completes the proof of Lemma \ref{LH46}. \qquad $\Box$

\begin{Lemma}\label{LH47}
Under the assumptions Proposition \ref{LH41},  there exists a positive constant $C_\infty\geq1$ relying on $C_0$, $M:=\|n_{in}\|_{L^1}$ and $\|n_{in}\|_{X\cap L^\infty}$ such that
\[
\|n\|_{L^{\infty}L^{\infty}}\leq 2C_\infty.
\]
\end{Lemma}
{\bf Proof.}  We can establish the desired  $L^\infty L^\infty$ estimate of $n$ by a Moser-Alikakos iteration, which is almost parallel  to Lemma 4.5 in Zhang et al. \cite{ZZZ}. We omit the details here,  but we would like to remark that in the current setting, the constant $C_\infty$ can be independent of the solution component  $c$ due to our universal constant in the estimate of $\big\|\na c\big\|_{L^4}$ (see Lemma \ref{LH35} and Lemma \ref{LH36}).  \qquad $\Box$
\\
\\
{\bf Proof of Proposition \ref{LH1}.} Collecting Lemma \ref{LH42} - Lemma \ref{LH47} and choosing 
\[
A_0=\max\big\{A_2, A_3, A_4, A_5, A_6\big\},
\]
it is easy to end the proof of Proposition \ref{LH1}.         \qquad $\Box$

\section{Proof of Theorem \ref{global regular}}\label{proof}

\quad In this section, we will prove our main theorem by combining the bootstrap estimates and the standard continuity argument  with the help of the following  local well-posedness result, which can be proven through standard argument.
\begin{Theorem}\label{local}
Assume that the initial data satisfy $n_{in}\in L^1\cap X \cap L^\infty({\mathbb T}\times{\mathbb R})$ with  $\partial_xn_{in}\in X$ and  ${\bf u}_{in}\in H^1({\mathbb T}\times{\mathbb R})$ with $\omega_{in}\in X.$ Then for each $A>0$,  there exist $T^* \in(0,+\infty]$ and a unique triple $(n, c, \mathbf{u})$ with $n\ge 0$ and $c\ge 0$ solving system \eqref{3} in ${\mathbb T}\times{\mathbb R}\times(0, T^*)$. 
\end{Theorem}

For simplicity, we first introduce the following weighted energy functional:
\[
E(T) = \f{1}{A}\|\nabla n_{\neq}\|_{L^2X}^2 + \|e^{\epsilon_0\lambda_A t} n_{\neq}\|_{L^\infty X} + \|\partial_x n_{\neq}\|_{L^\infty X}^2 + \f{1}{A}\|\nabla \omega_{\neq}\|_{L^2X}^2 + \|e^{\epsilon_0\lambda_A t} \omega_{\neq}\|_{L^\infty X}. 
\]
Then the assumptions  (A-1) - (A-6) are equivalent to the hypotheses:
\begin{equation}\label{etbound}
E(T)\leq 4C_{\neq}       \qquad \mathrm{and} \qquad  
\|n\|_{L^\infty L^\infty}\le 4C_{\infty},
\end{equation}
where $C_{\neq}$ is defined by 
\[
C_{\neq}:=\Big(\|(n_{in})_{\neq}\|_X^2 + C_0\|(n_{in})_{\neq}\|_X + \|(\partial_xn_{in})_{\neq}\|_X^2 + \|(\omega_{in})_{\neq}\|_X^2 + A^{-\f{7}{16}} + C_0\big(\|(n_{in})_{\neq}\|_X + A^{-\f{7}{16}}\big)\Big),
\]
while the refined bounds  (B-1) - (B-6) are equivalent to 
\begin{equation}\label{reetbound}
E(T)\leq 2C_{\neq}    \qquad \mathrm{and} \qquad  
\|n\|_{L^\infty L^\infty}\le 2C_{\infty}. 
\end{equation}

By this simplification,  Proposition \ref{LH1} can be restated as follows. 
\begin{Proposition}\label{LH11}
Assume that the initial data $(n_{in}, {\bf u}_{in})$ satisfy the assumptions of Theorem \ref{global regular}. If the solution $(n, c, {\bf u})$ of system \eqref{3} possesses the bounds \eqref{etbound}, then there exists a positive constant  $A_0$ depending only on $C_0$, $C_\infty$,  $\|n_{in}\|_{L^1\cap X\cap L^\infty}$,  $\|(\partial_xn_{in})_{\neq}\|_X$  and  $\|{\bf u}_{in}\|_{L^2}$ such that  the refined bounds \eqref{reetbound} hold  provided that  $A>A_0$ and $A^{\f{3}{4}}\|\omega_{in}\|_X\leq 1$. 
\end{Proposition} 

Now we are in the position to end the proof of our main result. 
\\
\\
{\bf Proof of Theorem \ref{global regular}.} Under the regularity assumptions on the initial data of  Theorem \ref{global regular}, it is easy to show the local well-posedness of system \eqref{3}(Theorem \ref{local}), which together with Lemma \ref{LH41} and Proposition \ref{LH11} will yield the global local well-posedness  by a continuity argument.     \qquad $\Box$
\\
\\
{\bf Acknowledgements.} This project is partially supported by National Key R\&D Program of China, Project Number 2021YFA1001200. Z. Xiang was supported by the NNSF of China (no. 11971093) and the Special Funds for Local Scientific and Technological Development Guided by the Central Government (no. 2022ZYD0007). The work of X.Xu is supported by NNSF of China Youth program (no. 12101278), and Kunshan Shuangchuang Talent Program (no. kssc202102066).


{\small
\begin{thebibliography}{50}


\bibitem{A} N. D. Alikakos,
\it $L^p$ bounds of solutions of reaction-diffusion equations,
\rm Commun. Partial Differ. Equ. 4 (1979), 827-868.

\bibitem{BH} J. Bedrossian, S. He,
\it Suppression of blow-up in Patlak-Keller-Segel via shear flows,
\rm SIAM J. Math. Anal. 49 (2017) 4722-4766.

\bibitem{BC} A. Blanchet, J. A. Carrillo, N. Masmoudi,
\it Infinite time aggregation for the critical Patlak-Keller-Segel model in $\R^2$,
\rm Commun. Pure Appl. Math. 61 (2008) 1449-1481.

\bibitem{BD} A. Blanchet, J. Dolbeault, B. Perthame,
\it Two-dimensional Keller-Segel model: optimal critical mass and qualitative properties of the solutions,
\rm Electron. J. Differ. Equ. 2006 (2006) 1-32.

\bibitem{CKL} M. Chae, K. Kang, J. Lee,
\it Global existence and temporal decay in Keller-Segel models coupled to fluid equations,
\rm Commun. Partial Differ. Equ. 39 (2014) 1205-1235.


\bibitem{Coll} J. C. Coll, B. E Bowden, G. V. Meehan et al.,  
\it Chemical aspects of mass spawning in corals. I. Sperm-attractant molecules in eggs of the scleractinian coral Montipora digitata, 
\rm Mar. Biol.  118 (1994) 177-182.


\bibitem{CPZ} L. Corrias, B. Perthame, H. Zaag,
\it Global solutions of some chemotaxis and angiogenesis systems in high space dimensions,
\rm Milan J. Math. 72 (2004) 1-28.


\bibitem{ZEW} M. Coti Zelati, T. M. Elgindi, K. Widmayer,
\it Enhanced dissipation in the Navier-Stokes equations near the Poiseuille flow,
\rm Commun. Math. Phys. 378 (2020) 987-1010.

\bibitem{DL} R. Duan, A. Lorz, P. Markowich,
\it Global solutions to the coupled chemotaxis-fluid equations,
\rm Commun. Partial Differ. Equ. 35 (2010) 1635-1673.

\bibitem{ES} E. Espejo, T. Suzuki, 
\it Reaction terms avoiding aggregation in slow fluids, 
\rm Nonlinear Anal. Real World Appl. 21 (2015) 110-126. 

\bibitem{GhH} S. Ghorai,  N. Hill,
\it Development and stability of gyrotactic plumes in bioconvection, 
\rm J. Fluid Mech. 400 (1999) 1-31. 


\bibitem{GH} Y. Gong,  S. He,
\it  On the $8\pi$-subcritical mass threshold of a Patlak-Keller-Segel-Navier-Stokes system,
\rm SIAM J. Math. Anal.  53 (2021) 2925-2956. 

\bibitem{H} S. He,
\it Suppression of blow-up in parabolic-parabolic Patlak-Keller-Segel via strictly monotone shear flows,
\rm Nonlinearity 31 (2018) 3651-3688.


\bibitem{H2} S. He,
\it Enhanced dissipation and blow-up suppression in a chemotaxis-fluid system, 
\rm 	arXiv: 2207.13494. 

\bibitem{HT} S. He, E. Tadmor,
\it Suppressing chemotactic blow-up through a fast splitting scenario on the plane,
\rm Arch. Ration. Mech. Anal. 232 (2019) 951-986.

\bibitem{HTZ} S. He, E. Tadmor, A. Zlatos,
\it On the fast spreading scenario,
\rm Comm. Amer. Math. Soc. 2 (2022) 149-171. 

\bibitem{IXZ} G. Iyer, X. Xu, and A. Zlatos,
\it Convection-induced singularity suppression in the Keller-Segel and other non-linear PDEs, 
\rm Trans. Amer. Math. Soc. 374 (2021) 6039-6058.


\bibitem{WL} W. J{\"a}ger, S. Luckhaus,
\it On explosions of solutions to a system of partial differential equations modelling chemotaxis,
\rm Trans. Am. Math. Soc. 329 (1992) 819-824.

\bibitem{KS} E. F. Keller, L. A. Segel,
\it Initiation of slime mold aggregation viewed as an instability,
\rm J. Theor. Biol. 26 (1970) 399-415.


\bibitem{Ke} J. O. Kessler, 
\it Hydrodynamic focussing of motile algal cells,
\rm  Nature 313 (1985) 218-220. 


\bibitem{KR} A. Kiselev, L. Ryzhik, 
\it Biomixing by chemotaxis and enhancement of biological reactions, 
\rm Comm. Partial Differential Equations 37 (2012) 298-312.


\bibitem{KX} A. Kiselev, X. Xu,
\it Suppression of chemotactic explosion by mixing,
\rm Arch. Ration. Mech. Anal. 222 (2016) 1077-1112.


\bibitem{KMS2} H. Kozono, M. Miura, Y. Sugiyama,
\it Time global existence and finite time blow-up criterion for solutions to the Keller-Segel system coupled with the Navier-Stokes fluid,
\rm J. Differ. Equ. 267 (2019) 5410-5492.


\bibitem{LWZ}  C. Lai, J. Wei, Y. Zhou, 
\it Global existence of free-energy solutions to the 2D Patlak-Keller-Segel-Navier-Stokes system with critical and subcritical mass, 
\rm Indiana Univ. Math. J. to appear. 



\bibitem{LWZ} T. Li, D. Wei, Z. Zhang,
\it Pseudospectral and spectral bounds for the Oseen vortices operator,
\rm arXiv preprint arXiv:1701.06269 (2017).

\bibitem{LL} J. Liu, A. Lorz, 
\it A coupled chemotaxis-fluid model: Global existence, 
\rm Ann. Inst. H. Poincar\'e Anal. Non Lin\'eaire 28 (2011) 643-652.

\bibitem{LO} A. Lorz,
\it Coupled chemotaxis fluid model,
\rm Math. Models Methods Appl. Sci. 20 (2010) 987-1004.

\bibitem{LA} A. Lorz,
\it A coupled Keller-Segel-Stokes model: global existence for small initial data and blow-up delay
\rm Commun. Math. Sci. 10 (2012) 555-574.

\bibitem{N} T. Nagai,
\it Blow-up of radially symmetric solutions to a chemotaxis system,
\rm Adv. Math. Sci. Appl. 5 (1995) 581-601.

\bibitem{OK} S. A. Orszag, L. C. Kells, 
\it Transition to turbulence in plane Poiseuille and plane Couette flow, 
\rm Journal of Fluid Mechanics  96 (1980)  159-205.

\bibitem{P} C. S. Patlak,
\it Random walk with persistence and external bias,
\rm Bull. Math. Biophys. 15 (1953) 311-338.


\bibitem{T} I. Tuval, L. Cisneros, C. Dombrowski, C. W. Wolgemuth, J. O. Kessler, R. E. Goldstein,
\it Bacterial swimming and oxygen transport near contact lines,
\rm Proc. Nati. Acad. Sci. USA. 102 (2005) 2277-2282.


\bibitem{VMG} F. J. Vermolen, M. M. Mul, A. Gefen, 
\it Semi-stochastic cell-level computational modeling of the immune system response to bacterial infections and the effects of antibiotics, 
\rm Biomech Model Mechanobiol 13 (2014) 713-734. 



\bibitem{WWX1} Y. Wang, M. Winkler, Z. Xiang, 
\it Global classical solutions in a two-dimensional chemotaxis-Navier-Stokes system with subcritical sensitivity, 
\rm Ann. Sc. Norm. Super. Pisa Cl. Sci. XVIII (2018) 421-466.


\bibitem{WWX2} Y. Wang, M. Winkler, Z. Xiang, 
\it Local energy estimates and global solvability in a three-dimensional chemotaxis-fluid system with prescribed signal on the boundary, 
\rm Commun. Partial Differ. Equ. 46 (2021) 1058-1091.


\bibitem{W1} M. Winkler,
\it Global large-data solutions in a chemotaxis-(Navier-)Stokes system modeling cellular swimming in fluid drops,
\rm Commun. Partial Differ. Equ. 37 (2012) 319-351.

\bibitem{Win2013} M. Winkler,
\it Finite-time blow-up in the higher-dimensional parabolic-parabolic Keller-Segel system, 
\rm J. Math. Pures Appl. 100 (2013) 748-767.

\bibitem{W2} M. Winkler,
\it Stabilization in a two-dimensional chemotaxis-Navier-Stokes system,
\rm Arch. Ration. Mech. Anal. 211 (2014) 455-487.

\bibitem{W3} M. Winkler,
\it Global weak solutions in a three-dimensional chemotaxis-Navier-Stokes system,
\rm Ann. Inst. H. Poincar\'e Anal. Non Lin\'eaire 33 (2016) 1329-1352.


\bibitem{W4} M. Winkler,
\it Small-mass solutions in the two-dimensional Keller-Segel system coupled to the Navier-Stokes equations, 
\rm SIAM J. Math. Anal. 52 (2020) 2041-2080. 


\bibitem{ZZZ} L. Zeng, Z. Zhang, R. Zi,
\it Suppression of blow-up in Patlak-Keller-Segel-Navier-Stokes system via the Couette flow,
\rm J. Funct. Anal. 280 (2021) 108967.

\bibitem{ZZ} Q. Zhang, X. Zheng,
\it Global well-posedness for the two-dimensional incompressible chemotaxis-Navier-Stokes equations,
\rm SIAM J. Math. Anal. 46 (2014) 3078-3105.
\end{thebibliography}}
\end{document}